\documentclass[nonblindrev]{informs3} 

\OneAndAHalfSpacedXI 

\usepackage{endnotes}
\let\footnote=\endnote
\let\enotesize=\normalsize
\def\notesname{Endnotes}%
\def\makeenmark{$^{\theenmark}$}
\def\enoteformat{\rightskip0pt\leftskip0pt\parindent=1.75em
  \leavevmode\llap{\theenmark.\enskip}}

\allowdisplaybreaks[4]

\usepackage{natbib}
 \bibpunct[, ]{(}{)}{,}{a}{}{,}%
 \def\bibfont{\small}%

\usepackage{natbib}
\usepackage{subfigure}
\usepackage{epstopdf}
\usepackage{booktabs}
\usepackage{color}
\usepackage{multirow}
\usepackage{dcolumn}
\usepackage{ulem}
\usepackage{cancel}

\usepackage[ruled,linesnumbered]{algorithm2e}
\TheoremsNumberedThrough     
\ECRepeatTheorems

\EquationsNumberedThrough    

\begin{document}


\RUNAUTHOR{Lim, Wei, Zhang}

\RUNTITLE{Partial Backorder Inventory Control with Demand Learning}

\TITLE{Partial Backorder Inventory System: Asymptotic Optimality and Demand Learning}

\ARTICLEAUTHORS{%
\AUTHOR{Andrew E. B. Lim}
\AFF{Institute of Operations Research and Analytics, National University Singapore, Singapore, 117602, \EMAIL{andewlim@nus.edu.sg}}
\AUTHOR{Zhao-Xuan Wei}
\AFF{Institute of Operations Research and Analytics, National University Singapore, Singapore, 117602, \EMAIL{zhaoxuanwei@u.nus.edu}}
\AUTHOR{Hanqin Zhang}
\AFF{Department of Analytics and Operations, NUS Business School , National University of Singapore, 117602, \EMAIL{bizzhq@nus.edu.sg}}
} 

\ABSTRACT{%
We develop a stochastic inventory system which accounts for the limited patience of backlogged customers. While limited patience is a feature that is closer to the nature of unmet demand, our model also unifies the classic backlogging and lost-sales inventory systems which are special cases of the one we propose. We establish the uniform (asymptotic) optimality of the base-stock policy when both demand and patience distributions are known. When the backlogged demands become unobservable, we introduce a novel policy family that operates without backlogged demands information, and prove that it can approach the cost efficiency of the optimal policy in the system when the demand and patience distributions are known. Finally, we consider an online inventory control problem in which backlogged demand is unobservable and demand and patience distributions are also not known, and develop a UCB-type algorithm that yields a near-optimal policy. The regret bounds given by the algorithm are provably tight within the planning horizon, and are comparable to the state-of-the-art results in the literature, even in the face of partial and biased observations and weaker system ergodicity.
}%


 \HISTORY{}

\maketitle

%


\section{\sf Introduction}\label{intro}
In most stochastic inventory management systems, the unmet demands are assumed to be either completely lost or backordered. 
However, the reality is somewhere between these two extremes: backlogged customers may be willing to wait but  can leave without notice when their patience expires  (see \cite{corstengruen2004} for more discussion on this matter and how  patience varies across different industries). Moreover, the pool of waiting customers is often not observable and estimating its size, which is critical for the ordering decision, is difficult.


 We consider the problem of inventory control where customers with unmet demand are impatient. In particular, we consider a periodic-review random-demand inventory problem in which, at each period, unsatisfied customers decide whether to stay in the system (backorder) or abandon (lost-sale). Specifically, we assume that  every backlogged customer flips a coin at each epoch, stays if she flips heads (with probability $p\in[0, 1]$) but otherwise leaves the system. A larger value of $p$ corresponds to a more patient customer. We assume that the decision to stay or leave is independent across every backlogged customer. Our model is a hybrid of the classical  backorder and lost-sales models, which correspond to setting $p=1$ and $p=0$, respectively. In the following, we call it the {\it partial backorder inventory system}.

For an intuitive overview of the partial backorder system, consider the case where lead time is positive and inventory is being replenished according to a base-stock policy with level $s$.
For the backorder  system, if the initial inventory position (on-hand inventory and the inventory on-order net the backorder) is below $s$, then the inventory position at the beginning of each period will always stay below level $s$. This also holds for the lost-sales  system. However, for the partial backorder system, this no longer holds when the lead time is positive because backordered customers can lose patience before the order that was placed for them arrives, causing the inventory position to jump  above the level $s$. This  {\it overshooting} feature makes the system dynamics more complicated and difficult to be analyzed.

The following three questions naturally arise for the partial backorder inventory system:
\begin{itemize}
\item[(i)] 
What is an asymptotically optimal policy when all parameters are known and the backorders  are observable?
\item[(ii)] 
What is asymptotically optimal when  backordered demand is not observable but model parameters are known?
\item[(iii)] 
Is there an online algorithm for computing the asymptotically optimal policy when model parameters are not known and backorders are not observable?
\end{itemize}

\subsection{\sf Main Results and Our Approach}\label{intro-1}
We first consider the case in which the demand distribution and the customers' patience parameter $p$ are available to the decision maker, and  the system states including the on-hand inventory level, the backorders and the inventory on-order (when the lead time is positive) can be observed by the decision maker. When the lead time is zero, we  use  dynamic programming to show that the value function  is quasi-concave and that the base-stock policy is optimal.
Our analysis shows that for all values of the patience parameter $p \in [0, 1]$, a  $p$-dependent base-stock policy is {\it uniformly optimal} for  partial backorder inventory systems, unifying two well-known results for backorder and lost-sales systems (\cite{karlin1958}, and \cite{zipkin2000}) to the partial backorder case. When the lead time is positive, however, the results in Zipkin (2008a;b) for the lost-sales system suggest that it is difficult to characterize the structure of the optimal policy for partially backordered systems using dynamic programming. For this reason, we show instead that a $p$-dependent base-stock policy ($p\in[0, 1]$) is uniformly asymptotically optimal when the sale price gets large. The dependence  on $p$ shows us that  customer patience is an important parameter in the partially backordered system.

We prove the asymptotic optimality of the base-stock policy using a  different approach from the one used by \cite{huhetal2009b} for the lost-sales system.
In particular, since the Markov chain associated with the partially backordered system cannot be shown to be uniformly ergodic, we  work directly with the ratio of the base-stock level to the sale price instead of the shortage to holding costs as was done in \cite{huhetal2009b}.

We next consider the case where the demand distribution and patience parameter $p$ are known, but backordered demands are unobservable. Unobservability occurs naturally in the real-world systems because of demand censoring, which hides the number of customers who missed out during a stockout, and customers losing patience and leaving unannounced. We call this the {\it partially observed partial backorder inventory system}.  Our decision at each period now only depends on  observations
of the past and current on-hand inventory levels and the inventory on-order (when the lead time is positive). Building on the analysis of the perfect observation case, we propose a family of $(s,q)$-policies  where $q$ can be thought of as an order for the unobserved backordered demand. We show that the $(s,q)$-policy is asymptotically optimal.

To establish the asymptotic optimality of the $(s,q)$-policy, we propose a new method that combines (a) the asymptotic optimality of the base-stock policy for the fully observed and partial backorder systems discussed above; (b) sample-path analysis; and (c) the Lyapunov function defined on the state space of the Markov chain associated with the partially observed  partial backorder system under the $(s,q)$-policy. The sample-path analysis provides an upper bound for the backorders given by the backorder of the backlogging system. The Lyapunov function not only guarantees the Markov chain is ergodic but convergence of the expected Lyapunov function, which can be used to give the long-run average backorder cost for the partially observed  partial backorder systems.

Finally, we consider the case when the demand distribution and patience parameter are both not known by the decision maker and backorders are not observed. Now, the decision maker must initiate the inventory decision from scratch, engaging in a dynamic learning process to develop an (near) optimal policy while balancing profit on the fly. This exploration-exploitation trade-off has been extensively studied in the multi-arm bandit (MAB) literature and efficient methods have been developed to handle the trade-off. For example, upper confidence bound (UCB) algorithms choose actions optimistically using the upper confidence bound as a proxy of the true parameter (\cite{aueretal2002}, \cite{auer2002}), the $\varepsilon$-greedy algorithm would force an exploration by randomly choosing actions once in a while (\cite{aueretal2002}), and the Thompson Sampling (TS) selects actions based on samples from the posterior (\cite{agrawal2012}). These algorithms have been proven to be asymptotically optimal in the standard MAB setting. Inspired by the MAB literature and the literature about learning algorithms in inventory control (\cite{huhetal2009a}, \cite{agrawal2019}, and \cite{lyuetal2021}), we design a UCB-type online algorithm based on the geometric ergodicity results established in this paper, which allows the decision maker to learn the optimal $(s,q)$-policy with the regret in the order of $O(N^{\frac{1}{2}+\delta}\log N)$ for arbitrary small $\delta>0$. This additional $\delta>0$ can be taken as a ``penalty" of the degeneration from uniform to geometric ergodicity (with this stronger notion of ergodicity, the line-search, UCB and stochastic-gradient-decent type algorithms should be able to achieve $O(\sqrt{N}\log N)$-regret, as shown in \cite{agrawal2019}, \cite{zhangetal2020}, and \cite{lyuetal2021}).

The approach used to obtain the regret bound of our UCB algorithm is similar to the method used by \cite{huhetal2009a}, which is also based on the convergence rate of the Markov chain characterizing the system state under the $(s,q)$-policy. However, in contrast to \cite{huhetal2009a}, we cannot prove the uniform ergodicity of the corresponding Markov chain due to the overshooting feature in the partial backorder system. Uniform ergodicity breaks down because the state space for the lost-sales system under the base-stock policy is compact while the state space for the partial backorder system under the $(s,q)$-policy is unbounded. Nevertheless, using results in \cite{meyntweedie1994} on the geometric ergodicity of Markov chains, we construct a petite set and prove the irreducibility and aperiodicity  to obtain the geometric ergodicity for the Markov chains corresponding to any $(s,q)$-policy.

\subsection{\sf Brief Literature Review}\label{intro-2}
Our paper relates to three streams of literature:  stochastic inventory management, online inventory control, and multi-armed bandits.

\noindent
{\sf Stochastic Inventory Management:} \ For zero-lead time, periodic review systems, \cite{karlin1958} and \cite{zipkin2000} show that the base-stock policy is optimal for backorder and lost-sales systems by a dynamic programming (DP) approach. This DP approach can be carried over to the system with positive lead times for backorder systems. However, it fails for the lost-sales system with positive lead times. Moreover, \cite{karlin1958}, and Zipkin (2008a;b) show that it is difficult to get such structured optimal policy in lost-sales systems when lead time is positive. Even though the exact optimal policy is complex in the lost-sales systems with positive lead times, effective heuristics has been developed and the effectiveness has been proved in different limit regimes. \cite{huhetal2009b} show that the base-stock policy performs well in the large lost-sales penalty regime. \cite{goldbergetal2016} and \cite{xingoldberg2016} show the constant order policy is asymptotically optimal in the large lead time regime and the optimality gap has been shown to be exponentially decay as the lead time increases.

There are also some literature that focus the structure property of the base-stock policy. \cite{janakiramenroundy2004} show that the long-run average cost under base-stock policy is convex with respect to the base-stock level. \cite{janakiramanetal2007} compare the cost of the backorder and lost-sales systems when both systems use the base-stock policy.

Even with the large body of literature, all the literature mentioned above focuses on the models where customers' reactions to the stock out are deterministic, that is, either infinity patient (backorder model) or infinity impatient (lost-sales model), and the customers' patience behavior is completely ignored. One exception is a recent paper by \cite{buetal2023} who allow unmet demands to abandon with an i.i.d. fraction abandonment. This is different from our paper where each unmet demand makes its own abandonment decision each period. At the same time, our paper assumes that the backorder demand is unobservable to the decision maker. 

\noindent
{\sf  Online Inventory Control:} \ The body of the literature developing efficient algorithms for online inventory control grows rapidly in recent years. For single-product lost-sales model with positive lead time, \cite{huhetal2009a} first develop a gradient-descent type algorithm and prove an $O(N^{2/3})$ regret. \cite{zhangetal2020} then introduce another gradient-descent type algorithm with $O(\sqrt{N})$ regret, and  \cite{agrawal2019} improve the regret by showing a linear dependence on the lead time. The benchmark in \cite{huhetal2009a}, \cite{zhangetal2020}, and \cite{agrawal2019} are the optimal base-stock policy. \cite{lyuetal2021} proposes a UCB-type algorithm incorporating KM-estimator and simulation techniques. In \cite{lyuetal2021}, both optimal base-stock and capped base-stock policy are considered as benchmark and in both cases the regret is proved in the order of $O(\sqrt{N})$.

Online inventory control are also considered some related inventory systems including inventory systems with multi-products and warehouse-capacity constraints (\cite{shietal2016}), perishable inventory systems (\cite{zhangetal2018}), including joint pricing and inventory control (\cite{chenetal2019}, \cite{chenetal2021}, \cite{chenetal2022}).

\noindent
{\sf  Multi-Armed Bandit:} \
The MAB problem has been extensively studied in the literature as one of the most classic model studying the famous exploration-and-exploitation trade-off. Early study including \cite{robbins1952}, \cite{lairobbins1985}. More recently, different efficient algorithms have been developed. For example, \cite{auer2002} and \cite{aueretal2002} develop the UCB algorithm following the frequentist while in \cite{agrawal2012}, a Thompson sampling algorithm is proposed under the Bayesian methodology. We refer readers to \cite{lattimoreszepesvari2020} for a comprehensive review.

\subsection{\sf Paper Structure and Notations}\label{notation}

\noindent
{\sf Paper Structure:} \ The paper is organized as follows: in the next section, we introduce the partial backorder inventory model, and make a connection with the backlogging and lost-sales inventory models. Section \ref{policy} establishes the asymptotic optimality of the base-stock policy for the partial backorder system. In Section \ref{partial-ob}, we consider the asymptotic optimality of the $(s,q)$-policy for the partially observed  partial backorder system. Section \ref{algorithm} concerns online inventory control for the partial backorder system with a UCB-algorithm.
In Section \ref{num} conducts  numerical experiments for the base-stock policy, $(s,q)$-policy, and UCB algorithm. The paper is concluded in Section \ref{conclusion}. All the proofs except for the theorems can be found in the Appendix.\\[-0.20in]

\noindent
{\sf Notation:} Throughout this paper, let $\mathbb{N}$ and $\mathbb{N}_+$ be the sets of all integer and nonnegative integer numbers, respectively.
For a sequence of real numbers $\{a_i:i\geq 1\}$, we use $a_{[i,k]}$ to represent $\sum_{j=i}^ka_j$ with the convention that $a_{[i,k]}=0$ when $i>k$.
 For any real numbers $x$ and $y$, let $x\wedge y=\min(x,y)$,
$x\vee y=\max\{x,y\}$, $x^+=\max\{x,0\}$, and $x^-=\max\{-x,0\}$. $\mathbb{I}_{\{A\}}$ is an indicator function with  $\mathbb{I}_{\{A\}} = 1$ if $A$ is true and
$\mathbb{I}_{\{A\}}=0$ otherwise. $|A|$ is the cardinality of set $A$.

\section{\sf Modeling Framework and Preliminaries}\label{model}
In this section, we introduce our model,  show its relation to existing backlogging and lost-sales inventory models, and explore its  sample-path properties.
To begin we build a model to describe the partial backorder feature for the inventory problem.

\subsection{\sf Modeling Framework}\label{model-1}
We consider a periodic-review, single-product inventory system with a constant lead time $\tau\in \mathbb{N}_+$. The demand at period-$i$ is denoted by $D_i$, and
$\{D_i: i\geq 1\}$ are assumed to be i.i.d. nonnegative integer-valued random variables with mean $d$. Let $I_i$ be the on-hand inventory at the beginning of period-$i$, and $U_i$ the unmet demand at period-$i$ after demand $D_i$ is realized. Of the unmet demand $U_i$, a quantity $B_i$ is backordered and $L_i$ is lost. To characterize this partial backorder/lost-sales feature, we introduce a sequence of i.i.d. mean $p$ Bernoulli random variables $\{\xi_{i;\ell}: \ell\geq 1\}$ for period-$i$ where
\begin{eqnarray}
B_i=\xi_{i; [1,U_i]} \ \ \mbox{and} \ \ L_i=U_i-\xi_{i; [1,U_i]}.\label{zhang-partial-1}
\end{eqnarray}
 We also assume that the sequences $\{\xi_{i;\ell}: \ell\geq 1\}$ for $i\geq 1$ are independent of each other and independent of demands $\{D_i: i\geq 1\}$,  though this is clearly not always the case. 
Let $Q_i$ be the order quantity at period-$i$ delivered at the start of period-$(i+\tau)$ where lead time $\tau$ is a non-negative constant.

The sequence of events in period-$i$ for the {\it system with partial backorder} is described as follows:
\begin{itemize}
    \item [(i)] The on-hand inventory $I_i$ and backorder $B_{i-1}$ are observed where $I_i\cdot B_{i-1}=0$;
    \item [(ii)] The replenishment order due in period-$i$, $Q_{i-\tau}$, arrives;
    \item [(iii)] Based on $(I_i,B_{i-1},Q_{i-\tau},\ldots,Q_{i-1})$, the order quantity $Q_i$ is determined;
    \item [(iv)] Demand $D_i$ is realized;
    \item [(v)] Demand $D_i$ plus backorder $B_{i-1}$ will be satisfied by $I_i+Q_{i-\tau}$;
    \item [(vi)] Unmet demand $U_i=(B_{i-1} +D_i-I_i-Q_{i-\tau})^+$, the backorder at period-$(i+1)$, $B_{i}=\xi_{i;[1,U_i]}$, and lost-sales at period-$i$,
    $L_{i}={U_i}-\xi_{i;[1,U_i]}$, are realized.
\end{itemize}

Starting with an initial state $X_1 \triangleq (I_1, B_0,Q_{1-\tau},\ldots,Q_0)$, $N_i \triangleq I_i+Q_{i-\tau}-B_{i-1}$
and ${IP}_i \triangleq I_i+Q_{[i-\tau,i-1]}-B_{i-1}$ are respectively the {\it net inventory} and {\it inventory position} in period-$i$. Observe that
\[
{IP}_i=N_i+Q_{[i-\tau+1,i-1]}, \  I_{i+1}=(N_i-D_i)^+,    \ i\geq 1.
\]
Since $B_i=\xi_{i;[1,U_i]}$ and each $\xi_{i;\ell}$ is Bernoulli, the unmet demand is {\it partially backordered} in  that each backordered customer reassesses each period whether she continues to wait (backorder) or leaves (lost-sale). The patience times of the unmet demands are independent and follow a geometric distribution with parameter $p$. There is only one decision at each period, the order quantity.  At period-$i$ ($\geq 1$), the order quantity $Q_i$ depends only on $X_i\triangleq (I_i,B_{i-1},Q_{i-\tau},\ldots,Q_{i-1})$. The sequence of such order quantities, $(Q_1,Q_2,\ldots)$ denoted by $\pi$, is called an {\it admissible order policy}. The set of all such admissible order policies is denoted by ${\cal A}$. The order quantity, performance measures, and states all depend on the order policy. To highlight the policy dependence, the superscript ``$\pi$" is added  when policy $\pi$ is implemented, for instance, the backorder at period-$i$ is written by $B_i^\pi$. For each $\pi\in {\cal A}$,  independence of the demands $\{D_i:i\geq 1\}$ and the partial backorder variables $\{\xi_{i;\ell}:\ell\geq 1\}$ over different periods implies  that  the sequence of systems states
$\{X^\pi_i: i\geq 1\}$ is a Markov chain with state space ${\cal X}^\pi\subseteq \mathbb{N}_+^{\tau+2}$.

The {\it base-stock policy} with {\it level $s$} is an admissible policy where the order quantity in period-$i$ 
\begin{align}
Q^s_i = (s-{IP}^s_i)^+,\label{bsp-1}
 \end{align}
tries to keep the inventory position ${IP}^s_i$ at a constant level $s$. The state space of the corresponding Markov chain $\{X^s_i:i\geq 1\}$ is ${\cal X}^s=\{(x_1,x_2,\ldots, x_{\tau+2}): x_i\in\mathbb{N}_+ \ \mbox{for $1\leq i\leq \tau+2$ and $x_1x_2=0$}\}$.

In each period, we receive revenue and pay costs according to the number of units that are sold and leftover.   Let
$r$ and $h$ be the sale price and holding cost per unit and per period, respectively. Under policy $\pi\in {\cal A}$, the expected profit at period-$i$ given
the backorders $B_{i-1}^\pi$ at period-$(i-1)$,  on-hand inventory $I^\pi_i$, and order quantity $Q^\pi_{i-\tau}$ at period-$(i-\tau)$,  is given by
\begin{align}
&\Pi^\pi_i(I^\pi_i,B^\pi_{i-1},Q^\pi_{i-\tau}) =\mathbb{E}_{(I^\pi_i, B^\pi_{i-1},Q^\pi_{i-\tau})}\Big[r\cdot \left((I^\pi_i+Q^\pi_{i-\tau})\wedge (B^\pi_{i-1}+D_i)\right)-h\cdot \left(I^\pi_i+Q^\pi_{i-\tau}-B^\pi_{i-1}-D_i\right)^+\Big].\label{1-period-profit}
\end{align}
Here $\mathbb{E}_{ (I^\pi_i, B^\pi_{i-1},Q^\pi_{i-\tau})}$ is the expectation with respect to $D_i$ for given $(I^\pi_i, B^\pi_{i-1},Q^\pi_{i-\tau})$. Given the initial state $X_1=(I_1,B_0, Q_{1-\tau}, Q_{2-\tau},\cdots,Q_0)$ denoted by ${\bf x}$, the long-run average profit denoted by ${\cal C}^\pi( r,h; {\bf x})$ under policy $\pi$ is defined as
\begin{equation}
    \begin{aligned}
     {\cal C}^\pi(r,h; {\bf x})=\liminf_{N\rightarrow \infty} \frac{1}{N}\mathbb{E}\Big( \sum_{i=1}^N \Pi^\pi_i(I^\pi_i,B^\pi_{i-1},Q^\pi_{i-\tau})\Big| X_1={\bf x}\Big). \label{hanqin-partial-2}
    \end{aligned}
\end{equation}
Our objective is to find a policy $\pi^*\in {\cal A}$ to maximize the long-run average profit.

\subsection{\sf Relation with Backlogging and Lost-sales Inventory Models}\label{model-2}
The partially backordered inventory system introduced above ({\bf P}-system for shorthand) is quite general in the following sense. When the mean $p$ of Bernoulli random variable $\xi_{i;1}$ is one, the unmet demands become infinitely patient and never abandon; there are no lost-sales. In this case, $B_i=U_i$ for $i\geq 1$ and the {\bf P}-system degenerates into the traditional backlogging inventory system, which we abbreviate as {\bf B}-system. On the other hand,  the unmet demands become impatient and   abandon instantaneously when $p=1$, there is no backorder. Now, $L_i=U_i$ for $i\geq 1$ and the {\bf P}-system degenerates into the traditional lost-sales inventory system ({\bf L}-system for shorthand). In the following, when the {\bf P}-system degenerates into {\bf B}-system ($p=1$) or {\bf P}-system ($p=0$), we add a ``bar" or ``hat", respectively, to the  corresponding system performance measures and order policy,  e.g., $\overline{Q}_i$ and $\widehat{Q}_i$ for the order quantities at period-$i$ for the {\bf B}- and {\bf L}-systems, respectively.

Compared to the  {\bf L}- and {\bf B}-systems, the analysis for the {\bf P}-system is more challenging.  Consider the base-stock policy with level $s$ (simply call $s$-policy in the rest of the paper) for a given initial state $X_1=(I_1,B_0,Q_{1-\tau},\ldots,Q_0)$ with $B_0=0$ and $I_1+Q_{[1-\tau,0]}=s$. The {\bf L}-system always has $\widehat{IP}^s_i+\widehat{Q}^s_i=\widehat{I}^s_i+\widehat{Q}^s_{[i-\tau,i]}=s$, while for the {\bf B}-system, $\overline{IP}^s_i+\overline{Q}^s_i=\overline{I}^s_i+\overline{Q}^s_{[i-\tau,i]}-\overline{B}^s_{i-1}=s$ when $i\geq 1$ with the convention that $\widehat{I}^s_1=\overline{I}_1^s=I_1$, $\overline{Q}^s_k=\widehat{Q}^s_k=Q_k$ for $1-\tau\leq k\leq 0$, and $\overline{B}^s_0=B_0$.
However, the {\bf P}-system with positive lead time $\tau$ is quite different. When the base-stock policy $s$ is implemented, if there are backorders, some orders are reserved for backlogged customers. However, as the lead time is positive, the customers may not be patient enough to wait for their order to arrive and  the system can ``overshoot" in the sense that the inventory position at period-$i$, given by  ${IP}_i^s = I_i^s+Q^s_{[i-\tau+1,i-1]}-B^s_{i-1}$ with $I^s_1=I_1$ and $Q^s_{k}=Q_k$ for $1-\tau\leq k\leq 0$, can be much larger than $s$. Notice that this situation mainly comes from the {\it random patience time} and never happens in the {\bf L}- and {\bf B}-systems where the patience time is either zero ({\bf L}-system) or infinity ({\bf B}-system).

Consequently, the main difference between the {\bf B}-, {\bf L}-, and {\bf P}-systems under the $s$ policy is that the {\bf B}- and {\bf L}-systems satisfy 
$\overline{IP}^s_i+\overline{Q}^s_{i}=\widehat{IP}^s_i+\widehat{Q}^s_{i} =s$, whereas we only have ${IP}^s_i+Q^s_{i}\ge s$ for the {\bf P}-system.
The reason is that $IP_{i-1}^s\le s$ at period-$(i-1)$ so the order $Q^s_{i-1}=s-IP_{i-1}^s$ will pull the inventory position back to $s$; however, if at period-$(i-1)$, $N^s_{i-1}<0$ and $N^s_{i-1}+Q^s_{[i-\tau,i-1]}=s$, and $L_{i-1}^s= D_{i-1}-N^s_{i-1}$, which essentially means that all unmet demand after the realization of $D_{i-1}$ loses patience and leaves, then $I^s_{i}=B^s_{i-1}=0$, which implies $N_{i}^s=I^s_{i}+Q^s_{i-\tau}-B^s_{i-1}=Q^s_{i-\tau}$ and $Q^s_{i}=\big(s-IP_i^s\big)^+=\big(s-(N^s_{i}+Q^s_{[i-\tau+1,i-1]}\big)^+=\big(s-Q^s_{[i-\tau,i]}\big)^+=0$, and hence $IP^s_i+Q^s_{i}=N^s_{i}+Q^s_{[i-\tau+1, i]}=Q^s_{[i-\tau,i]}=s-N^s_i>s$, which results in overshooting. Formally, the overshooting level under the $s$-policy at period-$i$ is defined by
\begin{eqnarray}
O^s_i=(IP^s_i-s)^+. \label{12-14-1}
\end{eqnarray}

With consideration of the overshooting issue for the {\bf P}-system, now we establish a connection between the net inventory of the {\bf B}-system and {\bf P}-system under the same base-stock policy.

\begin{proposition}\label{Lemma_Connect P-system and B-system_Dynamic}
Suppose {\bf B}- and {\bf P}-systems start with the same state $(I_1,B_0,Q_{1-\tau},\ldots,Q_0)\in {\cal X}^s$ at the beginning of the initial period, and are both operated under the same base-stock policy with level $s$. Then we have that with probability one,
\begin{equation}
\begin{aligned}
\overline{N}^s_i\le N_i^s\le \overline{N}_i^s+L^s_{[(i-\tau)\vee 1, i-1]}+O_{(i-\tau)^+}^s \ \mbox{with $O_0^s=0$}, \ \mbox{and} \ B^s_i+L^s_i\le\overline{B}_i^s \ \ \mbox{for all $i\ge 1$}.
\end{aligned}\label{connect}
\end{equation}
\end{proposition}

Next we explore the sample-path relation between the overshooting and backorders for the {\bf P}-system.

\subsection{\sf Some Sample-path Properties}\label{model-3}
 While overshooting can happen because backlogged customers are impatient, its level can be controlled by the size of the backorder. 
\begin{proposition}\label{estimation of the overshooting_lemma}
Consider {\bf P}-system with positive lead time $\tau>0$. Let the system operate under the base-stock policy with level $s$, and start with the state $X_1^s=(I^s_1,B^s_0,Q^s_{1-\tau},\ldots,Q^s_0)\in {\cal X}^s$  at the beginning of the initial period. For each $i\geq 1$, define
$\ell(i)=\max\{k\leq i: I^s_k+Q^s_{[k-\tau,k-1]}-B^s_{k-1}\leq s\}$ with the convention in which $\ell(i)=1$ if the set is empty,
then with probability one,
    \begin{align*}
O_i^s\leq  \min\Big\{B_{\ell(i)-1}^s, \ \left(B_{\ell(i)-1}^s-D_{[\ell(i)+\tau,i-1]}\right)^+\cdot \mathbb{I}_{\{\ell(i)<i-\tau\}}+B_{\ell(i)-1}^s\cdot \mathbb{I}_{\{\ell(i)\geq i-\tau\}}\Big\}.
    \end{align*}
\end{proposition}

Next we discuss the (asymptotic) optimality about the base-stock policy for the {\bf P}-system.
\section{\sf (Asymptotically) Optimal Policy for the {\bf P}-system}\label{policy}

 It is well-known that when the lead time $\tau$ is zero, the base-stock policy is optimal for both {\bf L}- and {\bf B}-system (see, e.g., \cite{zipkin2000}). We first show a similar result for the {\bf P}-system.

\subsection{\sf Optimal Policy for the {\bf P}-system: Zero Lead Time}\label{optimal-zero-leadtime}

When the lead time is zero, the state space reduces to the two-dimensional nonnegative-integer value. Formally, at period-$i$, the first coordinate is the on-hand inventory ($I_i$); and second one is the backorder ($B_{i-1}$). Since $(N_i)^+=I_i$ and $(N_i)^-=B_{i-1}$, we can use $N_i$ as the system state for the zero-lead-time {\bf P}-system, which changes the state space from $\mathbb{N}_+^2$ into $\mathbb{N}$. The system dynamics under order policy $\pi=(Q_1,Q_2,\ldots)$ are given by
\[
N_i=N_{i-1}+Q_{i-1}-D_{i-1}+(U_{i-1}-\xi_{i-1;[1,U_{i-1}]}) \ \ \mbox{with $U_{i-1}=(D_{i-1}-N_{i-1}-Q_{i-1})^+$}, \ i\geq 2,
\]
where $N_1$ is the initial state. 

As usual, we use the vanishing discount approach to find the optimal policy for the problem given by
(\ref{hanqin-partial-2}). To this end, consider the discount version of (\ref{hanqin-partial-2}) with discount factor $\alpha$:
\begin{align*}
{\cal C}^{\alpha}(r,h; x) &= \max_{\pi\in {\cal A}}\mathbb{E} \Big\{\sum_{i=1}^\infty  \alpha^{i-1}\Big[
\mathbb{E}_{N_i}\Big( r \left[ Q_i\wedge (N_i)^- +(N_i+Q_i)^+\wedge D_i\right] - h(N_i+Q_i-D_i)^+\Big)\Big]\Big|N_1=x\Big\}.
\end{align*}

Define the one-period profit by
\begin{align}
{\cal L}(x,Q)=\mathbb{E}\Big(r \left[ Q\wedge (x)^- +(Q+x)^+\wedge D_1\right]-h\cdot(x+Q-D_1)^+\Big).
\label{hanqin-3}
\end{align}
Then we have ${\cal C}^\alpha(r,h; x)= \max_{\pi\in {\cal A}}\mathbb{E} \Big\{\sum_{i=1}^\infty  \alpha^{i-1}{\cal L}(N_i,Q_i) \Big|N_1=x\Big\}$. 
This gives us the corresponding dynamic programming (DP) equation:
\begin{eqnarray}
{\cal C}^\alpha(r,h; x)=\sup_{Q \in \mathbb{N}_+}\Big\{{\cal L}(x,Q)+\alpha \cdot \mathbb{E}{\cal C}^\alpha\Big(r,h; (x+Q-D_1)^+-\xi_{1; [1,(D_1-x-Q)^+]}\Big)\Big\}. \label{hanqin-4}
\end{eqnarray}

Using the one-period profit given by (\ref{hanqin-3}), the above DP equation can be equivalently written as
\begin{eqnarray}
&& {\cal C}^\alpha(r,h;x)-rx^-=\sup_{Q \in \mathbb{N}_+}\Big\{-r\cdot (Q+x)^{-}+\mathbb{E}\Big(r \left[(Q+x)^+\wedge D_1\right]-h\cdot(x+Q-D_1)^+\nonumber\\
&& \  \ \  +\alpha r\cdot\mathbb{E}\Big[\Big((x+Q-D_1)^+-\xi_{1;[1,(D_1-x-Q)^+]}\Big)^-\Big]\label{hanqin-5}\\
&&  \ \ \ +\alpha \cdot \mathbb{E}\Big[ {\cal C}^\alpha\Big(r,h; (x+Q-D_1)^+-\xi_{[1;(D_1-x-Q)^+]}\Big)-r\cdot \Big((x+Q-D_1)^+-\xi_{1;[1,(D_1-x-Q)^+]}\Big)^-\Big]\Big\}. \nonumber
\end{eqnarray}
Letting ${\cal C}^{\alpha,o}(r,h; x)={\cal C}^\alpha(r,h;x)-rx^-$ and ${\cal L}^o(x)=-rx^-+\mathbb{E}[r (x^+\wedge D_1)-h\cdot(x -D_1)^+ + \alpha r\cdot ((x-D_1)^+-\xi_{1;[1,(D_1-x)^+]})^-]$, then the DP equation (\ref{hanqin-5}) can be further written as
\begin{eqnarray}
 {\cal C}^{\alpha,o}(r,h; x)=\sup_{y\geq x}\Big\{{\cal L}^o(y)+\alpha\cdot \mathbb{E}\left[ {\cal C}^{\alpha,o}\Big(r,h; (y-D_1)^+-\xi_{1;[1,(D_1-y)^+]}\Big)\right]\Big\}.
 \label{hanqin-6}
 \end{eqnarray}
 This modified DP equation can be completely solved by the following proposition.

\begin{proposition}\label{optimal policy of zero lead time_Theorem}
Both ${\cal L}^o(\cdot)$ and ${\cal C}^{\alpha,o}(r,h; \cdot)$ are quasi-concave.
Further, their global maximizers are identical, and given by $s^{*,\alpha}=\arg\max\limits_{x\in\mathbb{N}}{\cal L}^o(x)$.
\end{proposition}

Note that $s^{*,\alpha}=\min\Big\{x: \mathbb{P}[D_1\leq x]\geq \frac{(1-\alpha p)r }{(1-\alpha p)r + h}\Big\}$.
Using the vanishing discount approach (see \cite{puterman2014}), asking $\alpha$ to approach one, we immediately have the following result.

\begin{theorem}\label{zero-leadtime}
The base-stock policy with level $s^*$ given by $\min\Big\{x: \mathbb{P}[D_1\leq x]\geq \frac{(1-p)r}{(1-p)r + h}\Big\}$
is optimal for the {\bf P}-system with zero lead time.
\end{theorem}

\subsection{\sf Asymptotically Optimal Policy: Positive Lead Time} \label{Positive lead time_sec}


When the lead time is positive, the structured optimal policy is hard to find.  On the other hand, \cite{huhetal2009b} show that the base-stock policy, which is a special case of our model, is asymptotically optimal for the {\bf L}-system in the limiting regime where the penalty cost dominates the holding cost.

Recalling the overshooting phenomenon for the {\bf P}-system in Section \ref{model}, it follows that the {\bf P}-system with positive lead time is fundamentally different from the {\bf L}-, {\bf B}-, and  zero lead time {\bf P}-systems.
Nevertheless even with this critical difference, we can show that the base-stock policy works well uniformly for all {\bf P}-systems with $p\in[0,1]$. 
We focus on the case where $p<1$, as optimality of the base-stock policy  for the {\bf B}-system which is obtained when $p=1$ is well known.  As a first step to prove asymptotic optimality, we show the long-run average overshooting is bounded and asymptotically vanishing.

\begin{proposition}\label{Uniform bounded overshooting_lemma}
Assume that there exists $\gamma>0$ such that $\mathbb{E}e^{\gamma D_1}<\infty$. Under a fixed base-stock policy $s$, and an initial state
$(I^s_1,B^s_0,Q^s_{1-\tau},\ldots,Q^s_0)\in {\cal X}^s$,
$\limsup\limits_{N\to\infty}\frac{1}{N}\mathbb{E}\Big( O^s_{[1,N]}\Big)<\infty$, and $\limsup\limits_{s\to\infty}\limsup\limits_{N\to\infty}\frac{1}{N}\mathbb{E}\Big( O^s_{[1,N]}\Big)=0$.
\end{proposition}

We can then apply Proposition \ref{Uniform bounded overshooting_lemma} to show that the long-run average profit for the {\bf P}-system given by (\ref{hanqin-partial-2}) under the $s$-policy
can be bounded from below by another {\bf B}-system with the same $s$-policy, the same unit sale price and holding cost, but the backorder cost charged by
$(r+\tau h)$. The intuition is that the backorder cost in the {\bf B}-system is so high that it is never profitable to serve the backorders who would be lost in the {\bf P}-system. To highlight this cost change, for a given initial state ${\bf x}$, the long-run average profit for this {\bf B}-system is denoted by
$\overline{\cal C}^{s}(r,h,r+\tau h; {\bf x})$.
Formally, the two long-run average profits for these two systems have the following relation:

\begin{proposition}\label{upper bound of base-stock P-system_prop}
For any initial state ${\bf x}\in {\cal X}^s$, we have
\[
{\cal C}^{s}(r,h; {\bf x}) \ge \overline{\cal C}^s(r, h, r+\tau h; {\bf x})-\limsup\limits_{N\to\infty}\frac{1}{N}h\cdot\mathbb{E}O^s_{[1,N]}.
\]
\end{proposition}

Next, we construct an upper bound for the optimal long-run average profit of the {\bf P}-system. To do this, we couple the {\bf P}-system with a {\bf B}-system  that has a per-period unit backorder cost that is so small that it is profitable to serve each backordered customer who would otherwise have abandoned in the {\bf P}-system.  

Here, we extend the technique introduced in \cite{janakiramanetal2007} that enables us to simultaneously match the on-hand inventory, backlogged demand, and lost-sales in the {\bf P}-system with {\bf B}-system. We will divide the backlogged demand in the {\bf B}-system into two groups, which we refer to as ``regular backlogged demand" and ``lost backlogged demand", which will match the backlogged demand and lost-sales in the {\bf P}-system. This enables us to show that the optimal long-run average profit of the {\bf P}-system can be upper bounded by the long-run average profit for some {\bf B}-system.

\begin{proposition}\label{lower bound of optimal P-system_prop}
For any initial state ${\bf x}\in {\cal X}^c$, let $\pi^*=\argmax_{\pi \in {\cal A}}{\cal C}^\pi(r,h; {\bf x})$. Then
${\cal C}^{\pi^*}(r,h; {\bf x})\leq \max_{s} \overline{\cal C}^s(r,h,\frac{(1-p)r}{2(1+\tau)};{\bf x})$.
\end{proposition}

With the help of Propositions \ref{upper bound of base-stock P-system_prop} and \ref{lower bound of optimal P-system_prop}, we can obtain the
asymptotic optimality for the base-stock policy in the {\bf P}-system.

\begin{theorem}\label{symptotic-s-policy} $($\textbf{Asymptotic Optimality of Base-Stock Policies}$)$
Assume that $\mathbb{E}e^{\gamma D_1}<\infty$ for some $\gamma>0$. Then for any initial state ${\bf x}$,
\begin{align*}
\frac{\max\limits_{s}{\cal C}^s(r,h;{\bf x})}{\max\limits_{\pi\in{\cal A}}{\cal C}^\pi(r,h;{\bf x})}  \rightarrow 1 \ \ \mbox{as} \
 r \rightarrow \infty.
 \end{align*}
\end{theorem}

\noindent
{\it Proof}: First by Proposition \ref{upper bound of base-stock P-system_prop},
\begin{align*}
 \max_{s}{\cal C}^s(r,h; {\bf x})\geq \bar{\cal C}^{\overline{s}^*_1}(r, h, r+\tau h; {\bf x})-\limsup\limits_{N\to\infty}\frac{1}{N}h\cdot\mathbb{E} O^{\overline{s}^*_+}_{[1,N]},
 \end{align*}
where $\overline{s}^*_1=\argmax_{\overline{s}}\bar{\cal C}^{\overline{s}}(r, h, r+\tau h; {\bf x})$.
This combining with Proposition \ref{lower bound of optimal P-system_prop} gives
\begin{align}
\frac{\max_{s}{\cal C}^s(r,h; {\bf x})}{\max_{\pi\in {\cal A}}{\cal C}^\pi(r,h; {\bf x})}
&\geq \frac{\bar{\cal C}^{\overline{s}^*_1}(r, h, r+\tau h; {\bf x})-\limsup\limits_{N\to\infty}\frac{1}{N}h\cdot\mathbb{E} O^{\overline{s}^*_1}_{[1,N]}}{\max_{s} \overline{\cal C}^s(r,h,\frac{(1-p)r}{2(1+\tau)};{\bf x})}\nonumber\\
&=\frac{\bar{\cal C}^{\overline{s}^*_1}(r, h, r+\tau h; {\bf x})-\limsup\limits_{N\to\infty}\frac{1}{N}h\cdot\mathbb{E} O^{\overline{s}^*_1}_{[1,N]}}{\bar{\cal C}^{\overline{s}^*_2}(r,h,\frac{(1-p)r}{2(1+\tau)};{\bf x})},\label{hanqin-31}
\end{align}
where $\overline{s}^*_2=\arg\max_{s}\overline{\cal C}^s(r,h,\frac{(1-p)r}{2(1+\tau)};{\bf x})$. Recall that for the {\bf B}-system, we know that
\begin{align}
&\overline{s}^*_\ell=\inf\Big\{y: \mathbb{P}\Big(D_{[1,\tau+1]}\geq y\Big)\geq \frac{h}{h+b_\ell}\Big\}, \ \ \ell=1,2; \label{12-18-1}\\
&b_1=r+\tau h, \ \ \mbox{and} \ \ b_2=\frac{(1-p)r}{2(1+\tau)}.\label{12-18-6}
\end{align}
Note that  $\lim_{r\rightarrow \infty} \overline{s}^*_\ell= \inf\{y: \mathbb{P}(D_{[1,\tau+1]}\leq y)=1\}$.
Thus, it follows from Proposition \ref{Uniform bounded overshooting_lemma} that
 $\lim_{r\rightarrow \infty}\lim_{N\rightarrow \infty}\frac{1}{N}\mathbb{E} O^{\overline{s}^*_1}_{[1,N]}=0$.
 In view of (\ref{hanqin-31}), to prove the theorem, it suffices to show that
 \begin{align}
\frac{\bar{\cal C}^{\overline{s}^*_1}(r, h, r+\tau h; {\bf x})}{\bar{\cal C}^{\overline{s}^*_2}(r,h,\frac{(1-p)r}{2(1+\tau)};{\bf x})} \to 1 \quad\text{as $r\to \infty$}.\label{12-18-2}
\end{align}
Recall that
\begin{align*}
\bar{\cal C}^{\overline{s}^*_\ell}\Big(r,h,b_\ell;{\bf x}\Big)
=r\cdot \mathbb{E}D_1-h\cdot \mathbb{E}\Big(\overline{s}^*_2-D_{[1,\tau+1]}\Big)^+-b_\ell \cdot\mathbb{E}\Big( D_{[1,\tau+1]}-\overline{s}^*_\ell\Big)^+, \ \ell=1,2.
\end{align*}
Because $\mathbb{E}D_1>0$ and by (\ref{12-18-1}),
\[
\lim_{r\rightarrow \infty}\mathbb{E}\Big( D_{[1,\tau+1]}-\overline{s}^*_1\Big)^+=\lim_{r\rightarrow \infty}\mathbb{E}\Big( D_{[1,\tau+1]}-\overline{s}^*_2\Big)^+=0,
\]
to prove (\ref{12-18-2}), from $\mathbb{E}\Big(\overline{s}^*_\ell-D_{[1,\tau+1]}\Big)^+\leq \overline{s}^*_\ell$ with $\ell=1,2$, it is sufficient to prove
\begin{align}
\lim_{r\rightarrow \infty}\frac{\overline{s}^*_1}{r}=\lim_{r\rightarrow \infty}\frac{\overline{s}^*_2}{r}=0. \label{12-18-3}
\end{align}
As $\mathbb{E}e^{\gamma D_{[1,\tau+1]}}<\infty$, then the Markov inequality provides us
\begin{align}
\mathbb{P}\Big(D_{[1,\tau+1]}\geq y\Big)\leq e^{-\gamma y}\mathbb{E}e^{\gamma D_{[1,\tau+1]}}.
\label{12-18-7}
\end{align}
According to the definition of $\overline{s}^*_1$, by the decreasing property of the functions of $y$ given by both sides of (\ref{12-18-7}),
\begin{align}
\overline{s}^*_\ell&=\inf\Big\{y: \mathbb{P}\Big(D_{[1,\tau+1]}\geq y\Big)\geq \frac{h}{h+b_\ell}\Big\}\nonumber\\
&\leq \inf\Big\{y:e^{-\gamma y}\mathbb{E}e^{\gamma D_{[1,\tau+1]}} \geq \frac{h}{h+b_\ell}\Big\}\nonumber\\
& =\frac{1}{\gamma} \ln\Big(\frac{h+b_\ell}{h}\mathbb{E}e^{\gamma D_{[1,\tau+1]}}\Big).\label{12-18-8}
\end{align}
(\ref{12-18-3}) directly follows from (\ref{12-18-6}) and (\ref{12-18-8}).
Therefore, we have the theorem. \hfill$\Box$

\begin{remark}
 From the proof of Theorem \ref{symptotic-s-policy}, observe that if all what we care about is the long-run average profit when the unit sale price $r$ is high, we can blindly choose to implement $\overline{s}^*$-policy, which is given by optimizing a {\bf B}-system with the same sale price and holding cost but the backorder cost  $(1-p)r/(2(1+\tau))$. This is actually intuitive. As shown in Theorem 2 of \cite{huhetal2009b}, under a base-stock policy, when the unit sale price is high, the {\bf L}-system will be ``close" to the {\bf B}-system. It is then not surprise that our {\bf P}-system with $p\in[0,1]$ as a ``in-between" system should be sandwiched by {\bf L}- and {\bf B}-systems  and they are all ``close" to each other in the high sale-price limit regime.
 However, our result is based on the assumption $\mathbb{E}e^{\gamma D_1}<\infty$, which may not satisfy what \cite{huhetal2009b} need.
 On the other hand, as we will show in our numerical experiment (cf. Tables \ref{tab 2:optimality_gap_basestock_policy_zhaoxuan_0.7} and \ref{tab 4:optimality_gaps_basestock_policy_zhaoxuan_0.3}), when the price is comparable with respect to the holding cost, the loss of profit could be large if the behavior of the customers is ignored and the policy $\overline{s}^*$ is blindly chosen.

 It seems that the value of considering the customers' abandonment behavior should vanish in the limiting regime with high sale-price. However, from the numerical experiments shown in
 Tables \ref{tab 1:avg_profit_and_levels_basestock_policy_zhaoxuan_0.7} and \ref{tab 3:avg_profit_and_levels_basestock_policy_zhaoxuan_0.3}, we see that even in the regime where the unit sale price dominates the holding cost, choosing the base-stock policy $\overline{s}^*$ requires more thought. Since the decisions under the optimal  $s^*$-policy  could still be different from the decisions given by the $\overline{s}^*$-policy. For example, our numerical results in Tables \ref{tab 1:avg_profit_and_levels_basestock_policy_zhaoxuan_0.7} and \ref{tab 3:avg_profit_and_levels_basestock_policy_zhaoxuan_0.3} show that even in asymptotic regime ($h/r=1/64\approx 1.6\%$), the relative difference between $s^*$ and $\overline{s}^*$, defined by $(\overline{s}^*-s^*)/s^*$, could still be around $7\%-10\%$. This will then provide the system more flexibility under the optimal base-stock policy $s^*$ compared to the $\overline{s}^*$-policy.
\end{remark}

\section{\sf Partially Observed {\bf P}-System}\label{partial-ob}
We showed in Section \ref{policy} that the base-stock policy is asymptotically optimal for the {\bf P}-system when the on-hand inventory and backorders are observable to the decision maker. However, it is commonly the case that the decision maker is unable to perfectly observe new customer arrivals, abandonments and the number of backorders during a stock-out period. As a result, the inventory position is not known and the base-stock policy given by (\ref{bsp-1}) is unimplementable. We refer to this as a {\it partially observed {\bf P}-system} ({\bf POP}-system for short hand).

For the {\bf POP}-system with a given initial state $X_1=(I_1,B_0, Q_{1-\tau},\ldots,Q_0)$, based on $(I_1, Q_{1-\tau},\ldots,Q_0)$ obtained from $X_1$ by taking out the backorder $B_0$, we need to determine the order quantity $Q_1$, which is delivered at period-$(1+\tau)$. In general, for state $X_i=(I_i,B_{i-1},Q_{i-\tau},\ldots,Q_{i-1})$ at period-$i$, only based on $(I_i,Q_{i-\tau},\ldots,Q_{i-1})$, we determine the order quantity $Q_i$ which is delivered at period-$(i+\tau)$. Although the backorders over periods are not observable, our goal for the {\bf POP}-system is still to find an order policy to maximize the long-run average profit given by (\ref{hanqin-partial-2}).

Because of the information lost about the backorders, the optimal long-run average profit of the {\bf POP}-system should not be larger that of the {\bf P}-system.
In view of the asymptotic optimality  of the base-stock policy for the {\bf P}-system, we develop an analog for the {\bf POP}-system.
First observe that $I_i=0$ tells us there may be backorders happen in the previous period, that is, even though $B_{i-1}$ is not observable, $B_{i-1}>0$ happens only when $I_i=0$.
Consider the base-stock policy given by (\ref{bsp-1}) and recall that $IP_i=I_i+Q_{[i-\tau,i-1]}-B_{i-1}$. Thus if $\mathbb{E}(B_{i-1}|I_i=0)$ is available, then we can use
$\Big(s-I_i-Q_{[i-\tau,i-1]}+ \mathbb{E}[B_{i-1}|I_i=0]\Big)$ as a good proxy for the base-stock policy in the {\bf P}-system. Hence, for the {\bf POP}-system, we will consider the following $(s,q)$-policy by modifying the base-stock policy given by (\ref{bsp-1}): for a given initial state $X^{(s,q)}_1=(I_1^{(s,q)},B^{(s,q)}_0,Q^{(s,q)}_{1-\tau},\ldots,Q_0^{(s,q)})$,
\begin{align}
Q^{(s,q)}_i=\Big(s+q\cdot \mathbb{I}_{\{I^{(s,q)}_i=0\}}-I^{(s,q)}_i-Q^{(s,q)}_{[i-\tau,i-1]}\Big)^+ \ \ \mbox{for $i\geq 1$}.\label{hanqin-35}
\end{align}
Note that the above $(s,q)$-policy is also an admissible order policy for the {\bf P}-system. So each $(s,q)$-policy defined by (\ref{hanqin-35}) belongs to ${\cal A}$.
Further, the sequence of the states generated by the $(s,q)$-policy, $\{X_i^{(s,q)}=(I^{(s,q)}_i,B^{(s,q)}_{i-1},Q^{(s,q)}_{i-\tau},\ldots,Q^{(s,q)}_{i-1}): i\geq 1\}$,
is a Markov chain with the state space
\begin{align*}
{\cal X}^{(s,q)}=\Big\{\Big(x_1,x_2,\ldots,x_{\tau+2}\Big):  x_1+x_{[3,\tau+2]}\leq s+q, \ x_1\cdot x_2=0, \ \mbox{and} \ x_j\in \mathbb{N}_+ \ \mbox{for $1\leq j\leq \tau+2$}\Big\}.
\end{align*}
The ergodicity of this Markov chain is given by the following proposition.

\begin{proposition}\label{ergodic base-stock policy_thm}
Assume that $\mathbb{P}[D_1=0]=\alpha_0>0$ and $\mathbb{P}[D_1=1]=\alpha_1>0$. Under the $(s,q)$-policy, the Markov chain generated by the system dynamics of the {\bf POP}-system with $p<1$, $\Big\{X_i^{(s,q)}=\Big(I^{(s,q)}_i,B_{i-1}^{(s,q)},Q_{i-\tau}^{(s,q)},\ldots,Q^{(s,q)}_{i-1}\Big): i\geq 1\Big\}$, is ergodic with the steady state distribution given by $X_\infty^{(s,q)}$. Further, the Lyapunov function can be constructed by ${\cal V}({\bf x})=p\cdot (x_2)^\delta+1$ for any $\delta\in (0, 1]$, and
$\lim_{i\rightarrow \infty}\mathbb{E}{\cal V}(X_i^{(s,q)})=\mathbb{E}{\cal V}(X_\infty^{(s,q)})$.
\end{proposition}

 Here we intuitively explain why such Lyapunov function works. Because the $(s,q)$-policy is implemented, the order quantity at period-$i$ is given by $Q_i^{(s,q)}=(s+r\cdot\mathbb{I}_{\{I_i^{(s,q)}=0\}}-I_i^{(s,q)}-Q^s_{[i-\tau,i-1]})^+$ which can be controlled $(s+q)$. At the same time, the on-hand inventory $I^{(s,q)}_i$ mainly comes the order quantities. With this observation, we know that Markov chain $\{X^{(s,q)}_i: i\geq 1\}$ is actually controlled by its second coordinate. In turn, the Lyapunov function should be type of ${\cal V}(x)=c\cdot (x_2)^\gamma+1$ for some $c,\gamma\in (0, \infty)$. The convergence of $\{\mathbb{E}{\cal V}(X_i^{(s,q)}): i\geq 1\}$ directly follows from Theorem 15.0.1 in \cite{meyntweedie2012}.

In the following two subsections, we show
the asymptotic optimality for the $(s,q)$-policy. The asymptotic optimality shows that even though the backorders get unobservable to the decision maker for the {\bf POP}-system, the $(s,q)$-policy can guarantee the asymptotic optimality and performs as good as the base-stock policy which is applicable only when the whole system information is fully observable. We first start with the zero lead time model.

\subsection{\sf Asymptotic Optimality of $(s,q)$-Policy: Zero Lead Time Case}\label{analysis of (S,r)-policy_subsection}
For the zero lead time case, in view of Theorem \ref{zero-leadtime}, it is sufficient to compare the $(s,q)$-policy with the base-stock level policy. Given an initial state $X_1=(I_1,B_0)$, consider two policies: the base-stock policy with level $s$, and the $(s,q)$-policy. Then
\begin{align*}
\begin{array}{lll}
&Q^s_1=(s-I_1+B_0)^+, &Q^s_i=(s-I_i^s-B^s_{i-1})^+ \ \mbox{for $i\geq 2$};\\
&Q^{(s,q)}_1=\Big(s+q\cdot \mathbb{I}_{\{I_1=0\}}-I_1\Big)^+, &Q^{(s,q)}_i=\Big(s+q\cdot \mathbb{I}_{\{I^{(s,q)}_{i}=0\}}-I_i^{(s,q)}\Big)^+ \ \mbox{for $i\geq 2$}.
\end{array}
\end{align*}
By (\ref{1-period-profit}), their profits at period-$i$ can be written by
\begin{align}
&\Pi^s_i(I^s_i,B^s_{i-1})=\mathbb{E}_{(I^s_i,B^s_{i-1})}\Big[r\cdot \left((I^s_i+Q^s_{i-\tau})\wedge (B^s_{i-1}+D_i)\right)-h\cdot \left(I^s_i+Q^s_{i-\tau}-B^s_{i-1}-D_i\right)^+\Big];\label{2-1-period-profit}\\
&\Pi^{(s,q)}_i(I^{(s,q)}_i,B^{(s,q)}_{i-1})=\mathbb{E}_{(I^{(s,q)}_i, B^{(s,q)}_{i-1})}\Big[r\cdot \left((I^{(s,q)}_i+Q^{(s,q)}_{i-\tau})\wedge (B^{(s,q)}_{i-1}+D_i)\right)\nonumber\\
& \ \ \ \ \ \ \ \ \ \ \ \ \ \  \ \ \ \ \ \ \ \ \ \  -h\cdot \left(I^{(s,q)}_i+Q^{(s,q)}_{i-\tau}-B^{(s,q)}_{i-1}-D_i\right)^+\Big].\label{1-1-period-profit}
\end{align}
Simply, let ${\cal R}^s_i$ and ${\cal R}^{(s,q)}_i$ denote the above profits before taking the expectations  $\mathbb{E}_{(I^s_i,B^s_{i-1})}$ and $\mathbb{E}_{(I^{(s,q)}_i,B^{(s,q)}_{i-1})}$, respectively.
Then we can make a sample path comparison between their profits.

\begin{proposition}\label{couple the diff for positive on-hand inv}
 Given any initial state $X_1=(I_1,B_0)$, for the $s$-policy and $(s,q)$-policy, we have the following comparisons$:$\\
 {\rm (i)} If $0<I^{(s,q)}_i\le s$ for $i= i_0,\ldots,i_0+k$ with $i_0\ge2$ and $k\ge 0$, then the sale amounts under these two policies are same at periods
 through $i_0+1$ to $i_0+(k\vee 1)$. Furthermore, we have $I^{(s,q)}_i = I^{s}_i$ for $i=i_0+1,\ldots, i_0+k+1$,  and
 ${\cal R}^s_{[i_0,i_0+k]} -{\cal R}^{(s,q)}_{[i_0, i_0+k]} \leq r\cdot B^{s}_{i_0-1};$\\
 {\rm (ii)} If $I_i^{(s,q)}=0$ for $i= i_0,\ldots, i_0+k$ with $i_0\ge2$ and $k\ge 0$, then we have
 ${\cal R}^s_{[i_0,i_0+k]} -{\cal R}^{(s,q)}_{[i_0,i_0+k]}\leq r\cdot B^s_{[i_0-1,i_0+k-1]} -(h+r)\cdot I^s_{[i_0+1,i_0+k]}-krq+hq;$\\
 {\rm (iii)} If $s< I_i^{(s,q)}\le s+q$ for $i= i_0,\ldots, i_0+k$ with $i_0\ge2$ and $k\ge 0$,  we then have
  ${\cal R}^s_{[i_0,i_0+k]} -{\cal R}^{(s,q)}_{[i_0, i_0+k]}\leq (k+1)(qh+r(q-s)^+)$.
 \end{proposition}

With the help of Proposition \ref{couple the diff for positive on-hand inv}, we get the asymptotic optimality of the $(s,q)$-policy.

\begin{theorem}\label{gap between benchmark and optimal base-stock}
Assume that $\mathbb{P}[D_1=0]=\alpha_0>0$ and $\mathbb{P}[D_1=1]=\alpha_1>0$.  For any initial state $X_1=(I_1,B_0)$, there exists a constant $\nu \le 1$ such that the gap between the long-run average profits given by the optimal $(s,q)$-policy and the optimal base-stock policy can be expressed as
\begin{align*}
\max_{s}{\cal C}^s(r,h; I_1,B_0)-\max_{(s,q)} {\cal C}^{(s,q)}(r,h; I_1,B_0)\leq \nu \cdot r \mathbb{E}(s^*-D_{[1,\tau+1]})^+ \ \mbox{with} \ s^*=\arg\max_{s}{\cal C}^s(r,h; I_1,B_0).
\end{align*}
Moreover, we have the asymptotical optimality about the $(s,q)$-policy
\begin{align*}
\lim_{r\rightarrow \infty}\frac{\max_{s}{\cal C}^s(r,h; I_1,B_0)-\max_{(s,q)}{\cal C}^{(s,q)}(r,h;I_1,B_0)}{\max_{s}{\cal C}^s(r,h; I_1,B_0)}=0.
\end{align*}
\end{theorem}

\noindent
{\it Proof}: For any positive integer $N\geq 3$, let
\begin{align*}
& {\cal I}_1(s,q; N)=\Big\{i: 2\leq i\leq N \ \mbox{and} \  0<I^{(s,q)}_i\le s\Big\};\\
& {\cal I}_2(s,q; N)=\Big\{i: 2\leq i\leq N \ \mbox{and} \  I^{(s,q)}_i=0\Big\};\\
& {\cal I}_3(s,q; N)=\Big\{i: 2\leq i\leq N \ \mbox{and} \  s<I^{(s,q)}_i\le s+q\Big\}.
\end{align*}
By Proposition \ref{couple the diff for positive on-hand inv},
\begin{align}
{\cal R}^s_{[1,N]} -{\cal R}^{(s,q)}_{[1,N]}&= \Big({\cal R}^s_{1}-{\cal R}^{(s,q)}_{1}\Big)+ \Big({\cal R}^s_{[2,N]}-{\cal R}^{(s,q)}_{[2,N]}\Big)\nonumber\\
&= \Big({\cal R}^s_{1}-{\cal R}^{(s,q)}_{1}\Big)+\sum_{i\in {\cal I}_1(s,q; N)}\Big({\cal R}^s_{i}-{\cal R}^{(s,q)}_{i}\Big)\nonumber\\
& \ +\sum_{i\in {\cal I}_2(s,q;N)}\Big({\cal R}^s_{i}-{\cal R}^{(s,q)}_{i}\Big)+\sum_{i\in {\cal I}_3(s,q;N)}\Big({\cal R}^s_{i}-{\cal R}^{(s,q)}_{i}\Big)\nonumber\\
& \le  \Big({\cal R}^s_{1}-{\cal R}^{(s,q)}_{1}\Big)+r\cdot \sum_{i\in {\cal I}_1(s,q;N)\cup {\cal I}_2(s,q;N)}B^s_{i-1}\nonumber\\
& \  + \Big|{\cal I}_2(s,q;N)\cup {\cal I}_3(s,q;N)\Big| \cdot hq+ \Big| {\cal I}_3(s,q;N)\Big|\cdot r(q-s)^+. \label{hanqin-40}
\end{align}
Under the $s$-policy, for any initial state $X_1=(I_1,B_0)$, we have that $B_i^s$ follows the distribution $\xi_{1,[1,U^s]}$ with $U^s=(D_{[1,\tau+1]}-s)^+$ for large $i$. Hence, using (\ref{hanqin-40}), for the $(s,q)$-policy and $s$-policy,
\begin{align}
{\cal C}^s(r,h;I_1,B_0)-{\cal C}^{(s,q)}(r,h;I_1,B_0)&\leq \Big(\nu_1(s,q)+\nu_2(s,q)\Big)r\cdot \mathbb{E}(s-D_{[1,\tau+1]})^+\nonumber\\
& +\Big(\nu_2(s,q)+\nu_3(s,q)\Big)hq+\nu_3(s,q) \cdot r(q-s)^+,\label{hanqin-41}
\end{align}
where $\nu_j(s,q)=\lim_{N\rightarrow \infty }|{\cal I}_j(s,q;N)|/{N}$ for $j=1,2,3$. The existence of these three limits  is warranted by Proposition
\ref{ergodic base-stock policy_thm}.
 Note that
\begin{align}
\max_{s} {\cal C}^s(r,h;I_1,B_0)-\max_{s,q}{\cal C}^{(s,q)}(r,h;I_1,B_0)
& \leq {\cal C}^{s^*}(r,h;I_1,B_0)-{\cal C}^{(s^*,0)}(r,h;I_1,B_0)\nonumber\\
& \le \Big(\nu_1(s^*,0)+\nu_2(s^*,0)\Big)r\cdot \mathbb{E}(s^*-D_{[1,\tau+1]})^+.\label{hanqin-42}
\end{align}
This gives the first part of the theorem. The second part of the theorem directly follows (\ref{hanqin-41})-(\ref{hanqin-42}), and $s^*\rightarrow ({\cal D}_{sup}\wedge \infty)$, which comes from Theorem \ref{zero-leadtime}.\hfill$\Box$

\subsection{\sf Asymptotic Optimality of $(s,q)$-Policy: Positive Lead Time}\label{Analysis of $(S,r)$-Policy: Model with Positive Lead Time_section}
 From the (asymptotic) optimality analysis for the base-stock policy for zero and no zero lead times in Subsections \ref{optimal-zero-leadtime} and \ref{Positive lead time_sec}, we see  some fundamental differences between the zero lead time system and the positive lead time system since the backlogged customers may lose their patience before they receive order. Thus, it is very difficult to analyze the optimality gap between the optimal policy and the $(s,q)$-policy by the sample-path method as what we have done in Subsection \ref{analysis of (S,r)-policy_subsection}. Instead, we use the result about the ergodicity established in Proposition \ref{ergodic base-stock policy_thm} to prove the asymptotic optimality of  the $(s,q)$-policy under the same regime as Theorem \ref{symptotic-s-policy}.

\begin{theorem} \label{symptotic-s-q-policy} $(${\bf Asymptotic Optimality of $(s,q)$-Policy}$)$
Assume that $\mathbb{P}[D_1=0]=\alpha_0>0$ and $\mathbb{P}[D_1=1]=\alpha_1>0$.  For any fixed initial state $X_1=(I_1,B_0,Q_{1-\tau},\ldots,Q_0):={\bf x}$,
\begin{align*}
\frac{\max_{(s,q)} {\cal C}^{(s,q)}(r,h; {\bf x})}{\max_{\pi\in {\cal A}}{\cal C}^\pi(r,h; {\bf x})}\rightarrow 1 \ \ \mbox{as $r\rightarrow \infty$}.
\end{align*}
\end{theorem}

\noindent
{\it Proof}: Consider a {\bf B}-system with the sale price $r$, the holding cost $h$, and the backorder cost $r+\tau h$. Let
$\overline{s}^*=\arg\max_s \overline{\cal C}^s(r,h,r+\tau h;x)$.
In view of the proof of of Theorem \ref{symptotic-s-policy}, it suffices to prove that
\begin{align}
\frac{{\cal C}^{(\overline{s}^*,0)}(r,h; {\bf x})-\bar{\cal C}^{\overline{s}^*}(r, h, r+\tau h; {\bf x})}{\bar{\cal C}^{\overline{s}^*}(r, h, r+\tau h; {\bf x})}\rightarrow 0 \ \ \mbox{as $r\rightarrow \infty$}.\label{hanqin-57}
\end{align}
By the ergodicity of the backlogging inventory system with any base-stock policy, the steady states of the order quantity, the on-hand inventory, and the backorder exist, and let $\bar{Q}^{\overline{s}^*}_\infty$, $\bar{I}^{\overline{s}^*}_\infty$, and $\bar{B}^{\overline{s}^*}_\infty$ represent their steady states, respectively. Then
\begin{align}
 \bar{\cal C}^{\overline{s}^*}(r, h, r+\tau h; {\bf x}) =
\mathbb{E} \Big( r\cdot \bar{Q}^{\overline{s}^*}_\infty-h \cdot \bar{I}^{\overline{s}^*}_\infty
-(r+\tau h) \bar{B}^{\overline{s}^*}_\infty\Big). \label{12-12-4}
\end{align}
Similarly, using Proposition \ref{ergodic base-stock policy_thm}, let ${Q}^{\overline{s}^*}_\infty$, ${I}^{\overline{s}^*}_\infty$, and ${L}^{\overline{s}^*}_\infty$, and ${N}^{\overline{s}^*}_\infty$ be the steady states of the order quantity, the on-hand inventory, the lost-sales, and the net inventory  under the $\overline{s}^*$-policy for the {\bf P}-system. Then
\begin{align}
 {\cal C}^{(\overline{s}^*,0)}(r, h; {\bf x}) =
\mathbb{E} \Big( r\cdot {Q}^{(\overline{s}^*,0)}_\infty-h \cdot {I}^{(\overline{s}^*,0)}_\infty\Big). \label{12-12-5}
\end{align}
Notice that for any $i>\tau$,
\begin{align}
\bar{I}_i^{\overline{s}^*}+\bar{Q}_{[i-\tau,i]}^{\overline{s}^*}-\bar{B}_{i-1}^{\overline{s}^*}=\overline{s}^*
=I_i^{(\overline{s}^*,0)}+Q_{[i-\tau,i]}^{(\overline{s}^*,0)}.  \label{12-12-9}
\end{align}
This implies that
\begin{align}
\bar{I}_\infty^{\overline{s}^*}+(\tau+1)\bar{Q}_\infty^{\overline{s}^*}-\bar{B}_\infty^{\overline{s}^*}
=I_\infty^{(\overline{s}^*,0)}+(\tau+1)Q_\infty^{(\overline{s}^*,0)}. \label{12-12-6}
\end{align}
It follows from (\ref{12-12-4})-(\ref{12-12-6}) that
\begin{align}
&\bar{\cal C}^{\overline{s}^*}(r, h, r+\tau h; {\bf x})-{\cal C}^{(\overline{s}^*,0)}(r, h; {\bf x}) \nonumber\\
& \ \ \ =\mathbb{E} \Big((r+(\tau+1)h)( \bar{Q}^{\overline{s}^*}_\infty-{Q}^{(\overline{s}^*,0)}_\infty)
-(r+(\tau+1)h)\bar{B}^{\overline{s}^*}_\infty\Big).\label{12-12-7}
\end{align}
For the {\bf P}-system, recall that for any $i>\tau$, the net inventory has the relation:
$N_{i+1}^{(\overline{s}^*,0)}=N_{i}^{(\overline{s}^*,0)}+Q_{i+1-\tau}^{(\overline{s}^*,0)}-D_i+L_{i}^{(\overline{s}^*,0)}$.
This together with $\mathbb{E}\bar{Q}^{\overline{s}^*}_\infty= \mathbb{E}D_1$ gives us
\begin{align}
\mathbb{E}\bar{Q}^{\overline{s}^*}_\infty-\mathbb{E} {Q}^{(\overline{s}^*,0)}_\infty=\mathbb{E} {L}^{(\overline{s}^*,0)}_\infty.\label{12-12-8}
\end{align}
For the {\bf P}-system, the unmet demand at period-$i$ with $i>2\tau+1$,
\begin{align*}
U_i^{(\overline{s}^*,0)}&\leq \Big(B^{(\overline{s}^*,0)}_{i-\tau-1}+D_{[i-\tau,i]}-I^{(\overline{s}^*,0)}_{i-\tau}-Q^{(\overline{s}^*,0)}_{[i-2\tau,i-\tau]}\Big)^+
\nonumber\\
&= \Big(B^{(\overline{s}^*,0)}_{i-\tau-1}+D_{[i-\tau,i]}-\overline{s}^*\Big)^+ \ \ \mbox{(by (\ref{12-12-9}))}\nonumber\\
&\le B^{(\overline{s}^*,0)}_{i-\tau-1}+\Big(D_{[i-\tau,i]}-\overline{s}^*\Big)^+.
\end{align*}
From this, we have $\mathbb{E} {B}^{(\overline{s}^*,0)}_\infty+\mathbb{E} {L}^{(\overline{s}^*,0)}_\infty\leq \mathbb{E} {B}^{(\overline{s}^*,0)}_\infty+\mathbb{E} \Big(D_{[i-\tau,i]}-\overline{s}^*\Big)^+$. Hence
\begin{align}
\mathbb{E} {L}^{(\overline{s}^*,0)}_\infty\leq \mathbb{E} \Big(D_{[i-\tau,i]}-\overline{s}^*\Big)^+.\label{12-12-10}
\end{align}
Using (\ref{12-12-8})-(\ref{12-12-10}), from (\ref{12-12-7}), we have $\Big| \bar{\cal C}^{\overline{s}^*}(r, h, r+\tau h; {\bf x})-{\cal C}^{(\overline{s}^*,0)}(r, h; {\bf x})\Big| \le 2(r+(\tau+1)h)
\mathbb{E}\bar{B}^{\overline{s}^*}_\infty$.
(\ref{hanqin-57}) directly follows from $\lim_{r\rightarrow \infty}\mathbb{E}\bar{B}^{\overline{s}^*}_\infty=0$. Thus we have the theorem. \hfill$\Box$

\section{\sf {\bf POP}-System with Unknown Demand Distribution}\label{algorithm}
In the last section,  asymptotic optimality of the $(s,q)$-policy is proved in the setting where the demand distribution is known but  backorders are not observable ({\bf POP}-system).
We now develop an online algorithm for finding the optimal $(s,q)$-policy  when the demand distribution is also not known. Following the usual procedure, we first carry out the concentration analysis for the $(s,q)$-policy. Compared with the existing approach (e.g. \cite{auer2002}, \cite{huhetal2009a}) which mainly focuses on the i.i.d. random variables, here we need a refined geometric ergodicity for the Markov chains to handle the case in which the demand distribution is not fully known and the backorders are not observable.

\subsection{\sf Concentration Analysis of $(s,q)$-Policy}\label{Concentration Analysis of $(S,r)$-Policy_Section}
To analyze for the concentration given by the $(s,q)$-policy, based on Proposition \ref{ergodic base-stock policy_thm}, first we need to specify the geometric convergence rate  for the Markov chain which characterizes the system dynamics. We use the technique developed in Markov chains (\cite{meyntweedie1994}) to identify this rate.

\begin{proposition}\label{ergodic-s-q-demand}
Assume that $\mathbb{P}[D_1=0]=\alpha_0>0$ and $\mathbb{P}[D_1=1]=\alpha_1>0$.
Under the $(s,q)$-policy, the Markov chain generated by the system dynamics of the {\bf POP}-system with $p<1$, $\Big\{X_i^{(s,q)}=\Big(I^{(s,q)}_i,B_{i-1}^{(s,q)},Q_{i-\tau}^{(s,q)},\ldots,Q^{(s,q)}_{i-1}\Big): i\geq 1\Big\}$, is geometrically ergodic  with the steady state $X_\infty^{(s,q)}$.  In particular, for all ${\bf x} \in {\cal X}^{(s,q)}$,
\begin{align}
\sum_{{\bf y}\in {\cal X}^{(s,q)}}\Big|\mathbb{P}\Big[X_i^{(s,q)}={\bf y}\Big|X_1^{(s,q)}={\bf x}\Big]-\mathbb{P}\Big[
 X_\infty^{(s,q)}={\bf y}\Big]\Big|\leq {\cal V}({\bf x})\frac{(1+\beta(s,q))\rho^{i+1}(s,q)}{\rho(s,q)-\theta(s,q)}\label{geo-rate}
 \end{align}
 where ${\cal V}({\bf x})=p\cdot (x_2)^\delta+1$ for ${\bf x}=(x_1,x_2,\ldots,x_{2+\tau})\in {\cal X}^{(s,q)}$, $\beta(s,q)$, $\theta(s,q)$ and $\rho(s,q)$ are three positive constants that only depend on $\alpha_0,\alpha_1$, $p$, and $(s,q)$ with $\theta(s,q)\in (0,1)$ and $\rho(s,q)\in (\theta(s,q), 1)$.
\end{proposition}

For ${\bf x}=(x_1,x_2,x_3,\ldots,x_{2+\tau})\in {\cal X}^{(s,q)}$, define
\begin{align}
&\mathcal{R}({\bf x})=\mathbb{E}\Big[ r\cdot\Big( (x_1+x_3)\wedge (D_1+x_2)\Big)-h\cdot (x_1+x_3-x_2-D_1)^+\Big].\label{dem-19}
\end{align}
 With the help of Proposition {\rm \ref{ergodic-s-q-demand}}, now we give the concentration bound.
\begin{theorem}\label{dem-bound}
Additional assumptions in Proposition {\rm \ref{ergodic-s-q-demand}}, suppose $\mathbb{P}[D_1\leq\Lambda]=1$.
Then for any $\delta,\delta_N\in (0, 1)$, any given $(s,q)$-policy and initial state $X^{(s,q)}_1=(I^{(s,q)}_1,B_0^{(s,q)},Q^{(s,q)}_{1-\tau},\ldots,Q^{(s,q)}_0)\in {\cal X}^{(s,q)}$, with at least probability $(1-\delta_N)$,
\begin{align*}
&\Big| \sum_{i=1}^N \Big({\cal R}(X^{(s,q)}_i) - {\cal C}^{(s,q)}\Big) \Big|\le  \Big (\chi_1(X_1^{(s,q)})+\chi_2N^{\delta}\Big)\Big(1+\sqrt{2(N-1)\ln \frac{2}{\delta_N}}\Big) \ \ \mbox{for $N\geq 2$},\\
&\chi_1(X_1^{(s,q)})= \frac{r(s+q)(1+\lambda_0+b_0)(1+\beta(s,q))\rho^2(s,q)}{(\rho(s,q)-\theta(s,q))(1-\rho(s,q))}\Big(1+p(B_0^{(s,q)})^{\delta}\Big),\\
& \chi_2=\frac{pr(s+q)(1+\lambda_0+b_0)(1+\beta(s,q))\rho^2(s,q)}{(\rho(s,q)-\theta(s,q))(1-\rho(s,q))}\Lambda^{\delta},
\end{align*}
$\rho(s,q)$, $\beta(s,q)$, and $\theta(s,q)$ are given in Proposition {\rm \ref{ergodic-s-q-demand}}, and $\lambda_0$ and $b_0$ are given in {\rm (\ref{12-14-2})}.
\end{theorem}

\noindent
{\it Proof}:  In view of the definition (\ref{dem-19}), we have $\mathcal{R}({\bf x})\leq r(x_1+x_3)\leq r(s+q)$ under the $(s,q)$-policy.
First consider the bias function $h(\cdot)$ on ${\cal X}^{(s,q)}$ given by
\begin{align}
        h({\bf x})=\lim_{N\to\infty}\mathbb{E}\Big[\sum_{i=1}^N \Big(\mathcal{R}(X_i^{(s,q)})-\mathcal{R}(X^{(s,q)}_\infty)\Big)\Big|X^{(s,q)}_1={\bf x}\Big].\label{dem-10}
\end{align}
Then, by Proposition \ref{ergodic-s-q-demand},
\begin{align}
|h({\bf x})|&=\lim_{N\to\infty}\Big|\mathbb{E}\Big[\sum_{i=1}^N \Big(\mathcal{R}(X_i^{(s,q)})-\mathcal{R}(X^{(s,q)}_\infty)\Big)\Big|X^{(s,q)}_1={\bf x}\Big] \Big|\nonumber\\
        &=\lim_{N\to\infty}\Big| \sum_{i=1}^N \sum_{{\bf y}\in {\cal X}^{(s,q)}} \mathcal{R}({\bf y})
        \Big( \mathbb{P}\Big[X_i^{(s,q)}={\bf y}\Big|X^{(s,q)}_1={\bf x}\Big]- \mathbb{P}[X^{(s,q)}_\infty)={\bf y}]\Big)\Big|\nonumber  \\
        & \leq   \lim_{N\to\infty}  \sum_{i=1}^N \sum_{{\bf y}\in {\cal X}^{(s,q)}} \mathcal{R}({\bf y})
        \Big| \mathbb{P}\Big[X_i^{(s,q)}={\bf y}\Big|X^{(s,q)}_1={\bf x}\Big]- \mathbb{P}[X^{(s,q)}_\infty)={\bf y}] \Big|\nonumber  \\
       & \leq  r(s+q) \lim_{N\to\infty}  \sum_{i=1}^N \sum_{{\bf y}\in {\cal X}^{(s,q)}}
        \Big| \mathbb{P}\Big[X_i^{(s,q)}={\bf y}\Big|X^{(s,q)}_1={\bf x}\Big]- \mathbb{P}[X^{(s,q)}_\infty)={\bf y}] \Big|\nonumber  \\
       &\leq  r(s+q) \lim_{N\to\infty}  \sum_{i=1}^N  {\cal V}({\bf x})(1+\beta(s,q))\frac{\Big(\rho(s,q)\Big)^{i+1}}{\rho(s,q)-\theta(s,q)}\nonumber\\
       &=  r(s+q)\frac{(1+\beta(s,q))\rho^2(s,q){\cal V}({\bf x})}{\Big(\rho(s,q)-\theta(s,q)\Big)(1-\rho(s,q))}.\label{dem-11}
    \end{align}
Consider the martingale difference $\Delta_{i+1}=h(X^{(s,q)}_{i+1})-\mathbb{E}\Big(h(X^{(s,q)}_{i+1})\Big|X^{(s,q)}_i\Big)$.
It follows from (\ref{dem-5}) and (\ref{dem-11}) that
\begin{align}
\Big| \mathbb{E}\Big(h(X^{(s,q)}_{i+1})\Big|X^{(s,q)}_i\Big)\Big| &\leq \mathbb{E}\Big( r(s+q)\frac{(1+\beta(s,q))\rho^2(s,q){\cal V}(X^{(s,q)}_{i+1})}{\Big(\rho(s,q)-\theta(s,q)\Big)(1-\rho(s,q))} \Big|X^{(s,q)}_i\Big)\Big| \nonumber\\
&  = r(s+q)\frac{(1+\beta(s,q))\rho^2(s,q)}{\Big(\rho(s,q)-\theta(s,q)\Big)(1-\rho(s,q))} \cdot \mathbb{E}\Big({\cal V}(X^{(s,q)}_{i+1})\Big|
X^{(s,q)}_i\Big)\nonumber\\
&  \leq r(s+q)(\lambda_0+b_0)\frac{(1+\beta(s,q))\rho^2(s,q)}{\Big(\rho(s,q)-\theta(s,q)\Big)(1-\rho(s,q))} \cdot {\cal V}(X^{(s,q)}_{i}). \label{dem-12}
\end{align}
As the demand is capped by $\Lambda$, with probability one, $B_i^{(s,q)}\leq B^{(s,q)}_0+i \Lambda$ holds. This gives us that with probability one,
\begin{align}
{\cal V}(X^{(s,q)}_i)\leq p\cdot \Big(B^{(s,q)}_0+(i-1) \Lambda\Big)^\delta+1,\label{dem-14}
\end{align} which, by (\ref{dem-11})-(\ref{dem-12}), implies that
\begin{align}
\Delta_{i+1} \leq \chi_1 (B_0^{(s,q)}) +\chi_2 i^\delta   \ \ \mbox{for $i\geq 1$}.\label{dem-13}
\end{align}
By the Hoeffding-Azuma inequality (Theorem 3 on p476 in \cite{grimmettstirzaker2001}), and (\ref{dem-13}), for any $\delta_N\in (0, 1)$,
\begin{align}
\mathbb{P}\Big[ \Big| \sum_{i=1}^{N-1} \Delta_{i+1}\Big| \geq (\chi_1(B_0^{(s,q)})+\chi_2N^\delta)\sqrt{2(N-1)\ln \frac{2}{\delta_N}}\Big]
\le \delta_N. \label{dem-16}
\end{align}
By Poisson's equation, ${\cal R}(X_i^{(s,q)})-{\cal C}^{(s,q)}=h(X^{(s,q)}_{i})-\mathbb{E}\Big(h(X^{(s,q)}_{i+1})\Big|X^{(s,q)}_i\Big) $, we have
 \begin{align}
 \Big|\sum_{i=1}^N \Big({\cal R}(X^{(s,q)}_i) - {\cal C}^{(s,q)}\Big) \Big|
 &= \Big|\sum_{i=1}^N\Big(h(X^{(s,q)}_i)-\mathbb{E}(h(X^{(s,q)}_{i+1})|X^{(s,q)}_i)\Big)\Big| \nonumber\\
  & =\Big|\Big[h(X_1^{(s,q)})-\mathbb{E}\Big(h(X_{N+1}^{(s,q)})|X_N^{(s,q)}\Big)\Big]\Big| +\Big|\sum_{i=1}^{N-1} \Delta_{i+1}\Big|.\label{dem-15}
  \end{align}
Using again (\ref{dem-11})-(\ref{dem-14}), similar to (\ref{dem-13}),
\begin{align}
&|h(X_1^{(s,q)})|+\Big|\mathbb{E}\Big(h(X_{N+1}^{(s,q)})|X_N^{(s,q)}\Big)\Big|\nonumber\\
& \ \ \ \le r(s+q)\frac{(1+\beta(s,q))\rho^2(s,q)}{\Big(\rho(s,q)-\theta(s,q)\Big)(1-\rho(s,q))} \Big[{\cal V}( X_1^{(s,q)})+(\lambda_0+b_0){\cal V}( X_N^{(s,q)})\Big]\nonumber\\
& \ \ \ \leq \chi_1 (B^{(s,q)}_0)+\chi_2 N^\delta.\label{dem-17}
\end{align}
If follows from (\ref{dem-15})-(\ref{dem-17}) that when $\Big| \sum_{i=1}^{N-1} \Delta_{i+1}\Big| < \Big(\chi_1(B^{(s,q)}_0)+\chi_2N^\delta\Big)\sqrt{2(N-1)\ln \frac{2}{\delta_N}}$,
\begin{align*}
\Big|\sum_{i=1}^N \Big({\cal R}(X^{(s,q)}_i) - {\cal C}^{(s,q)}\Big) \Big|& \leq \Big(\chi_1(B^{(s,q)}_0) +\chi_2 N^\delta\Big) +\Big (\chi_1(B^{(s,q)}_0)+\chi_2N^\delta\Big)\sqrt{2(N-1)\ln \frac{2}{\delta_N}}\\
&\le  \Big (\chi_1(B^{(s,q)}_0)+\chi_2N^\delta\Big)\Big(1+\sqrt{2(N-1)\ln \frac{2}{\delta_N}}\Big).
\end{align*}
This, by (\ref{dem-16}), shows the theorem holds.\hfill$\Box$

Here we separately put two terms $\chi_1(\cdot)$ and $\chi_2 N^\delta$ as we try to clearly show the penalty of the degeneracy from the uniform ergodic Markov chain to the Markov chain only being geometrically ergodic, where the second term $\chi_2 N^\delta$ can be taken as such penalty. This degeneracy is mainly from the dimension of the backorder part. When $p=0$, the system is nothing but the classic lost-sales system, for which, we can follow the idea in \cite{huhetal2009a} showing that the Markov chain under each $(s,q)$-policy is uniformly ergodic. By Theorem 16.0.2 in \cite{meyntweedie2012}, this is equivalent to Proposition \ref{ergodic-s-q-demand}  holds for ${\cal V}(\cdot)=1$ and thus we can choose some $\chi_1$ to bound the martingale difference $\Delta_i$ uniformly and also make the above bound holds. However, here as we can see that the chain we consider here may not always be uniformly ergodic but only the geometrically ergodicity can be guaranteed. This degeneracy will then introduce some loss of convergence rate to the steady state of the Markov chain, as shown by the the term $\chi_2 N^\delta$.

\subsection{\sf Algorithm and Regret Analysis}\label{Regrest Analysis_SubSection}

Since demand is capped by $\Lambda$, we know that the optimal $(s^*,q^*)$-policy  satisfies  $s^*+q^*\leq \tau \Lambda$. Hence, we can restrict our policy searching to the bounded region: $[0,\underline{s}]\times [0, \underline{q}]$.
The algorithm will split the time horizon into several epochs, and in each epoch, we will choose one $(s,q)$-policy with the largest UCB-index to implement during the whole epoch. At the end of each epoch, we update the UCB-index for each $(s,q)$-policy and again choose the policy with the largest UCB-index to implement for the coming epoch.
Here, the construction of the UCB-index is not a trivial replication of the one for multi-arm bandits problem where the UCB-index is a direct result from the concentration of the i.i.d. rewards. Our UCB-index is constructed via the theory of the ergodicity of the Markov chains, which mainly applies the Lyapunov-type analysis to quantify the mixing time of the Markov chains.

\begin{algorithm}
\caption{Online UCB-index Algorithm for Partial Backorder Inventory}\label{online algorithm}
\KwData{The upper and lower bound: $[0,\underline{s}]$, $[0,\underline{q}]$; initial inventory $I_{0}=\underline{s}+\underline{q}$; $\varphi_0(s,q)=1$ and $T_0(s,q)=1$  for all $(s,q)\in[0,\underline{s}]\times[0,\underline{q}]$; $\eta_0 = 0$; $H_0=0$, $\delta_\ell=\frac{1}{2^{2\ell}}$;}
Initialization: At period-$0$, observe the sales $\min(I_0,D_0)$, then for all $(s,q)$, calculate
\begin{align*}
\breve{\cal R}_0^{(s,q)}=\widetilde{\cal R}_0^{(s,q)}&= r\cdot \min\Big(s+q, \min(I_0,D_0)\Big)-h\Big(s+q-\min(I_0,D_0)\Big)^+;\\
    I^{(s,q)}_1&=(s+q-D_0)^+;
\end{align*}
\For{$1 \le \ell \le T_{\cal E}(N) = \min\{k: \sum_{i=1}^k {\cal E}_i \le N \}$}{For all $(s,q)$-policy, update ${\cal N}_\ell(s,q)$ and their UCB-index by
\begin{align*}
       & {\cal N}_\ell(s,q):=\Big\{(s',q'):|s'-s|\le\frac{\varepsilon_1}{1\vee \sqrt{\eta_{\ell-1}}},|q'-q|\le\frac{\varepsilon_2}{1\vee \sqrt{\eta_{\ell-1}}}, T_{\ell-1}(s',q')\ge T_{\ell-1}(s,q) \Big\}\\
&F_\ell^{(s,q)}= \frac{1}{|{\cal N}_\ell(s,q)|}\sum_{(s',q')\in {\cal N}_\ell(s,q)}\Big\{ \widetilde{\cal R}^{(s',q')}_{\ell-1}+\frac{H_{\ell-1}}{T_{\ell-1}(s',q')}\Big(\sqrt{2(T_{\ell-1}(s',q')-1)\ln\frac{2}{\delta_\ell}}+1\Big)\Big\}
    \end{align*}

Choose $(s^*_\ell,q^*_\ell)\in\arg\max\limits_{(s,q)}\{F^{(s,q)}_\ell\}$ (break ties arbitrarily)\;
Update $\varphi_\ell(s^*_\ell,q^*_\ell)= \varphi_{\ell-1}(s^*_\ell,q^*_\ell)+ 1$ and $\varphi_\ell(s,q)= \varphi_{\ell-1}(s,q)$ for all $(s,q)\neq(s^*_\ell,q^*_\ell)$\;
Update ${\cal E}_\ell = 2^{\varphi_\ell(s^*_\ell,q^*_\ell)}$\;
Implement the $(s^*_\ell,q^*_\ell)$-policy from $\eta_{\ell-1}+1$ to $\eta_\ell:= \min\{\eta_{\ell-1} + {\cal E}_\ell ,N\}$ to generate
$\breve{\cal R}_k^{(s^*_\ell,q^*_\ell)}$\;
\If{$\{\eta_{\ell-1}+1 \le k \le \eta_\ell: I_k+q_{[k-\tau,k]} \le s^*_\ell+q^*_\ell \}\neq\emptyset$}{
Update $\nu_\ell = \min\{k: I_k+q_{[k-\tau,k]} \le s^*_\ell+q^*_\ell \ \mbox{and} \ \eta_{\ell-1}+1 \le k \le \eta_\ell\}$\;
Update $T_\ell(s^*_\ell,q^*_\ell)=(\eta_{\ell}-\nu_\ell+1)$ and $T_\ell(s,q) = T_{\ell-1}(s,q)$ for all $(s,q)\neq (s^*_\ell,q^*_\ell)$\;
Update $\widetilde{\cal R}_\ell^{(s^*_\ell,q^*_\ell)}=\frac{1}{T_\ell(s^*_\ell,q^*_\ell)}\sum_{k=\nu_\ell}^{\eta_\ell}\breve{\cal R}^{(s^*_\ell,q^*_\ell)}_k$, and
$\widetilde{\cal R}_\ell^{(s,q)}=\widetilde{\cal R}_{\ell-1}^{(s_,q)}$  for all $(s,q)\neq (s^*_\ell,q^*_\ell)$;}
\Else{Update   $T_\ell(s,q) = T_{\ell-1}(s,q)$ and $ \widetilde{\cal R}_\ell^{(s,q)}=\widetilde{\cal R}_{\ell-1}^{(s_,q)}$ for all $(s,q)$;}
Update $H_\ell=\chi_1+\chi_2\cdot \eta_\ell^\delta$, and $\ell\to \ell+1$.}

\end{algorithm}

According to the online UCB-index algorithm developed in the following, the number of the policies to be explored is bounded by $\underline{s}\underline{q}$,  and each epoch takes the number of the periods to explore according to the power-2 growth. Let $L_{\cal E}(N)$ represents the total number of the epochs by $N$ periods, among $L_{\cal E}(N)$ epochs, let $L_{\cal E}^{(s,q)}(N)$ be the number of the epochs we use the $(s,q)$-policy to explore. Then $\sum_{(s,q)} \sum_{\ell=1}^{L_{\cal E}^{(s,q)}(N)} 2^\ell \leq N$. This implies that $L_{{\cal E}}(N)\leq \underline{s}\underline{q}\log_2N$.
In each epoch, say epoch-$\ell$, let $(s^*_\ell,q^*_\ell)$ represent the policy with the largest UCB-index, and $\breve{{\cal R}}_k^{(s^*_\ell,q^*_\ell)}$ represent the reward we observe (sample) at period-$k$ under $(s^*_\ell,q^*_\ell)$-policy. At the beginning of epoch-$\ell$, we may need to take few periods ({\it switching periods} from period-$(\eta_{\ell-1}+1)$ to period-$(\nu_\ell-1)$) to make that the system state becomes valid, that is, $I_k+q_{[k-\tau,k-1]}\leq s^*_\ell+q^*_\ell$. As the probability that demand is at least one is $1-\mathbb{P}[D_1=0]=1-\alpha_0$, and starting with any state in which the on-hand inventory plus the inventory on-order is below $\underline{s}+\underline{q}$, the expected number of the periods needed is then bounded by $\frac{1-(1-\alpha_0)^{\underline{s}+\underline{q}}}{\alpha_0 (1-\alpha_0)^{\underline{s}+\underline{q}}}$. This further gives us the expected regret incurred by the switching periods:
\begin{align}
\mathbb{E}\Big[\sum_{k={\eta_{\ell-1}+1}}^{\nu_\ell-1} \Big({\cal C}^{(s^*,q^*)}-\breve{\cal R}_k^{(s^*_\ell,q^*_\ell)}\Big)\Big] \leq r (\underline{s}+\underline{q}) \frac{1-(1-\alpha_0)^{\underline{s}+\underline{q}}}{\alpha_0 (1-\alpha_0)^{\underline{s}+\underline{q}}},\label{hanqin-80}
\end{align}
where $(s^*,q^*)$ is the optimal one among the whole $(s,q)$-policies.
Now look at the concentration generated by the empirical rewards in the algorithm at epoch-$\ell$ for each $(s,q)$-policy with $(s,q)\in [0,\underline{s}]\times[0,\underline{q}]$, namely,
\begin{align*}
C_\ell = \Big\{\Big|\widetilde{\cal R}_{\ell}^{(s,q)}-{\cal C}^{(s,q)}\Big|\le \frac{H_\ell}{T_\ell(s,q)}\Big(\sqrt{2(T_\ell(s,q)-1)\ln \frac{2}{\delta_\ell}}+1\Big) \ \mbox{for all \ } 0\le s\le \underline{s} \ \mbox{and} \ 0\leq q\leq \underline{q}\Big\}.
\end{align*}

In the algorithm, $T_\ell(s,q)$ is random but less than $2^\ell$, and $H_\ell=\chi_1+\chi_2 \eta_\ell^\delta \geq \chi_1+\chi_2 (T_\ell(s,q))^\delta$, thus by Theorem \ref{dem-bound},
for each $(s,q)$,
\begin{align*}
&\mathbb{P}\Big[\Big|\widetilde{\cal R}_{\ell}^{(s,q)}-{\cal C}^{(s,q)}\Big|> \frac{H_\ell}{T_\ell(s,q)}\Big(\sqrt{2(T_\ell(s,q)-1)\ln \frac{2}{\delta_\ell}}+1\Big)\Big]\\
& \ \ \ \le \sum_{n=1}^{2^\ell} \mathbb{P}\Big[\Big|\widetilde{\cal R}_{\ell}^{(s,q)}-{\cal C}^{(s,q)}\Big|> \frac{H_\ell}{n}\Big(\sqrt{2(n-1)\ln \frac{2}{\delta_\ell}}+1\Big)\Big]
\le 2^\ell \times \delta_\ell=\frac{1}{2^\ell}.
\end{align*}
Hence, $\mathbb{P}[C_\ell]\geq 1-  \underline{s}\underline{q}/2^\ell$. The total expected regret on the complement of $C_\ell$ (denoted by $C_{\ell,c}$),
\begin{align}
\mathbb{E}\Big[\sum_{\ell=1}^{T_{\cal E}(N)} \sum_{k=\nu_\ell}^{\eta_\ell} \Big({\cal C}^{(s^*,q^*)}-  \breve{\cal R}_k^{(s^*_\ell,q^*_\ell)}\Big)\mathbb{I}_{\{C_{\ell,c}\}}\Big]
\leq r(\underline{s}+\underline{q})T_{\cal E}(N) \mathbb{P}[C_{\ell,c}]\times 2^\ell \leq r \underline{s}\underline{q}(\underline{s}+\underline{q})\log_2N.\label{hanqin-81}
\end{align}
With the help of Theorem \ref{dem-bound}, we can analyze the regret on set $C_{\ell}$. Then we have the following theorem about the regret bound from the above algorithm. Its proof is relegated to Appendix.
\begin{theorem}\label{regretbound}
Assume that $\mathbb{P}[D_1\leq\Lambda]=1$, $\mathbb{P}[D_1=0]>0$ and $\mathbb{P}[D_1=1]>0$. Then there exists a constant $c\le r(\underline{s}+\underline{q})$ such that
\begin{align*}
{\sf Rg}(N)=&\sum_{\ell=1}^{L_{\cal E}(N)}\sum_{k=\eta_{\ell-1}+1}^{\eta_\ell} \Big({\cal C}^{(s^*,q^*)}- \breve{\cal R}_k^{(s^*_\ell,q^*_\ell)}\Big)
\le 5(2+\sqrt2)\underline{s}\underline{q}\Big(\underline{\chi}_1+\underline{\chi}_2 N^\delta\Big) \sqrt{\frac{1}{2}N\ln2+\ln N}\\
& \ \ \ \ \ \ \ \ \ \ \ \ \ \ \ \ \ \ \ \ \ \ \ \ \ \ \ \ \ \ \ \ \ \ \ \ \
\ \ \ \ \ \  +2(\sqrt 2+1)(\varepsilon_1+\varepsilon_2)c\underline{s}\underline{q}\sqrt{N}+ o(\sqrt{N}),\\
&\underline{\chi}_1= \max_{(s,q)\in [0,\underline{s}]\times[0,\underline{q}]}\frac{r(s+q)(1+\lambda_0+b_0)(1+\beta(s,q))\rho^2(s,q)}{(\rho(s,q)-\theta(s,q))(1-\rho(s,q))},\\
&\underline{\chi}_2=\max_{(s,q)\in [0,\underline{s}]\times[0,\underline{q}]}\frac{pr(s+q)(1+\lambda_0+b_0)(1+\beta(s,q))\rho^2(s,q)}{(\rho(s,q)-\theta(s,q))(1-\rho(s,q))}\Lambda^\delta.
\end{align*}
That is, the regret under the algorithm in the first $N$ periods ${\sf Rg}(N)=O(N^{\frac{1}{2}+\delta})$ with arbitrary small $\delta>0$.
\end{theorem}

\section{\sf Numerical Experiments}\label{num}
In this section, we conduct numerical experiments for the {\bf P}-system to show the performance of the base-stock policy developed in Section \ref{policy}, the $(s,q)$-policy  studied in Section \ref{partial-ob}, and the online algorithm proposed in Section \ref{algorithm}. In all numerical experiments, we set $h=1$, and consider that demands follow either Poisson or Binomial distribution, and the patience of the unmet demands has two scenarios: high with $p=0.7$, and low with $p=0.3$. In order to make the comparison with the existing insights developed for the {\bf L}-systems (\cite{huhetal2009b}, and  \cite{xingoldberg2016}), we let the lead time $\tau$ change from the shorter one ($\tau=2)$ to the longer one $(\tau=10)$, and the sale price change from the lower one $(r=4)$ to the higher one $(r=64)$.

\subsection{\sf Performance of Base-Stock Policies in {\bf P}-System}\label{num-base}

Now we look at the performance given by the best base-stock policy for the {\bf P}-system. From Theorem \ref{symptotic-s-policy}, we know its asymptotic optimality when the sale price $r$ gets large. To measure its performance, the optimality gap is defined by
 \[
 \Delta^*_{BS}:=1-\frac{\max_s {\cal C}^s(r,h)}{\max_{\pi\in {\cal A}}{\cal C}^\pi(r,h)}.
 \]
 For the {\bf P}-system, we do not know its optimal policy, and even the global search is too time consumable to numerically get the optimal policy. Thus, instead of directly computing $\Delta^*_{BS}$, we estimate the optimality gap $\Delta^*_{BS}$. By Propositions \ref{upper bound of base-stock P-system_prop} and \ref{lower bound of optimal P-system_prop}, for any initial state ${\bf x}$,
\begin{align*}
\bar{\cal C}^{\overline{s}^*} (r,h,r+\tau h;{\bf x})-\limsup_{N\to\infty}\frac{1}{N}\mathbb{E}h\cdot O^{\overline{s}^*}_{[1,N]}\le {\cal C}^{\overline{s}^*} (r,h;{\bf x}) \le \max_{\pi\in {\cal A}}{\cal C}^\pi(r,h;{\bf x})\le \max_s \overline{\cal C}^s \Big(r,h,\frac{(1-p)r}{2(1+\tau)};{\bf x}\Big),
\end{align*}
where $\overline{s}^* = \argmax_{s}\overline{\cal C}^{s} (r,h,r+\tau h;{\bf x})$. Then we obtain an upper bound of $\Delta^*_{BS}$, which is denoted by $\overline{\Delta}$, and is  referred to as the upper optimality gap:
\begin{align*}
    \Delta^*_{BS} 
\le 1 - \frac{\bar{\cal C}^{\overline{s}^*} (r,h,r+\tau h)-\limsup_{N\to\infty}\frac{1}{N}\mathbb{E}h\cdot O^{\overline{s}^*}_{[1,N]}}{\max_s \overline{\cal C}^s (r,h,\frac{(1-p)r}{2(1+\tau)})}=:\overline{\Delta}.
\end{align*}
Let $s^*=\arg\max_{s}{\cal C}^s(r,h:{\bf x})$. Tables \ref{tab 1:avg_profit_and_levels_basestock_policy_zhaoxuan_0.7} and \ref{tab 3:avg_profit_and_levels_basestock_policy_zhaoxuan_0.3}
show that the long-run average profits under the best base-stock policy for the {\bf P}-system and $\overline{s}^*$-policy determined by the optimal base-stock policy from the {\bf B}-system. At the same time, the two tables also provide the concrete values of $(s^*,\overline{s}^*)$. We observe that the profit difference between ${\cal C}^{s^*}(r,h; {\bf x})$ and ${\cal C}^{\overline{s}^*}(r,h;{\bf x})$ gets closer as the sale price get larger. A similar observation for the policy difference between $s^*$ and $\overline{s}^*$  holds.
These two observations indicate that the base-stock policy is nearly optimal when the sales price $r$ dominates the holding cost $h$ and the optimal base-stock policy for the {\bf B}-system would be a good proxy of the optimal policy in this regime, which is consistent with our theoretical results in Theorem \ref{symptotic-s-policy}.

As shown above, it might be a good choice for managers to choose the base-stock level $\overline{s}^*$ given by {\bf B}-system when the sale price is high since the long-run average profit under policy $\overline{s}^*$ will converge to the optimal profit. However, it doesn't mean we can claim that considering the patience of the customers should be negligible in the limit regime. First, in the non-asymptotic regime, the $\overline{s}^*$ policy could perform poorly. Second, even in asymptotic regime, the decisions under the $\overline{s}^*$ could still be different from the best base-stock policy $s^*$. Specifically, in non-asymptotic regime, the best order-up-to level $s^*$ will shrink about $20\%-25\%$, and in asymptotic regime, where $h/r=1/64\approx 1.6\%$, the relative difference could still be around $7\%-10\%$. This shows the significant value of considering the customer abandonment.

\begin{table}[h]
    \centering
    \resizebox{1.0\textwidth}{!}{
    \begin{tabular}{cccccccccccccc}
        \toprule[1pt]
              &\multirow{2}{*}{Demand}&\multicolumn{5}{c}{\text{Average Profit $({\cal C}^{\overline{s}^*}(r,h), {\cal C}^{s^*}(r,h))$}} &\multicolumn{5}{c}{\text{ $(s^*, \overline{s}^*)$}}\\
              & & $r=4$ & $r=8$ & $r=16$ & $r=32$ & $r=64$ & $r=4$ & $r=8$ & $r=16$ & $r=32$ & $r=64$ \\ \cmidrule[0.5pt](l{.25em}r{.25em}){2-2}
              \cmidrule[0.5pt](l{.25em}r{.25em}){3-7}      \cmidrule[0.5pt](l{.25em}r{.25em}){8-12}

         \multirow{2}{*}{$\tau=2$} & Poisson(10)& (33.06, 35.84) & (72.07, 74.19) & (150.04, 152.50) & (308.87, 310.83) & (626.89, 629.30) & (29,36) & (32,37) & (34,39)& (37,41) & (39,42)\\

         &Binomial(10,0.5)& (16.56, 17.88)& (35.71, 37.11)& (75.42, 76.28)& (154.41, 155.47)& (313.81, 314.83)&  (15,18) & (16,19)& (17,19)& (18,20)& (19,21)\\&&&&&&&&&&&\\

         \multirow{2}{*}{$\tau=4$}&Poisson(10)& (30.06, 35.05) & (68.91, 73.03) & (147.05, 150.74) & (304.91, 308.48) & (623.55, 626.46) &  (47,59) & (50,60) & (55,62) & (58,64) & (61,66)\\

         &Bino(10,0.5)&(15.45, 17.48) & (34.74, 36.48) & (73.62, 75.43) & (152.84, 154.34) & (312.07, 313.50)&(24,29)& (26,30)& (27,31)& (29,32)& (30,33)\\&&&&&&&&&&&\\

         \multirow{2}{*}{$\tau=6$}&Poisson(10)& (28.02, 34.58) & (66.40, 72.19) & (144.35, 149.54) & (302.35, 306.81) & (620.70, 624.49) &  (65,81) & (70,83) & (74,85) & (78,87) & (82,89)\\

         &Bino(10,0.5)& (13.64, 17.23)& (33.55, 36.07)& (72.27, 74.79)& (151.64, 153.55)& (310.99, 312.40)&  (33,41)& (35,41)& (37,42)& (39,43)& (41,44)\\&&&&&&&&&&&\\

        \multirow{2}{*}{$\tau=8$}&Poisson(10)& (25.22, 34.22) & (64.20, 71.62) & (142.16, 148.60) & (300.49, 305.46) & (617.01, 622.80) &  (82,104) & (89,105) & (94,107) & (99,109) & (102,112)\\

        &Bino(10,0.5)& (12.63, 17.01)& (32.56, 35.74)& (71.53, 74.28)& (150.55, 152.80)& (310.09, 311.51) &  (41,52)& (44,52)& (47,53)& (49,54)& (51,55) \\&&&&&&&&&&&\\

       \multirow{2}{*}{$\tau=10$}&Poisson(10)& (23.44, 33.97) & (62.06, 71.17)& (140.11, 147.89)& (298.34, 304.45)& (614.33, 621.41)& (101,126)& (107,127)& (113,129)& (118,131)& (123,134)\\

       &Bino(10,0.5) & (11.72, 16.87)& (31.60, 35.53)& (70.59, 73.87)& (149.54, 152.23)& (308.40, 310.83) & (51,63)& (54,63)& (57,64)& (59,65)& (62,67)\\ \bottomrule[1pt]
    \end{tabular}}
    \caption{Average Profits and Base-Stock Levels of Base-Stock Policies ($h=1, p = 0.7$)}
    \label{tab 1:avg_profit_and_levels_basestock_policy_zhaoxuan_0.7}
\end{table}
\begin{table}[h]
    \centering
    \resizebox{0.75\textwidth}{!}{
    \begin{tabular}{cccccccccccccc}
        \toprule
              &\multirow{2}{*}{Demand} &\multicolumn{5}{c}{\text{Upper Optimality Gap $\overline{\Delta}$}}&\multicolumn{5}{c}{\text{$1-{\cal C}^{\overline{s}^*}(r,h;{\bf x})/{\cal C}^{s^*}(r,h;{\bf x})$}}\\
              & & $r=4$ & $r=8$ & $r=16$ & $r=32$ & $r=64$ & $r=4$ & $r=8$ & $r=16$ & $r=32$ & $r=64$ \\ \cmidrule[0.5pt](l{.25em}r{.25em}){2-2}
              \cmidrule[0.5pt](l{.25em}r{.25em}){3-7}      \cmidrule[0.5pt](l{.25em}r{.25em}){8-12}

         \multirow{2}{*}{$\tau=2$} & Poisson(10)&   19.12&  9.98& 5.00   &2.42& 1.19 & 7.49&  2.66& 1.29   & 0.61 & 0.20 \\

         &Bino(10,0.5)&  18.18&  9.29& 4.61 & 2.18& 1.02&  6.97&  3.52& 0.79 & 0.42& 0.25\\&&&&&&&&&&&\\

         \multirow{2}{*}{$\tau=4$}&Poisson(10) &  28.44&  14.47&  7.42&  3.68& 1.77&   14.04&  5.21&  2.22& 0.91 & 0.40\\

         &Bino(10,0.5)& 27.37&  13.65& 6.91&  3.41 & 1.62&  11.21&  5.02& 2.18& 0.95 & 0.41\\&&&&&&&&&&&\\

         \multirow{2}{*}{$\tau=6$}&Poisson(10)  &  36.60&  18.56&  9.59&  4.68& 2.26&  18.59&  8.11&  3.36&  1.40& 0.56\\

         &Bino(10,0.5)& 35.61&  18.05&  9.13&  4.47& 2.07&  20.22&  8.11&  3.36&  1.40& 0.56\\&&&&&&&&&&&\\

        \multirow{2}{*}{$\tau=8$}&Poisson(10)&  44.16&  22.03&  11.37&  5.50& 2.80&  26.70 &  10.32&  4.28&  1.62 & 0.75\\

        &Bino(10,0.5)& 40.75&  21.51&  10.72&  5.37& 2.52&  26.06&  8.82&  3.59& 1.38& 0.50\\&&&&&&&&&&&\\

       \multirow{2}{*}{$\tau=10$}&Poisson(10)& 50.20&  25.60&  12.98&  6.44& 3.18&  31.54&  12.67&  5.14&  1.94& 0.92\\

       &Bino(10,0.5) & 49.16& 24.44& 12.52& 6.20& 2.98&  31.11& 11.06& 4.38& 1.74 & 0.87\\ \bottomrule
    \end{tabular}}
    \caption{Optimality Gaps of Base-Stock Policies ($h=1, p = 0.7$)}
    \label{tab 2:optimality_gap_basestock_policy_zhaoxuan_0.7}
\end{table}


\begin{table}
    \centering
    \resizebox{1.0\textwidth}{!}{
    \begin{tabular}{cccccccccccccc}
        \toprule
              &\multirow{2}{*}{Demand}&\multicolumn{5}{c}{\text{Average Profit $({\cal C}^{\overline{s}^*}(r,h), {\cal C}^{s^*}(r,h))$}} &\multicolumn{5}{c}{\text{ $(s^*, \overline{s}^*)$}}\\
              & & $r=4$ & $r=8$ & $r=16$ & $r=32$ & $r=64$ & $r=4$ & $r=8$ & $r=16$ & $r=32$ & $r=64$ \\ \cmidrule[0.5pt](l{.25em}r{.25em}){2-2}
              \cmidrule[0.5pt](l{.25em}r{.25em}){3-7}      \cmidrule[0.5pt](l{.25em}r{.25em}){8-12}

         \multirow{2}{*}{$\tau=2$} & Poisson(10)& (32.72, 34.07) & (71.35, 72.14) & (149.51, 150.18) & (307.95, 308.20) & (625.68, 626.16) &  (31, 36)&  (34, 37) & (36, 39) & (39, 41) & (41, 42) \\

         &Bino(10,0.5)& (16.38, 17.04)& (35.44, 36.17)& (75.08, 75.27)& (154.29, 154.37)& (313.31, 313.51)&  (16, 18)&  (17, 19)& (18, 19) & (19, 20) & (20, 21)\\&&&&&&&&&&&\\

         \multirow{2}{*}{$\tau=4$}&Poisson(10)& (29.74, 33.29) & (68.46, 70.88) & (146.90, 148.32) & (304.44, 305.76) & (625.68, 623.17) &  (49, 59)&  (53, 60)&  (57, 62)& (60, 64)& (63, 66)\\

         &Bino(10,0.5)&(15.36, 16.60)& (34.32, 35.46) & (73.57, 74.27) & (152.27, 153.13) & (311.91, 311.97) &  (25, 29)& (27, 30) & (28, 31) & (30, 32) & (31, 33)\\&&&&&&&&&&&\\

         \multirow{2}{*}{$\tau=6$}&Poisson(10)& (27.48, 32.80) & (65.71, 70.06) & (143.66, 147.07) & (301.93, 304.01) & (619.01, 620.90) &  (67, 81)&  (72, 83)& (77, 85)& (80,87)& (84, 89)\\

         &Bino(10,0.5)& (13.50, 16.34) & (33.34, 35.03) & (72.32, 73.61)& (150.86, 152.23) & (310.37, 310.88) & (33, 41)&  (36, 41) & (38,42) & (40, 43) & (42, 44) \\&&&&&&&&&&&\\

        \multirow{2}{*}{$\tau=8$}&Poisson(10)& (25.04, 32.49) & (63.82, 69.47) & (141.77, 146.14) & (299.53, 302.64) & (617.63, 620.50) &  (84, 104)&  (90, 105)&  (96, 107)& (101, 109)& (105, 112)\\

        &Bino(10,0.5)& (12.51, 16.17) & (32.34, 34.71) & (71.41, 73.13)& (150.26, 151.50) & (309.16, 309.93) &  (42, 52)&  (45, 52) &  (48, 53) &  (51, 54) & (52, 55)\\&&&&&&&&&&&\\

       \multirow{2}{*}{$\tau=10$}&Poisson(10)& (22.87, 32.24)& (61.82, 69.02) & (139.97, 145.39) & (297.52, 301.66)& (614.73, 617.83) & (101,126)&  (109, 127)&  (116, 129)& (121, 131) & (126, 134)\\

       &Bino(10,0.5) & (11.57, 16.04) & (31.45, 34.46) & (70.43, 72.74) & (149.40, 150.95) & (307.50, 309.17) & (51, 63) & (55, 63) & (58, 64) & (61, 65) & (63, 67)\\ \bottomrule
    \end{tabular}}
    \caption{Average Profits and Base-Stock Levels of Base-Stock Policies ($h=1, p = 0.3$)}
    \label{tab 3:avg_profit_and_levels_basestock_policy_zhaoxuan_0.3}
\end{table}

\begin{table}
    \centering
    \resizebox{0.75\textwidth}{!}{
    \begin{tabular}{cccccccccccccc}
        \toprule
              &\multirow{2}{*}{Demand} &\multicolumn{5}{c}{\text{Upper Optimality Gap $\overline{\Delta}$}}&\multicolumn{5}{c}{\text{$1-{\cal C}^{\overline{s}^*}(r,h;{\bf x})/{\cal C}^{s^*}(r,h;{\bf x})$}}\\
              & & $r=4$ & $r=8$ & $r=16$ & $r=32$ & $r=64$ & $r=4$ & $r=8$ & $r=16$ & $r=32$ & $r=64$ \\ \cmidrule[0.5pt](l{.25em}r{.25em}){2-2}
              \cmidrule[0.5pt](l{.25em}r{.25em}){3-7}      \cmidrule[0.5pt](l{.25em}r{.25em}){8-12}

         \multirow{2}{*}{$\tau=2$} & Poisson(10)&  16.30&  7.92& 3.82 & 1.77 & 0.83& 4.49 & 1.14 & 0.47 & 0.22 & 0.04 \\

         &Bino(10,0.5)&  15.26&  7.42& 3.49 & 1.55& 0.70 & 4.10 & 1.92 & 0.21 & 0.09 & 0.08\\&&&&&&&&&&&\\

         \multirow{2}{*}{$\tau=4$}&Poisson(10) &  26.27&  12.72&  5.97&  2.82& 1.31& 10.90 & 3.37 & 1.23 & 0.45 & 0.17\\

         &Bino(10,0.5)&  24.71&  11.83& 5.58& 2.58 & 1.21 & 8.25 & 3.16 & 1.21 & 0.50 & 0.22\\&&&&&&&&&&&\\

         \multirow{2}{*}{$\tau=6$}&Poisson(10) &  34.78&  16.87& 8.07&  3.76& 1.83 & 16.00 & 6.26 & 2.33 & 0.83 & 0.30\\

         &Bino(10,0.5)&  33.43&  16.04&  7.65&  3.63& 1.67 & 17.60 & 4.97 & 1.77 & 0.64 & 0.25\\&&&&&&&&&&&\\

        \multirow{2}{*}{$\tau=8$}&Poisson(10)&  42.31&  20.44&  9.82&  4.75& 2.20 & 23.87 & 8.41 & 3.11 & 1.07 & 0.47\\

        &Bino(10,0.5)&  40.64&  19.95&  9.37&  4.37& 2.09& 23.07 & 7.04 & 2.55 & 0.86 & 0.30\\&&&&&&&&&&&\\

       \multirow{2}{*}{$\tau=10$}&Poisson(10) &49.80&  23.84&  11.41&  5.63 & 2.56 & 29.47 & 10.69 & 3.98 & 1.42 & 0.58\\

       &Bino(10,0.5) &47.26& 23.27& 10.79& 5.29& 2.46& 28.76 & 9.26& 3.38& 1.16 &0.62\\ \bottomrule
    \end{tabular}}
    \caption{Optimality Gaps of Base-Stock Policies ($h=1, p = 0.3$)}
    \label{tab 4:optimality_gaps_basestock_policy_zhaoxuan_0.3}
\end{table}

Tables \ref{tab 2:optimality_gap_basestock_policy_zhaoxuan_0.7} and \ref{tab 4:optimality_gaps_basestock_policy_zhaoxuan_0.3}
provide the upper optimality gaps and the relative difference between ${\cal C}^{\overline{s}^*}(r,h;{\bf x})$ and ${\cal C}^{s^*}(r,h;{\bf x})$, given by $1-{\cal C}^{\overline{s}^*}(r,h;{\bf x})/{\cal C}^{s^*}(r,h;{\bf x})$.
 The two tables show that the optimality gap  will shrink as the sale price $r$ grows. Further, consistent with \cite{xingoldberg2016}, as the lead time gets longer while the sale price $r$ is fixed to be small, the base-stock policy cannot perform better.

\subsection{\sf Performance of $(s,q)$-Policies in {\bf POP}-System}\label{num-sq}
We now numerically study the performance of the $(s,q)$-policy. To this end, let $(s^*_o,q^*_o)=\arg\max_{(s,q)} {\cal C}^{(s,q)}(r,h;{\bf x})$.
To see the performance given by $(s^*_o,q^*_o)$-policy, we compared it with the best base-stock policy with level $s^*$ which uses a full information about the backorder and demand.
 Let $B^{s^*}$ be the long-run average backorders under  $s^*$-policy, which is not observable under $(s^*_o,q^*_o)$-policy for the {\bf POP}-system.
Tables \ref{tab 5:performance_s_q_policy_zhaoxuan} and \ref{tab 6:performance_s_q_policy_zhaoxuan} list the results of our experiments. From these two tables, we can conclude that the $(s,q)$-policy is effective across all the instances and all the sale price regimes. And the main reason of the effectiveness is that under the $s^*$-policy, the long-run average backorders $B^{s^*}$ is almost zero. This then naturally implies that the system under the $s^*$-policy will be similar to the system under $(s^*_o,q^*_o)$-policy with $q^*_o\approx 0$.

\begin{table}
    \centering
    \begin{tabular}{cccccccccc}
        \toprule
              &\multirow{2}{*}{Demand}&\multicolumn{3}{c}{$(s^*_o,q^*_o); (s^*,B^{s^*})$} &\multicolumn{3}{c}{$1-\frac{{\cal C}^{(s_o^*,q^*_o)}(r,h;{\bf x})}{{\cal C}^{s^*}(r,h;{\bf x})}$}\\
              & & $r=4$ & $r=8$ & $r=16$ & $r=4$ & $r=8$ & $r=16$  \\ \midrule

         \multirow{2}{*}{$\tau=0$} & Poisson(10)& (10,2);(10,0.88) & (12,1);(12,0.36) & (13,1);(13,0.22) &  0.17&  0.13& 0.03  \\

         &Bino(20,0.5)& (10,2);(10,0.62) & (11,1);(11,0.32) & (12,1);(12,0.16) &  0.41&  0.13 & 0.07 \\[0.10in]


         \multirow{2}{*}{$\tau=2$} & Poisson(10)& (29,2);(29,1.36) & (34,1);(32,0.71) & (34,2);(34,0.42) &  0.61&  0.38& 0.01  \\

         &Bino(20,0.5)& (30,1);(29, 0.96) & (32,1);(31,0.55) & (33,1);(33,0.27) & 0.14 &  0.16 & 0.04 \\[0.10in]


       \multirow{2}{*}{$\tau=4$}&Poisson(10)& (48,1);(47,1.61)& (52,0);(50, 1.08)& (55,0);(55,0.46)& 0.60&  0.05 &
       0.02\\

       &Bino(20,0.5) & (48,1);(48,1.14)& (52,0);(50,0.78)& (53,2);(53, 0.35)&  0.19& 0.07& 0.06\\

       \bottomrule
    \end{tabular}
    \caption{Performance of $(s,q)$-Policy ($h=1, p = 0.7$)}
    \label{tab 5:performance_s_q_policy_zhaoxuan}
\end{table}

\begin{table}
    \centering
        \begin{tabular}{cccccccccc}
        \toprule
              &\multirow{2}{*}{Demand}&\multicolumn{3}{c}{$(s^*_o,q^*_o); (s^*,B^{s^*})$} &\multicolumn{3}{c}{$1-\frac{{\cal C}^{(s_o^*,q^*_o)}(r,h;{\bf x})}{{\cal C}^{s^*}(r,h;{\bf x})}$}\\
              & & $r=4$ & $r=8$ & $r=16$ & $r=4$ & $r=8$ & $r=16$  \\ \midrule

         \multirow{2}{*}{$\tau=0$} & Poisson(10)& (12,0);(12,0.16) & (13,1);(13,0.10) & (14,1);(15,0.04) &  0.10&  0.10& 0.02  \\

         &Bino(20,0.5)& (11,1);(11,0.14) & (12,1);(12,0.06) & (13,0);(13,0.03) &  0.02&  0.01 & 0.01 \\[0.10in]


         \multirow{2}{*}{$\tau=2$} & Poisson(10)& (31,0);(31,0.28) & (34,0);(34,0.18) & (36,0);(36,0.08) &  0.03&  0.00& 0.01  \\

         &Bino(20,0.5)& (31,0);(31,0.19) & (33,0);(33,0.09) & (34,0);(35,0.04) &  -0.03 &  -0.02 & 0.02 \\[0.10in]


       \multirow{2}{*}{$\tau=4$}&Poisson(10)& (49,0);(49,0.34)& (53,0);(53, 0.19)& (57,0);(57,0.10)& 0.00&  0.00 &  0.00\\

       &Bino(20,0.5) & (50,0);(49,0.25)& (53,0);(53,0.11)& (55,0);(55, 0.07)&  -0.07& -0.08& 0.02\\

       \bottomrule
    \end{tabular}
    \caption{Performance of $(s,q)$-Policy ($h=1, p = 0.3$)}
    \label{tab 6:performance_s_q_policy_zhaoxuan}
\end{table}

\subsection{\sf Performance of Online Algorithm}\label{num-algorithm}
Finally, in this subsection, we  numerically show the effectiveness of the online algorithm we proposed. Inspired by \cite{russovanroy2014}, and \cite{lyuetal2021}, we will tune the algorithm with a parameter $\eta$, which will be a multiplier of the UCB-index. The parameter $\eta$ needs to be chosen to further balance the exploration and exploitation trade-off. In our numerical experiments, we consider the instances where the lead time $\tau\in\{2,4,6\}$ and the unit sale price $r\in\{4,8,16\}$. For initialization, we set the upper bound $\underline{s}$ to be the $0.99$-quantile of the lead-time demand $D_{[1,\tau+1]}$ when demands follow Poisson distribution and to be the upper bound of the lead-time demand $D_{[1,\tau+1]}$ when demands follow binomial distribution. The upper bound $\underline{q}$ is set to be $2\tau$ and we set $\varepsilon_1=\varepsilon_2 = 64$. During the implementation, following \cite{russovanroy2014}, and \cite{lyuetal2021}, we tune the parameter $\eta$ such that for $\tau = 2, 4, 6$, $\eta H_{\ell} =10^{-2\tau}$ for all epochs $1\le \ell \le T_{\mathcal{E}}(N)$, and we replace $H_{\ell-1}$ by $\eta H_{\ell-1}$ when calculating $F^{(s,q)}_\ell$ for all epochs $1\le \ell \le T_{\mathcal{E}}(N)$. In Tables \ref{tab 7 : avg_profit_and_regret_low_patience_online_algorithm_zhaoxuan} and \ref{tab 8 : avg_profit_and_regret_high_patience_online_algorithm_zhaoxuan}, we report the average profit and the percentage of the regret:
\begin{align*}
    \kappa(N):=\frac{{\sf Rg}(N)}{\mathbb{E}{\cal R}^{(s^*,q^*)}_{[1,N]}}\times 100\%.
\end{align*}
Tables \ref{tab 7 : avg_profit_and_regret_low_patience_online_algorithm_zhaoxuan} and \ref{tab 8 : avg_profit_and_regret_high_patience_online_algorithm_zhaoxuan} provide the numerical evidence for the effectiveness of our online algorithm.

\begin{table}
    \centering
    \resizebox{0.9\textwidth}{!}{
        \begin{tabular}{ccccccccccccc}
        \toprule
              &\multirow{2}{*}{Demand}&\multirow{2}{*}{\text{Price}} & \quad&\multicolumn{4}{c}{\text{Average Profit under Algorithm}}&\multicolumn{4}{c}{\text{Percentage of Regret $\kappa(N)$}}\\
              & & & & $N=20$ & $N=200$ & $N=500$ & $N=1000$ & $N=20$ & $N=200$ & $N=500$ & $N=1000$   \\ \cmidrule[0.5pt](l{.25em}r{.25em}){2-2}
              \cmidrule[0.5pt](l{.25em}r{.25em}){3-3}
              \cmidrule[0.5pt](l{.25em}r{.25em}){4-8}
              \cmidrule[0.5pt](l{.25em}r{.25em}){9-12}

         \multirow{8}{*}{$\tau=2$} & \multirow{3}{*}{Poisson(10)}& $r=4$ & & 25.75 & 27.75 & 27.89 & 27.95 & 24.57 & 18.99 & 18.69 & 18.53   \\
         &&$r=8$&& 54.00& 59.13 & 59.88 & 60.14 & 25.52 & 17.92 & 16.90 & 16.54  \\
         &&$r=16$&& 118.64 & 132.76 & 134.02 & 134.39 & 20.88 & 11.51 & 10.66 & 10.40\\
         &&&&&&&&\\
         &\multirow{3}{*}{Bino(20,0.5)}& $r=4$ && 25.66 & 27.39 & 27.53 & 27.57& 29.27 & 23.99 & 23.58 & 23.45  \\
         &&$r=8$&& 56.34 & 61.60 & 62.19 & 62.42 & 24.68 & 17.64 & 16.82 & 16.53\\
         &&$r=16$&& 119.53 & 136.00 & 137.46 & 137.92 & 21.71 & 10.81 & 9.89 & 9.60\\
         &&&&&&&&\\

         \multirow{6}{*}{$\tau=4$} &\multirow{3}{*}{Poisson(10)}& $r=4$ && 23.13 &29.84 &30.31 & 30.48 & 32.70 & 10.79 & 9.49 & 8.95 \\
         &&$r=8$&& 56.25 & 62.21 & 62.66 & 62.86 & 24.96 & 12.85 & 12.14 & 11.84\\
         &&$r=16$&& 123.76 & 131.43 & 132.27 & 132.51 & 16.71 & 11.48 & 10.88 & 10.72\\
         &&&&&&&&\\
         &\multirow{3}{*}{Bino(20,0.5)}& $r=4$ && 24.01 & 30.68 & 31.15 & 31.30 & 34.05 & 13.33 & 11.86 & 11.37  \\
         &&$r=8$&& 57.63 & 63.20 & 63.64 & 63.78 & 22.68 & 14.12 & 13.51 & 13.34 \\
         &&$r=16$&& 128.63 & 134.58 & 135.27 & 135.48 & 17.34 & 11.58 & 10.84 & 10.58\\
         &&&&&&&&\\

         \multirow{6}{*}{$\tau=6$} &\multirow{3}{*}{Poisson(10)}& $r=4$ && 18.16 & 30.56 & 31.36 & 31.67 & 54.85 & 7.81 & 5.35 & 4.33  \\
         &&$r=8$&& 52.14 & 64.15 & 64.95 & 65.26 & 29.08 & 8.80 & 7.63 & 7.17\\
         &&$r=16$&& 123.59 & 133.84 & 134.64 & 134.96 & 17.47 & 9.26 & 8.68 & 8.46\\
         &&&&&&&&\\
         &\multirow{3}{*}{Bino(20,0.5)}& $r=4$ && 18.40 & 31.64 & 32.60 & 32.87 & 55.93 & 9.46 & 6.73 & 5.92  \\
         &&$r=8$&& 53.93 & 65.78 & 66.57 & 66.84 & 28.85 & 9.72 & 8.70 & 8.27\\
         &&$r=16$&& 124.75 & 135.79 & 136.68 & 136.98 & 18.48 & 10.41 & 9.82 & 9.58\\
       \bottomrule
    \end{tabular}}
    \caption{Average Profit and Percentage of the Regret of the Online Algorithm ($h=1, p = 0.3$)}
    \label{tab 7 : avg_profit_and_regret_low_patience_online_algorithm_zhaoxuan}
\end{table}

\begin{table}
    \centering
    \resizebox{0.9\textwidth}{!}{
        \begin{tabular}{ccccccccccccc}
        \toprule
              &\multirow{2}{*}{Demand}&\multirow{2}{*}{\text{Price}} & \quad&\multicolumn{4}{c}{\text{Average Profit under Algorithm}}&\multicolumn{4}{c}{\text{Percentage of Regret $\kappa(N)$}}\\
              & & & & $N=20$ & $N=200$ & $N=500$ & $N=1000$ & $N=20$ & $N=200$ & $N=500$ & $N=1000$   \\ \cmidrule[0.5pt](l{.25em}r{.25em}){2-2}
              \cmidrule[0.5pt](l{.25em}r{.25em}){3-3}
              \cmidrule[0.5pt](l{.25em}r{.25em}){4-8}
              \cmidrule[0.5pt](l{.25em}r{.25em}){9-12}

         \multirow{8}{*}{$\tau=2$} & \multirow{3}{*}{Poisson(10)}& $r=4$ & & 26.27 & 28.05 & 28.28 & 28.35 & 25.51 & 20.96 & 20.38 & 20.24   \\
         &&$r=8$&& 57.73 & 62.09 & 62.37 & 62.55 & 20.62 & 15.29 & 14.74 & 14.47 \\
         &&$r=16$&& 121.00 & 131.86 & 132.80 & 133.26 & 20.80 & 13.43 & 12.83 & 12.48\\
         &&&&&&&&\\
         &\multirow{3}{*}{Bino(20,0.5)}& $r=4$ && 26.07 & 27.57 & 27.74 & 27.81 & 28.99 & 24.84 & 24.50 & 24.39  \\
         &&$r=8$&& 56.99 & 61.78 & 62.31 & 62.46 & 24.65 & 18.24 & 17.58 & 17.43\\
         &&$r=16$&& 121.78 & 137.64 & 138.95 & 139.22 & 21.01 & 10.73 & 9.89 & 9.74\\
         &&&&&&&&\\

         \multirow{6}{*}{$\tau=4$} &\multirow{3}{*}{Poisson(10)}& $r=4$ && 24.45 & 30.92 & 31.38 & 31.55 & 31.73 & 11.60 & 10.29 & 9.81 \\
         &&$r=8$&& 56.25 & 62.21 & 62.66 & 62.86 & 21.12 & 12.66 & 11.97 & 11.70\\
         &&$r=16$&& 123.76 & 131.43 & 132.27 & 132.51 & 16.16 & 11.35 & 10.82 & 10.63\\
         &&&&&&&&\\
         &\multirow{3}{*}{Bino(20,0.5)}& $r=4$ && 25.15 & 31.17 & 31.61 & 31.76 & 33.47 & 14.20 & 13.05 & 12.66  \\
         &&$r=8$&& 59.30 & 63.77 & 64.15 & 64.29 & 22.02 & 14.81 & 14.32 & 14.15 \\
         &&$r=16$&& 125.29 & 133.87 & 134.93 & 135.31 & 16.63 & 12.25 & 11.80 & 11.69\\
         &&&&&&&&\\

         \multirow{6}{*}{$\tau=6$} &\multirow{3}{*}{Poisson(10)}& $r=4$ && 19.71 & 31.89 & 32.75 & 33.04 & 50.04 & 7.80 & 5.32 & 4.49  \\
         &&$r=8$&& 55.13 & 66.27 & 67.09 & 67.37 & 25.82 & 8.68 & 7.59 & 7.18\\
         &&$r=16$&& 126.33 & 137.07 & 137.89 & 138.27 & 15.64 & 7.94 & 7.40 & 7.12\\
         &&&&&&&&\\
         &\multirow{3}{*}{Bino(20,0.5)}& $r=4$ && 20.72 & 32.56 & 33.40 & 33.68 & 50.03 & 9.50 & 7.08 & 6.34  \\
         &&$r=8$&& 56.59 & 66.99 & 67.76 & 68.00 & 25.88 & 10.14 & 9.08 & 8.79\\
         &&$r=16$&& 128.57 & 137.64 & 138.60 & 139.00 & 16.25 & 9.75 & 9.09 & 8.89\\
       \bottomrule
    \end{tabular}}
    \caption{Average Profit and Percentage of the Regret of the Online Algorithm ($h=1, p = 0.7$)}
    \label{tab 8 : avg_profit_and_regret_high_patience_online_algorithm_zhaoxuan}
\end{table}

\section{\sf Conclusions}\label{conclusion}
We study a stochastic inventory system with partial backorders induced by customer abandonment behavior. Even with the ``overshooting" phenomenon induced by the random patience time, the base-stock policy is proved to be \textbf{uniformly} (asymptotically) optimal for all {\bf P}-systems, where the decision maker can perfectly observe the system states. We propose a base-stock type heuristic policy and numerical evidences showing the necessity of considering the patience of the customers.

The analysis is then extended to the model where the backorders, i.e., the customers who are patient during the stock-out periods, are hidden. In this partially observed system ({\bf POP}-system) with hidden backorders, leveraging on the optimality results we proved for the {\bf P}-system, we propose a new family of policies called $(s,q)$-policy. The asymptotic optimality of the $(s,q)$-policy is also provided with a heuristic choice of the parameters. Finally, we consider an online inventory control problem in the {\bf POP}-system, where the decision maker doesn't have perfect demand information and the system states are only partially observed. A UCB-type algorithm is developed based on the traditional multi-armed bandits algorithm and the geometrical ergodicity of the system state under $(s,q)$-policy developed in this work. We prove the regret of the algorithm is $O(N^{\frac{1}{2}+\delta})$.

\section*{\sf Acknowledgements:}

Andrew Lim is supported by the  Ministry of Education, Singapore, under its 2021 Academic Research Fund Tier 2 grant call (Award ref: MOE-T2EP20121-0014).



\newpage
\begin{APPENDICES}
\section{Proofs}
\renewcommand{\thetable}{A-\arabic{table}}
\renewcommand{\theequation}{A-\arabic{equation}}
\setcounter{equation}{0}
\setcounter{table}{0}

\noindent
{\it Proof of {\bf Proposition \ref{Lemma_Connect P-system and B-system_Dynamic}}}: As {\bf B}- and {\bf P}-systems start with the same state and implement the same base-stock policy,
\begin{align}
N^s_i &=I_1+Q_{[1-\tau, i-\tau]}-B_0-D_{[1,i-1]}+L^s_{[1,i-1]}\nonumber\\
& \ge I_1+Q_{[1-\tau, i-\tau]}-B_0-D_{[1,i-1]}=\overline{N}^s_i \ \ \mbox{for $i\leq \tau$}.\label{Prop-1}
\end{align}
 Note that for {\bf B}- and {\bf P}-systems, their order quantities at period-1 are same $(Q^s_1=\overline{Q}^s_1)$ as their initial states and base-stock levels are identical.  Further, at period-$i$ with $i\geq 2$, for the {\bf P}-system, the ordering quantity $Q^s_i$  should at most be $D_{i-1}$, which is exactly the ordering quantity for the {\bf B}-system. To see this, the inventory position For the {\bf P}-system will decrease at most the level of $D_{i-1}$ in period-$i$ (this happens when there is no loss of unmet demands in this period, $L^s_{i-1}=0$), which means, to push the inventory position back to $s$, the ordering quantity of {\bf P}-system is at most $D_{i-1}$. Thus,
\begin{align}
N_i^s &=I_1+Q_{[1-\tau,0]}+Q^s_{[1, i-\tau]}-B_0-D_{[1,i-1]}+L^s_{[1,i-1]}\nonumber\\
& =I_1+Q_{[1-\tau,0]}+Q^s_1+Q^s_{[2, i-\tau]}-B_0-D_{[1,i-\tau-1]}-D_{[i-\tau,i-1]}+L^s_{[1,i-1]}\nonumber\\
& \geq I_1+Q_{[1-\tau,0]}+Q^s_1-B_0-D_{[i-\tau,i-1]}+L^s_{[1,i-1]}\nonumber\\
& \geq I_1+Q_{[1-\tau,0]}+Q^s_1-B_0-D_{[i-\tau,i-1]}\nonumber\\
& = I_1+Q_{[1-\tau,0]}+\overline{Q}^s_1-B_0-D_{[i-\tau,i-1]}\nonumber\\
& =I_1+Q_{[1-\tau,0]}+\overline{Q}^s_1+\overline{Q}^s_{[2, i-\tau]}-D_{[1,i-\tau-1]}-B_0-D_{[i-\tau,i-1]} \nonumber\\
& =I_1+Q_{[1-\tau,0]}+\overline{Q}^s_{[1, i-\tau]}-D_{[1,i-1]}-B_0 =\overline{N}_i^s \ \ \mbox{for $i\geq 1+\tau$},\label{Prop-2}
\end{align}
which gives the first inequality of the first relation in (\ref{connect}).

Now we look at the second inequality in the first relation in (\ref{connect}). It holds for $i=1$ by the assumption that both systems have the same initial state.

For $i\geq 2$, observe that $N_i^s=N^s_{i-1}-D_{i-1}+L_{i-1}^s+Q^s_{i-\tau}$ with the convention in which $Q^s_k=Q_k$ for $1-\tau\leq k\leq 0$.

When $2\leq i\leq \tau$, by the telescope summing from $2$ to $i$, we have
\begin{align}
   N_i^s &= N^s_{1} + Q^s_{[2-\tau, i-\tau]} -D_{[1, i-1]}+L^s_{[1, i-1]}\nonumber \\
          &= I_1+Q_{[1-\tau,i-\tau]}-B_0-D_{[1, i-1]} +L^s_{[1, i-1]}\nonumber \\
          &=\overline{N}^s_i+L^s_{[1, i-1]}.\label{Prop-3}
\end{align}
This proves the second inequality of the first relation in (\ref{connect}) for $i\leq \tau$.

When $i\geq 1+\tau$, by the telescope summing from $(i-\tau+1)$ to $i$ for $N_i^s=N^s_{i-1}-D_{i-1}+L_{i-1}^s+Q^s_{i-\tau}$, we have
$N_i^s=N^s_{i-\tau} + Q^s_{[i-2\tau+1, i-\tau]} -D_{[i-\tau, i-1]}+L^s_{[i-\tau, i-1]}$. Under the base-stock policy with level $s$, we know that
$Q^s_{i-\tau}=(s-N^s_{i-\tau}-Q^s_{[i-2\tau+1,i-\tau-1]})^+$. Hence, by $O^s_{i-\tau}=(N^s_{i-\tau}+Q^s_{[i-2\tau+1,i-\tau-1]}-s)^+$,
\begin{align*}
   N_i^s  &= N^s_{i-\tau} + Q^s_{[i-2\tau+1, i-\tau]} -D_{[i-\tau, i-1]}+L^s_{[i-\tau, i-1]} \\
          &\leq s+O^s_{i-\tau}-D_{[i-\tau, i-1]} +L^s_{[i-\tau, i-1]} \\
          &=\overline{N}^s_i+L^s_{[i-\tau, i-1]} +O_{i-\tau}^s.
\end{align*}
This proves the second inequality of the first relation in (\ref{connect}) for $i\geq 1+\tau$.

The second relation $B^s_i+L^s_i\le\overline{B}_i^s$ of (\ref{connect}) directly follows from its first one $\overline{N}^s_i\le N_i^s$ by observing that $B_i^s+L^s_i=\big(D_i -N_i^s\big)^+\le(D_i-\overline{N}_i^s)^+=\overline{B}_i^s$.  So we have the proposition.\hfill$\Box$

\

\noindent
{\it Proof of {\bf Proposition \ref{estimation of the overshooting_lemma}}}: We first show that for all $i\ge1$, the overshooting $O_i^s$ is dominated by $B^s_{\ell(i)-1}$.
If $\ell(i)=i$, then there is no overshooting by the definition of $\ell(i)$, and the proposition holds.
Thus we just consider the case when $\ell(i)<i$, which then implies
that the inventory position at the periods from $\ell(i)+1$ to $i$ are all strictly above $s$. If so, then we will not place any order during the periods: $\ell(i)+1,\ldots,i$. During these periods, the only way to push the inventory position at the beginning of period-$i$, $I^s_i+Q^s_{i-\tau}-B^s_{i-1}$, higher is to make
\begin{align}
L^s_{[\ell(i)+1,k]}\geq D_{[\ell(i)+1,k]} \ \ \mbox{for $k=\ell(i)+1,\ldots,i-1$},
\label{hanqin-1}
\end{align}
which could equivalently be taken as during these periods, no demand is satisfied, and some of the backlogged demands from $B^s_{\ell(i)-1}$ also abandon the system. And the inventory position is exactly pulled up by the demands which is originally backordered but now lost. Hence, during periods $\ell(i)+1,\ldots,i-1$, the number of demands that contributes to the increase of the inventory position should be at most $B^s_{\ell(i)-1}$.

Finally, we prove that $O_i^s$ is also dominated by
\[
\left(B_{\ell(i)-1}^s-D_{[\ell(i)+\tau,i-1]}\right)^+\cdot \mathbb{I}_{\{\ell(i)<i-\tau\}}+B_{\ell(i)-1}^s\cdot \mathbb{I}_{\{\ell(i)\geq i-\tau\}}.
\]
 If $\ell(i)+\tau<i$,  there is no replenishment between periods $(\ell(i)+1)$ to $i$ and there is overshooting every period in between. At the end of period-$(\ell(i)+\tau)$, all  on-order inventory has been received and all backorders are cleared. From period $(\ell(i)+\tau+1)$ to $i$, the inventory position equals the on-hand inventory level and is always above $s$. It also implies there is no lost sale from period $(\ell(i)+\tau)$ to $(i-1)$. We then upper bound the overshooting as follows
 \begin{align*}
 O^s_i &=I_i^s-s \\
 &= I^s_{\ell(i)}+Q^s_{[\ell(i)-\tau,\ell(i)]}-B^s_{\ell(i)-1}
 -D_{[\ell(i),i-1]}+L^s_{[\ell(i),i-1]}-s\\
     &= L^s_{[\ell(i),\ell(i)+\tau-1]}-D_{[\ell(i),\ell(i)+\tau-1]}-D_{[\ell(i)+\tau,i-1]}\\
     &\le B_{\ell(i)-1}^s- D_{[\ell(i)+\tau,i-1]},
 \end{align*}
where the last inequality is from (\ref{hanqin-1}). Therefore, we have the proposition.\hfill$\Box$

\

\noindent
 {\it Proof of {\bf Proposition \ref{optimal policy of zero lead time_Theorem}}}: First, note that by Wald's equation, $\mathbb{E}\xi_{1;[1,(D_1-x)^+]}=p\cdot \mathbb{E}(D_1-x)^+$.  Then,
 \begin{eqnarray*}
 &&{\cal L}^o(x)=-rx^-+\mathbb{E}\Big[r \left(x^+\wedge D_1\right)-h\cdot(x -D_1)^+ + \alpha r\cdot \Big((x-D_1)^+-\xi_{1;[1,(D_1-x)^+]}\Big)^-\Big]\\
 && \ \ \ \ \ \ \ \ =-rx^-+\mathbb{E}\Big[r \left(x^+\wedge D_1\right)-h\cdot(x -D_1)^+ + \alpha r\cdot \xi_{1;[1,(D_1-x)^+]}\Big]\\
 && \ \ \ \ \ \ \ \ =-rx^-+\mathbb{E}\Big[r \left(x^+\wedge D_1\right)-h\cdot(x -D_1)^+ + \alpha pr\cdot(D_1-x)^+\Big]\\
 && \ \ \ \ \ \ \ \ =-rx^-+ rd +\mathbb{E}\Big[-r \left(D_1-x^+\right)^+-h\cdot(x -D_1)^+ + \alpha p r\cdot(D_1-x)^+\Big].
 \end{eqnarray*}
 The concavity of ${\cal L}^o(\cdot)$ directly follows from the concavity of the one-period profit for the classic newsvendor problem and the convexity of  $r\cdot x^-$.

 To obtain the quasi-concavity of ${\cal C}^{\alpha,o}(r,h; \cdot)$ with the global maximizer $s^{*,\alpha}$, the value iterating approach is applied. To this end, consider a sequence of function pairs defined by that for $k\ge1$,
\begin{eqnarray}
{\cal G}^{\alpha, (k-1)}(x) &=& {\cal L}^o(x) +\alpha\cdot \mathbb{E}\Big[{\cal C}^{\alpha,(k-1)}\Big(r,h; (x-D_1)^+-\xi_{1;[1,(D_1-x)^+]}\Big)\Big] \ \mbox{with} \ {\cal C}^{\alpha,(0)}(r,h;x)\equiv0,\label{hanqin-6}\\
{\cal C}^{\alpha,(k)}(x) &=& \max_{y\ge x}\Big\{ {\cal G}^{\alpha,(k-1)} (y) \Big\}. \label{hanqin-7}
\end{eqnarray}
We first establish that for each $k\geq 1$, ${\cal C}^{\alpha, (k)}(\cdot)$ is quasi-concave with the global maximizer $s^{*,\alpha}$.
In view of (\ref{hanqin-7}), if ${\cal G}^{\alpha,(k-1)} (\cdot)$ is quasi-concave with the global maximizer $s^{*,\alpha}$, then
\begin{eqnarray}
{\cal C}^{\alpha,(k)}(x) =\left\{
\begin{array}{ll}
{\cal G}^{\alpha,(k-1)} (s^{*,\alpha}), &\mbox{if $x\leq s^{*,\alpha}$},\\
{\cal G}^{\alpha,(k-1)} (x), &\mbox{if $x> s^{*,\alpha}$}.
\end{array}
\right.\label{hanqin-8}
\end{eqnarray}
This shows that the quasi-concavity of ${\cal G}^{\alpha, (k-1)}(\cdot)$, and its global maximizer $s^{*,\alpha}$ imply that  ${\cal C}^{\alpha,(k)}(\cdot)$ is quasi-concave with the global maximizer $s^{*,\alpha}$.
Thus, it suffices to show that for each $k\geq 0$, ${\cal G}^{\alpha, (k)}(\cdot)$ is quasi-concave with the global maximizer $s^{*,\alpha}$.

 By (\ref{hanqin-6}), the concavity of ${\cal L}^o(\cdot)$ proved above directly gives that ${\cal G}^{\alpha,(0)}(\cdot)$ is quasi-concave with the global maximizer $s^{*,\alpha}$.  By an induction, suppose that ${\cal G}^{\alpha, (k-1)}(\cdot)$ is quasi-concave with the global maximizer $s^{*,\alpha}$, we then try to show that ${\cal G}^{\alpha, (k)}(\cdot)$ is also quasi-concave with the same maximizer $s^{*,\alpha}$.
 Recall that for any $x<y$, with probability one,
$(x-D_1)^+-\xi_{1;[1,(D_1-x)^+]}\leq (y-D_1)^+-\xi_{1;[1,(D_1-y)^+]}$. Hence, we have that
${\cal L}^o(x) +\alpha\cdot \mathbb{E}\Big[{\cal C}^{\alpha,(k-1)}\Big((x-D_1)^+-\xi_{1;[1,(D_1-x)^+]}\Big)\Big]$ is also quasi-concave with the global maximizer $s^{*,\alpha}$. Thus, the quasi-concavity and global maximizer of ${\cal G}^{\alpha, (k)}(\cdot)$ directly follows from (\ref{hanqin-6}).

After obtaining that ${\cal C}^{\alpha,(k)}(\cdot)$ is quasi-concave, and has the global maximizer $s^{*,\alpha}$,
by a standard DP argument, ${\cal C}^{\alpha,o}(\cdot)=\lim_{k\to\infty}{\cal C}^{\alpha,(k)}(\cdot)$, which implies that ${\cal C}^{\alpha,o}(\cdot)$ is quasi-concave with the global minimizer $s^{*,\alpha}$.
Hence, the proposition is proved. \hfill$\Box$\\[-0.15in]

\

\noindent
{\it Proof of {\bf Proposition \ref{Uniform bounded overshooting_lemma}}}:
In view of Proposition \ref{estimation of the overshooting_lemma}, decompose the overshooting at each period into two parts,
\begin{align}
\limsup\limits_{N\to\infty}\frac{1}{N}\mathbb{E} O^s_{[1,N]}&=\limsup\limits_{N\to\infty}\frac{1}{N}\mathbb{E} O^s_{[4\tau ,N]}\nonumber\\
&\le\limsup\limits_{N\to\infty}\frac{1}{N}\sum_{i=4\tau}^N\mathbb{E}\Big[O_{i}^s\mathbb{I}_{\{\ell(i)<i-3\tau\}}\Big]+
    \limsup\limits_{N\to\infty}\frac{1}{N}\sum_{i=4\tau}^N\mathbb{E}\Big[O_{i}^s\mathbb{I}_{\{\ell(i)\geq i-3\tau\}}\Big].\label{hanqin-20}
\end{align}
By Propositions \ref{Lemma_Connect P-system and B-system_Dynamic} and \ref{estimation of the overshooting_lemma}, for the second term on the right-hand-side of (\ref{hanqin-20}),
\begin{align}
    \limsup\limits_{N\to\infty}\frac{1}{N} \sum_{i=4\tau}^N\mathbb{E}\Big[O_{i}^s\mathbb{I}_{\{\ell(i)\geq i-3\tau\}}\Big]&\le  \limsup\limits_{N\to\infty}\frac{1}{N}\sum_{i=4\tau}^N\mathbb{E}\Big[B^s_{\ell(i)-1}\mathbb{I}_{\{\ell(i)\geq i-3\tau\}} \Big]\nonumber \\
    &\le  \limsup\limits_{N\to\infty}\frac{1}{N}\sum_{i=4\tau}^N\mathbb{E}\Big[\max_{i-3\tau-1 \le k \le i-1}B^s_{k}\Big]\nonumber\\
    &\le  \limsup\limits_{N\to\infty}\frac{1}{N}\sum_{i=4\tau}^N\mathbb{E}\Big[\max_{i-3\tau-1 \le k \le i-1}\overline{B}^s_{k}\Big]\nonumber\\
    &\le  \limsup\limits_{N\to\infty}\frac{1}{N}\sum_{i=4\tau}^N\mathbb{E}\Big[\overline{B}^s_{[i-3\tau-1, i-1]}\Big]\nonumber\\
    &\le (3\tau+1)\mathbb{E}\Big(B_0+D_{[1,\tau]} - s\Big)^+,\label{hanqin-21}
\end{align}
and for the first term on the right-hand-side of (\ref{hanqin-20}), from Proposition \ref{estimation of the overshooting_lemma},
\begin{align}
    \limsup\limits_{N\to\infty}&\frac{1}{N}\sum_{i=4\tau}^N\mathbb{E}\Big[\Big(O^s_i\mathbb{I}_{\{\ell(i)<i-3\tau\}}\Big)\Big]\nonumber\\
    &\le\limsup\limits_{N\to\infty}\frac{1}{N}\sum_{i=4\tau}^N\mathbb{E}\Big[\Big(B^s_{\ell(i)-1}-D_{[\ell(i)+\tau,i-1]}\Big)^+\mathbb{I}_{\{\ell(i)<i-3\tau\}}\Big]\nonumber\nonumber\\
     &\le\limsup\limits_{N\to\infty}\frac{1}{N}\sum_{i=4\tau}^N\mathbb{E}\Big[\max\limits_{1\le k\le \tau}\Big(
     B^s_k-D_{[k+\tau, i-1]}\Big)^+\Big] \nonumber\\
    &\quad+\limsup\limits_{N\to\infty}\frac{1}{N}\sum_{i=4\tau+1}^N\mathbb{E}\Big[\max\limits_{\tau+1\le k\le i-3\tau}
    \Big(B^s_k-D_{[k+\tau,i-1]}\Big)^+\Big]\nonumber\\
    &\le\limsup\limits_{N\to\infty}\frac{1}{N}\sum_{i=4\tau}^N\mathbb{E}\Big[\max\limits_{1\le k\le \tau}\Big(
     \overline{B}^s_k-D_{[k+\tau, i-1]}\Big)^+\Big] \nonumber\\
    &\quad+\limsup\limits_{N\to\infty}\frac{1}{N}\sum_{i=4\tau+1}^N\mathbb{E}\Big[\max\limits_{\tau+1\le k\le i-3\tau}
    \Big(\overline{B}^s_k-D_{[k+\tau,i-1]}\Big)^+\Big].\label{hanqin-22}
\end{align}
For the first term on the right-hand side of (\ref{hanqin-22}),
\begin{align}
& \limsup\limits_{N\to\infty}\frac{1}{N}\sum_{i=4\tau}^N\mathbb{E}\Big[\max\limits_{1\le k\le \tau}\Big(
     \overline{B}^s_k-D_{[k+\tau, i-1]}\Big)^+\Big] \nonumber\\
& \ \ \ \le \limsup\limits_{N\to\infty}\frac{1}{N}\sum_{i=4\tau}^N\mathbb{E}\Big[\max\limits_{1\le k\le \tau}\Big(
     (B^s_0+D_{[1,\tau]})-D_{[k+\tau, i-1]}\Big)^+\Big] \nonumber\\
& \ \ \ \le \limsup\limits_{N\to\infty}\frac{1}{N}\sum_{i=4\tau}^N\mathbb{E}\Big(
     (B^s_0+D_{[1,\tau]})-D_{[2\tau, i-1]}\Big)^+,\label{hanqin-23}
\end{align}
and for the second term on the right-hand side of (\ref{hanqin-22}),
\begin{align}
& \limsup\limits_{N\to\infty}\frac{1}{N}\sum_{i=4\tau}^N\mathbb{E}\Big[\max\limits_{\tau+1\le k\le i-3\tau}
    \Big(\overline{B}^s_k-D_{[k+\tau,i-1]}\Big)^+\Big]\nonumber\\
& \ \ \ = \limsup\limits_{N\to\infty}\frac{1}{N}\sum_{i=4\tau}^N\mathbb{E}\Big[\max\limits_{\tau+1\le k\le i-3\tau}
    \Big(D_{[k-\tau,k]}-s-D_{[k+\tau,i-1]}\Big)^+\Big]\nonumber\\
& \ \ \ \le \limsup\limits_{N\to\infty}\frac{1}{N}\sum_{i=4\tau+1}^N\sum_{k=\tau+1}^{i-3\tau}\mathbb{E}\Big[
    \Big(D_{[k-\tau,k]}-s-D_{[k+\tau,i-1]}\Big)^+\Big]\nonumber\\
& \ \ \ = \limsup\limits_{N\to\infty}\frac{1}{N}\sum_{i=4\tau+1}^N\sum_{k=\tau+1}^{i-3\tau}\mathbb{E}\Big[
    \Big(D_{[1,\tau+1]}-s-D_{[\tau+2,i-k+1]}\Big)^+\Big].\label{hanqin-24}
\end{align}
The proposition directly follows from (\ref{hanqin-20})-(\ref{hanqin-24}), and the following lemma.
\hfill$\Box$

\begin{lemma}
    \label{tech-lem-1}
Assume that there exists $\gamma>0$ such that $\mathbb{E}e^{\gamma D_1}<\infty$. Under a fixed base-stock policy $s$, and an initial state
$(I^s_1,B^s_0,Q^s_{1-\tau},\ldots,Q^s_0)\in {\cal X}^s$, there exists $\gamma_0\in (0, \gamma)$ such that for any $i\geq 4\tau$,
\begin{align*}
&\sum_{k=\tau+1}^{i-3\tau}\mathbb{E}\Big[
    \Big(D_{[1,\tau+1]}-s-D_{[\tau+2,i-k+1]}\Big)^+\Big] \leq \frac{4}{\gamma\gamma_0d}  \exp\Big(-\gamma s+\frac{1}{2}\gamma_0\tau d \Big)  \Big(\mathbb{E}e^{\gamma D_1}\Big)^{\tau+1},\\
&\mathbb{E}\Big((B^s_0+D_{[1,\tau]})-D_{[2\tau, i-1]}\Big)^+\le \frac{2}{\gamma} e^{\gamma B^s_0} \times\exp\Big(-\frac{1}{2}\gamma_0[i-2\tau]d\Big)\Big(\mathbb{E}e^{\gamma D_1}\Big)^{\tau}.
\end{align*}
\end{lemma}

\noindent
{\it Proof}: First we prove the first inequality in the lemma.  For any positive integer $m$,
\begin{align}
     &\mathbb{P} \Big[\Big(D_{[1,\tau+1]}-s-D_{[\tau+2,i-k+1]}\Big)^+\ge m \Big]\nonumber\\
     &\qquad=\mathbb{P} \Big[D_{[1,\tau+1]}-s-D_{[\tau+2,i-k+1]} \ge m \Big]\nonumber\\
     &\qquad=\mathbb{P}\Big[ D_{[1,\tau+1]}-s-D_{[\tau+2,i-k+1]} \ge m \Big| D_{[\tau+2,i-k+1]}>\frac{1}{2}(i-k-\tau)d \Big]\nonumber\\
     &\qquad\quad \ \ \times \mathbb{P}\Big[ D_{[\tau+2,i-k+1]}>\frac{1}{2}(i-k-\tau)d \Big]\nonumber\\
     &\qquad\quad+\mathbb{P} \Big[ D_{[1,\tau+1]}-s-D_{[\tau+2,i-k+1]} \ge m \Big| D_{[\tau+2,i-k+1]}\leq \frac{1}{2}(i-k-\tau)d\Big] \nonumber\\
     &\qquad\quad \ \ \times  \mathbb{P}\Big[D_{[\tau+2,i-k+1]}\leq\frac{1}{2}(i-k-\tau)d \Big] \nonumber\\
     &\qquad \leq \mathbb{P}\Big[ D_{[1,\tau+1]}\ge m+s+\frac{1}{2}(i-k-\tau)d \Big| D_{[\tau+2,i-k+1]}>\frac{1}{2}(i-k-\tau)d \Big]\nonumber\\
     &\qquad\quad \ \ \times \mathbb{P}\Big[ D_{[\tau+2,i-k+1]}>\frac{1}{2}(i-k-\tau)d \Big]\nonumber\\
     &\qquad\quad+\mathbb{P} \Big[ D_{[1,\tau+1]} \ge m+s \Big| D_{[\tau+2,i-k+1]}\leq \frac{1}{2}(i-k-\tau)d\Big] \times  \mathbb{P}\Big[D_{[\tau+2,i-k+1]}\leq\frac{1}{2}(i-k-\tau)d \Big]\nonumber \\
     &\qquad=\mathbb{P}\Big[ D_{[1,\tau+1]}\ge m+s+\frac{1}{2}(i-k-\tau)d\Big]
     \times \mathbb{P}\Big[ D_{[\tau+2,i-k+1]}>\frac{1}{2}(i-k-\tau)d \Big]\nonumber\\
     &\qquad\quad+\mathbb{P} \Big[ D_{[1,\tau+1]}\geq m+s \Big] \times  \mathbb{P}\Big[D_{[\tau+2,i-k+1]}\leq\frac{1}{2}(i-k-\tau)d \Big]\nonumber \\
     & \qquad \le \mathbb{P}\Big[ D_{[1,\tau+1]}\ge m+s+\frac{1}{2}(i-k-\tau)d\Big]\nonumber\\
     & \qquad\quad+ \mathbb{P} \Big[ D_{[1,\tau+1]}\geq m+s \Big] \times  \mathbb{P}\Big[D_{[\tau+2,i-k+1]}\leq\frac{1}{2}(i-k-\tau)d \Big]\nonumber \\
     & \qquad \le \exp\Big(-\gamma \Big[m+s+\frac{1}{2}(i-k-\tau)d\Big]\Big)\Big(\mathbb{E}e^{\gamma D_1}\Big)^{\tau+1}\nonumber\\
     &  \qquad\quad+ \exp\Big(-\gamma[m+s]\Big)\Big(\mathbb{E}e^{\gamma D_1}\Big)^{\tau+1}\times \mathbb{P}\Big[D_{[\tau+2,i-k+1]}-(i-k-\tau)d\leq-\frac{1}{2}(i-k-\tau)d \Big].\label{hanqin-24-1}
\end{align}
Note that from Lemma 2.6.2 in \cite{durrett2010}, there exits $\gamma_0$ in the interval $(0, \gamma)$ such that
\begin{align*}
&\mathbb{P}\Big[D_{[\tau+2,i-k+1]}-(i-k-\tau)d\leq-\frac{1}{2}(i-k-\tau)d \Big]\\
& \ \ \ = \mathbb{P}\Big[(i-k-\tau)d-D_{[\tau+2,i-k+1]}\geq\frac{1}{2}(i-k-\tau)d \Big]\\
& \ \ \ \leq \exp\Big(-\frac{1}{2}\gamma_0 (i-k-\tau)d\Big).
\end{align*}
Hence, from (\ref{hanqin-24-1}),
\begin{align}
     &\mathbb{P} \Big[\Big(D_{[1,\tau+1]}-s-D_{[\tau+2,i-k+1]}\Big)^+\ge m \Big]\nonumber\\
     & \ \ \ \leq \exp\Big(-\gamma \Big[m+s+\frac{1}{2}(i-k-\tau)d\Big]\Big)\Big(\mathbb{E}e^{\gamma D_1}\Big)^{\tau+1}\nonumber\\
     & \ \ \ \ + \exp\Big(-\gamma[m+s]\Big)\times\exp\Big(-\frac{1}{2}\gamma_0[i-k-\tau]d\Big)\Big(\mathbb{E}e^{\gamma D_1}\Big)^{\tau+1}.\label{hanqin-25-1}
\end{align}
This implies that
\begin{align*}
& \mathbb{E} \Big(D_{[1,\tau+1]}-s-D_{[\tau+2,i-k+1]}\Big)^+\\
& \ \ \ =\sum_{m=1}^\infty  \mathbb{P} \Big[\Big(D_{[1,\tau+1]}-s-D_{[\tau+2,i-k+1]}\Big)^+\ge m \Big]\\
& \ \ \ \le 2 e^{-\gamma s} \Big(\mathbb{E}e^{\gamma D_1}\Big)^{\tau+1}\exp\Big(-\frac{1}{2}\gamma_0 [i-k-\tau]d\Big)
\cdot \Big(\sum_{m=1}^\infty e^{-\gamma m}\Big).
\end{align*}
From this, we have
\begin{align}
& \sum_{k=\tau+1}^{i-3\tau}\mathbb{E}\Big[
    \Big(D_{[1,\tau+1]}-s-D_{[\tau+2,i-k+1]}\Big)^+\Big]\nonumber\\
    & \ \ \ \leq \sum_{k=\tau+1}^{i-3\tau}2 e^{-\gamma s} \Big(\mathbb{E}e^{\gamma D_1}\Big)^{\tau+1}\exp\Big(-\frac{1}{2}\gamma_0 [i-k-\tau]d\Big)
\cdot \Big(\sum_{m=1}^\infty e^{-\gamma m}\Big)\nonumber\\
& \ \ \ =2 \exp\Big(-\gamma s+\frac{1}{2}\gamma_0 \tau d \Big)  \Big(\mathbb{E}e^{\gamma D_1}\Big)^{\tau+1}\cdot \Big(\sum_{m=1}^\infty e^{-\gamma m}\Big) \sum_{k=\tau+1}^{i-3\tau}\exp\Big(-\frac{1}{2}\gamma_0 [i-k]d\Big)\nonumber\\
& \ \ \ \leq \frac{4}{\gamma\gamma_0d}  \exp\Big(-\gamma s+\frac{1}{2}\gamma_0\tau d \Big)  \Big(\mathbb{E}e^{\gamma D_1}\Big)^{\tau+1}.\label{hanqin-26-1}
\end{align}
This completes the proof of the first inequality in the lemma.

For the second inequality in the lemma, going along the line of establishing (\ref{hanqin-25-1}),
\begin{align*}
     &\mathbb{P} \Big[\Big(D_{[1,\tau]}+B^s_0-D_{[2\tau,i-1]}\Big)^+\ge m \Big]\nonumber\\
     & \ \ \ \ \leq 2e^{\gamma B^s_0} \exp(-\gamma m)\Big(\mathbb{E}e^{\gamma D_1}\Big)^{\tau}\times\exp\Big(-\frac{1}{2}\gamma_0[i-2\tau]d\Big).
\end{align*}
Using this, similar to (\ref{hanqin-26-1}), we have the first inequality in the lemma. \hfill$\Box$

\

\noindent
{\it Proof of {\bf Proposition \ref{upper bound of base-stock P-system_prop}}}:
By Proposition \ref{Lemma_Connect P-system and B-system_Dynamic}, for all $i\ge\tau+1$, we have
\begin{align}
N_i^s\le \overline{N}_i^s+L^s_{[i-\tau, i-1]}+O_{i-\tau}^s=s-D_{[i-\tau,i-1]}+L^s_{[i-\tau, i-1]}+O_{i-\tau}^s. \label{hanqin-26}
\end{align}
For a given state $X_i^s=(I^s_i,B^s_{i-1},Q^s_{i-\tau},\ldots,Q^s_{i-1})$ at period-$i$ For the {\bf P}-system, the profit at period-$i$ can be written as
\begin{align*}
&{\cal P}^s_i(I^s_i,B^s_{i-1},Q^s_{i-\tau})\nonumber\\
& \ \ \ = r\cdot \left((I^s_i+Q^s_{i-\tau})\wedge (B^s_{i-1}+D_i)\right)-h\cdot \left(I^s_i+Q^s_{i-\tau}-B^s_{i-1}-D_i\right)^+\nonumber\\
& \ \ \ =r \cdot(B^\pi_{i-1}+D_i)-r\cdot\Big(B_{i-1}^s+D_i-I^s_i-Q_{i-\tau}^s\Big)^+-h\cdot(N_i^s-D_i)-h\cdot(D_i-N^s_i)^+\nonumber\\
& \ \ \ =r\cdot(B^s_{i-1}-B^s_i)+r\cdot D_i -h\cdot(N_i^s-D_i)-(h+r)\cdot L^s_i-h\cdot B^s_i.
\end{align*}
Here we use the facts: $B^s_i+L_i^s=(B_{i-1}^s+D_i-I^s_i-Q_{i-\tau}^s)^+$ and $N_i^s=I^s_i+Q_{i-\tau}^s-B_{i-1}^s$.
It follows from (\ref{hanqin-26}) that for $i\ge \tau+1$,
\begin{align}
 {\cal P}^s_i(I^s_i,B^s_{i-1},Q^s_{i-\tau}) & \ge r\cdot(B^s_{i-1}-B^s_i)+r\cdot D_i-h\cdot \Big(s-D_{[i-\tau,i-1]}+L^s_{[i-\tau, i-1]}+O_{i-\tau}^s-D_i\Big)\nonumber\\
 & \  -(h+r)\cdot L^s_i-h\cdot B^s_i\nonumber\\
 & =r\cdot(B^s_{i-1}-B^s_i)+r\cdot D_i -h \cdot (s-D_{[i-\tau,i]})-h\cdot\Big( L^s_{[i-\tau, i-1]}+O_{i-\tau}^s\Big)\nonumber\\
 & \  -(h+r)\cdot L^s_i-h\cdot B^s_i\nonumber\\
 & = r\cdot(B^s_{i-1}-B^s_i)+r\cdot D_i -h \cdot (s-D_{[i-\tau,i]})^+ + h \cdot(D_{[i-\tau,i]}-s)^+\nonumber\\
 & \  -h\cdot\Big( L^s_{[i-\tau, i-1]}+O_{i-\tau}^s\Big)-(h+r)\cdot L^s_i-h\cdot B^s_i\nonumber\\
 & \ge r\cdot(B^s_{i-1}-B^s_i)+r\cdot D_i -h \cdot (s-D_{[i-\tau,i]})^+  -h\cdot\Big( L^s_{[i-\tau, i-1]}+O_{i-\tau}^s\Big)  -r\cdot L^s_i\nonumber\\
 & \geq r\cdot(B^s_{i-1}-B^s_i)+r\cdot D_i -h \cdot (s-D_{[i-\tau,i]})^+-h\cdot \overline{B}^s_{[i-\tau, i-1]}-r\cdot \overline{B}^s_i-h\cdot O_{i-\tau}^s.\label{hanqin-27}
 \end{align}
Note that
\begin{align*}
r\cdot D_i=r\cdot \Big((\overline{I}^s_i+\overline{Q}^s_{i-\tau})\wedge (\overline{B}^s_{i-1}+D_i)\Big)-r\cdot(\overline{B}_{i-1}^s-\overline{B}_i^s).
\end{align*}
Hence, from (\ref{hanqin-27}), we finally have that for $i\geq \tau+1$,
\begin{align}
{\cal P}^s_i(I^s_i,B^s_{i-1},Q^s_{i-\tau}) & \ge r\cdot(B^s_{i-1}-B^s_i)-r\cdot(\overline{B}_{i-1}^s-\overline{B}_i^s)+r\cdot \Big((\overline{I}^s_i+\overline{Q}^s_{i-\tau})\wedge (\overline{B}^s_{i-1}+D_i)\Big)\nonumber\\
& -h \cdot  (s-D_{[i-\tau,i]})^+-h\cdot \overline{B}^s_{[i-\tau, i-1]}-r\cdot \overline{B}^s_i-h\cdot O_{i-\tau}^s.\label{hanqin-28}
 \end{align}
By Proposition \ref{Lemma_Connect P-system and B-system_Dynamic},
\begin{align*}
&\lim_{N\to \infty}\frac{1}{N}\sum_{i=1}^N r \cdot(B^s_{i-1}-B^s_i)=\lim_{N\to \infty}\frac{1}{N}r \cdot(B^s_0-B^s_N)\le\lim_{N\to \infty}\frac{1}{N}r \cdot \overline{B}^s_0=0; \\
&\lim_{N\to \infty}\frac{1}{N}\sum_{i=1}^N r \cdot(B^s_{i-1}-B^s_i)=\lim_{N\to \infty}\frac{1}{N}r \cdot(B^s_0-B^s_N)\ge\lim_{N\to \infty}\frac{1}{N}r \cdot (-\overline{B}^s_N)=0.
\end{align*}
This two relations give
\begin{align}
\lim_{N\to \infty}\frac{1}{N}\sum_{i=1}^N r \cdot(B^s_{i-1}-B^s_i)=0. \label{hanqin-29}
\end{align}
Similarly, we also have
 \begin{align}
\lim_{N\to \infty}\frac{1}{N}\sum_{i=1}^N r \cdot(\overline{B}^s_{i-1}-\overline{B}^s_i)=0. \label{hanqin-30-1}
\end{align}
For the {\bf B}-system, under the $s$-policy, we know that $\overline{B}^s_i$ follows the distribution given by $(D_{[1,\tau+1]}-s)^+$.
Hence, the terms through the third to the sixth of (\ref{hanqin-28}) together is the profit at period-$i$ for the {\bf B}-system under the $s$-policy when the unit sale price is $r$, the unit hold cost is $h$, and the unit backorder cost is $r+\tau h$. With this observation in hand, combining (\ref{hanqin-28})-(\ref{hanqin-30-1}) yields that for any initial state ${\bf x}\in {\cal X}^c$,
\begin{align*}
{\cal C}^{s}(r,h; {\bf x})&=\limsup\limits_{N\to\infty}\frac{1}{N}\mathbb{E}\Big(\sum_{i=1}^N \Pi^s_i(I^s_i,B^s_{i-1},Q^s_{i-\tau})\Big)\\
& =\limsup\limits_{N\to\infty}\frac{1}{N}\mathbb{E}\Big(\sum_{i=\tau+1}^N \Pi^s_i(I^s_i,B^s_{i-1},Q^s_{i-\tau})\Big)\\
& = \limsup\limits_{N\to\infty}\frac{1}{N}\mathbb{E}\Big(\sum_{i=\tau+1}^N \mathbb{E}_{(I^s_i,B^s_{i-1},Q^s_{i-\tau})}\Big[{\cal P}^s_i(I^s_i,B^s_{i-1},Q^s_{i-\tau})\Big]\Big)\\
& \ge \overline{\cal C}^s(r, h, r+\tau h; {\bf x})-\limsup\limits_{N\to\infty}\frac{1}{N}h\cdot\mathbb{E}O^s_{[1,N]}.
\end{align*}
This gives the proposition.\hfill$\Box$

\

\noindent
{\it Proof of {\bf Proposition \ref{lower bound of optimal P-system_prop}}}: Following the idea discussed before the proposition, let $(Q^{\pi^*}_1,Q^{\pi^*}_2,\ldots)$ be the order quantities over periods under the optimal policy $\pi^*$ for the {\bf P}-system.  Now we construct the corresponding order policy denoted by $\overline{\pi}^*$ for the {\bf B}-system as follows:
\begin{align*}
\bar Q^{\overline{\pi}^*}_1=Q^{\pi^*}_1, \ \ \bar Q^{\overline{\pi}^*}_i=Q^{\pi^*}_i+L^{\pi^*}_{i-1} \ \ \mbox{for $i\geq 2$}.
\end{align*}
This shows the above order policy for the {\bf B}-system just order the same quantity as the policy $\pi^*$ plus a mark-up order whose quantity is equal to the lost-sales from the previous period in the {\bf P}-system. At each period, we will partition the backlogged demands in the {\bf B}-system into two parts: the first part will be referred to as ``regular backlogged demand", and the rest will be referred to as ``lost backlogged demand". The way we do this partition is that we will select the same quantity of backlogged demand in the {\bf B}-system as the backlogged demand in the {\bf P}-system, and the rest of the backlogged demand in the {\bf B}-system will be classified into ``lost backlogged demand". The lost backlogged demand can be satisfied until the corresponding make-up order arrives, even there could possibly be some available inventory during the lead time. Notice that we don't use FIFO rules to fulfill the demand, which is actually the optimal allocation rule in the {\bf B}-system.

In sum, the policy we construct above for the {\bf B}-system not only determines the order policy $\overline{\pi}^*$ but also decides the rule about how to satisfy the backlogged demands.
With a little bit notation abuse, we still use $\overline{\pi}^*$ (order decision plus the allocation rule) to denote it.
 The policy $\overline{\pi}^*$ for the {\bf B}-system ensures that: (i) The on-hand inventory and the ``regular backlogged demand" will match the on-hand inventory and backlogged demand in the {\bf P}-system under the optimal policy $\pi^*$; (ii) The amount of the demands that arrive in period-$i$ and are classified as ``lost backlogged demand"
equals to the amount of the demands that arrive in period-$i$, are unsatisfied and eventually become lost-sales in the {\bf P}-system under the optimal policy $\pi^*$;
(iii) At each period, the sales amount is larger than the sales amount for the {\bf P}-system, and the increment part is the lost-sales part in the {\bf P}-system; (iv) Each unit of ``lost backlogged demand" contributes backlogged cost of $(\tau+1)\overline{b}$; (v) Each unit of ``regular backlogged demand" incurs the backorder cost $(\tau+1)\overline{b}$ at most;
and (vi) The ratio of ``regular backlogged demand" and ``lost backlogged demand" is $p/(1-p)$. With these properties in hand,
${\cal C}^{\pi^*}(r,h; {\bf x}) \leq \bar{\cal C}^{\overline{\pi}^*}(r,h,\frac{(1-p)r}{2(1+\tau)}; {\bf x})$.
As for the {\bf B}-system, its optimal policy is of base-stock level, then we have the proposition.\hfill$\Box$

\

\noindent
{\it Proof} of {\bf Proposition \ref{ergodic base-stock policy_thm}}:
For the given $(s,q)$-policy, consider the following subset of the state space ${\cal X}^{(s,q)}$:
\begin{align}
 \mathfrak{F}^{(s,q)}:=\left\{(x_1,x_2,\ldots,x_{\tau+2})\in {\cal X}^{(s,q)}: x_2\geq (s+q)\vee \big(\lfloor\frac{2p}{(1-p^\delta)^{1/\delta}}\left((\tau+1)\mathbb{E}D - (s+q)\right)^+\rfloor+1\big)\right\}.\label{hanqin-9}
\end{align}
As $x_2>0$, we know that $x_1=0$ in $\mathfrak{F}^{(s,q)}$.
We first show that, for any $\delta\in(0,1]$ and $\mathbf{X}^{(s,q)}_1={\bf x}_1\in \mathfrak{F}^{(s,q)}$,
   \begin{align}\label{hanqin-10}
   \mathbb{E}\Big[{\cal V} (\mathbf{X}^{(s,q)}_{\tau+2})\Big|\mathbf{X}^{(s,q)}_1={\bf x}_1\Big]-{\cal V}({\bf x}_1)
    \le -\frac{1}{2}\Big(1-\frac{1+(2^\delta-1)p^\delta}{2^\delta} \Big){\cal V} ({\bf x}_1).
    \end{align}
For concreteness, let ${\bf x}_1=(I_1^{(s,q)},B_0^{(s,q)}, Q^{(s,q)}_{1-\tau},\ldots,Q^{(s,q)}_0)$ and $\mathbf{X}^{(s,q)}_{2+\tau}=(I^{(s,q)}_{2+\tau},B^{(s,q)}_{1+\tau},Q^{(s,q)}_{2},\ldots, Q^{(s,q)}_{1+\tau})$. Under the
$(s,q)$-policy, the order quantity at period-1 is $Q^{(s,q)}_1=(s+q\cdot \mathbb{I}_{\{I^{(s,q)}_1=0\}}-I^{(s,q)}_1-Q^{(s,q)}_{[1-\tau,0]})^+=(s+q-Q^{(s,q)}_{[1-\tau,0]})$.
After the order decision is made, the total inventory on-order $Q^{(s,q)}_{[1-\tau,1]}$ will be received at the beginning of period-$(2+\tau)$.
Hence, the unmet demand in period-$(1+\tau)$ has the relation
$U^{(s,q)}_{1+\tau}\le (D_{[1,1+\tau]}+B^{(s,q)}_0-Q^{(s,q)}_{[1-\tau,1]}-I_1^{(s,q)})^+=(D_{[1,1+\tau]}+B^{(s,q)}_0-s-q)^+$.
This implies that $B^{(s,q)}_{1+\tau}\le \textcolor{blue}{\xi}_{1+\tau; [1,U^{(s,q)}_{1+\tau}]}$.
Hence, for each ${\bf x}_1\in \mathfrak{F}^{(s,q)}$,
\begin{align*}
    \mathbb{E}\Big[{\cal V}(\mathbf{X}^{(s,q)}_{2+\tau})\Big |\mathbf{X}^{(s,q)}_1={\bf x}_1\Big]-{\cal V}({\bf x}_1)&
    =\mathbb{E}\Big[p(B^{(s,q)}_{1+\tau})^\delta\Big |\mathbf{X}^{(s,q)}_1={\bf x}_1\Big] -p(x_2)^\delta\\
    &\le \mathbb{E}\Big[p(\xi_{1+\tau; [1,U^{(s,q)}_{1+\tau}]})^\delta\Big|\mathbf{X}^{(s,q)}_1={\bf x}_1\Big] -p(x_2)^\delta\\
    &\overset{(J)}{\le}p\Big[\mathbb{E}\Big(\xi_{1+\tau; [1,U^{(s,q)}_{1+\tau}]} \Big|\mathbf{X}^{(s,q)}_1={\bf x}_1\Big)^\delta -p(x_2)^\delta\\
    &\overset{(W)}{\leq} p\Big(p \mathbb{E}(D_{[1,1+\tau]}+ x_2- s-q)^+\Big)^\delta-p(x_2)^\delta\\
    &\overset{(d)}{=}p\Big(p \mathbb{E}(D_{[1,1+\tau]} + x_2- s-q)\Big)^\delta-p(x_2)^\delta\\
    &\overset{(d')}\le -\Big(1-\frac{1+(2^\delta-1)p^\delta}{2^\delta}\Big)p(x_2)^\delta\\
    & \le -\frac{1}{2}\Big(1-\frac{1+(2^\delta-1)p^\delta}{2^\delta} \Big){\cal V}({\bf x}_1),
\end{align*}
where $\overset{(J)}{\le}$ is from Jensen's inequality,  $\overset{(W)}{\le}$ comes from Wald's equation, $\overset{(d)}{\le}$ and $\overset{(d')}{\le}$ are both from the the definition of subset $\mathfrak{F}^{(s,q)}$ in (\ref{hanqin-9}), and the last inequality follows from the fact that $x_2\geq 1$ as $x_2$ is nonnegative integer and ${\bf x}\in \mathfrak{F}^{(s,q)}$.
 So (\ref{hanqin-10}) is established. With (\ref{hanqin-10}) in hand,  to show the geometric ergodicity, by Theorem 19.1.3 of \cite{meyntweedie2012}, it suffices to show:
\begin{eqnarray}
\left\{
\begin{array}{ll}
& \mbox{The complement of} \ \mathfrak{F}^{(s,q)} \ \mbox{from the state space ${\cal X}^{(s,q)}$ (denoted by $\mathfrak{F}^{(s,q)}_c$) is a petite set}, \\
&\{X_i^{(s,q)}=(I^{(s,q)}_i,B_{i-1}^{(s,q)},Q_{i-\tau}^{(s,q)},\ldots,Q^{(s,q)}_{i-1}): i\geq 1\} \ \mbox{is $\psi$-irreducible and aperiodic}.
\end{array}
\right.\label{hanqin-11}
\end{eqnarray}

According to the definition of the petite set (p117, \cite{meyntweedie2012}), to prove the first claim of (\ref{hanqin-11}), it suffices to show that there exist a period-$i_o$ after period-1 and a non-trivial measure ${\cal M}$ on ${\cal X}^{(s,q)}$ such that for any subset $B\subseteq {\cal X}^{(s,q)}$ and ${\bf x}=(I_1^{(s,q)},B_0^{(s,q)}, Q^{(s,q)}_{1-\tau},\ldots,Q^{(s,q)}_0)\in \mathfrak{F}^{(s,q)}_c$,
\begin{align}
\mathbb{P}\Big[\mathbf{X}^{(s,q)}_{i_o}\in B\Big |\mathbf{X}^{(s,q)}_1={\bf x}\Big] \ge{\cal M}(B).
\label{hanqin-13}
\end{align}

As $p\neq 1$, {\bf P}-system is not {\bf B}-system.
Conditioning on event ${\cal E}_1$ defined by the demand pattern $D_1=D_2=\cdots=D_{2+\tau}=0$ and the patience pattern in which the whole unmet demands $(B^{(s,q)}_0-I^{(s,q)}_1-Q^{(s,q)}_{1-\tau})^+$ in period-1 become impatient, then
\[
X^{(s,q)}_2=\Big(I^{(s,q)}_2,0,Q^{(s,q)}_{2-\tau},\ldots,Q^{(s,q)}_1\Big) \ \mbox{with} \ Q^{(s,q)}_1=\Big(s+q\cdot \mathbb{I}_{\{I^{(s,q)}_1=0\}} -
I^{(s,q)}_1-Q^{(s,q)}_{[1-\tau,0]}\Big)^+.
\]
Further, this demand pattern and the $(s,q)$-policy give
\begin{align*}
Q_2^{(s,q)}=\Big(s+q\cdot \mathbb{I}_{\{I^{(s,q)}_2=0\}} -I^{(s,q)}_2-Q^{(s,q)}_{[2-\tau,1]}\Big)^+, \ \mbox{and} \
Q_3^{(s,q)}=\cdots=Q_{2+\tau}^{(s,q)}=0.
\end{align*}
This implies that
\[
X^{(s,q)}_{3+\tau}=\Big(I^{(s,q)}_2+Q^{(s,q)}_{[2-\tau,2]},0,0,\ldots,0\Big).
\]
 Let $\kappa=I^{(s,q)}_2+Q^{(s,q)}_{[2-\tau,2]}-s$. As $s\leq I^{(s,q)}_2+Q^{(s,q)}_{[2-\tau,2]}\leq s+q$,  $0\leq \kappa \leq q$, and furthermore, conditioning on event ${\cal E}_2$ given by the demand pattern $D_{3+\tau}=\cdots=D_{2+\kappa+\tau}=1$ and $D_{3+\kappa+\tau}=\cdots=D_{2+\tau+q}=0$, then
\begin{align*}
X_{3+\tau+q}^{(s,q)}=(s,0,0,\ldots,0).
\end{align*}
Write $\kappa_1=\mathbb{E}(D_{[1,1+\tau]}-s)^++1$. Note that $\mathbb{P}[{\cal E}_1]\geq (\alpha_0)^{\tau+2}(1-p)^{\kappa_1}$ and $\mathbb{P}[{\cal E}_2]=(\alpha_0)^{q-\kappa}(\alpha_1)^{\kappa}$. Then
\begin{align}
\mathbb{P}\Big[X_{3+\tau+q}^{(s,q)}=(s,0,0,\ldots,0) \Big| X_{1}^{(s,q)}={\bf x}\Big]\geq (\alpha_0)^{\tau+2}(1-p)^{\kappa_1} \min_{0\leq \kappa \leq q}\Big\{(\alpha_0)^{q-\kappa}(\alpha_1)^{\kappa}\Big\}.\label{hanqin-40-1}
\end{align}
For each ${\bf y}\in {\cal X}^{(s,q)}$, let  ${\cal M}_1({\bf y})=\mathbb{P}[X_2^{(s,q)}={\bf y}|X_1^{(s,q)}=(s,0,\ldots,0)]$, and define
\begin{align}
{\cal M}({\bf y})= (\alpha_0)^{\tau+2}(1-p)^{\kappa_1} \min_{0\leq \kappa \leq q}\Big\{(\alpha_0)^{q-\kappa}(\alpha_1)^{\kappa}\Big\} \cdot {\cal M}_1({\bf y}), \
\mbox{and} \ i_o=4+\tau+q.
\label{hanqin-41-1}
\end{align}
It follows from (\ref{hanqin-40-1}) that
\begin{align}
&\mathbb{P}\Big[\mathbf{X}^{(s,q)}_{i_o}={\bf y}\Big |\mathbf{X}^{(s,q)}_1={\bf x}\Big]\nonumber\\
& \ \ \ =\mathbb{P}\Big[\mathbf{X}^{(s,q)}_{3+\tau+q}=(s,0,\ldots,0)\Big |\mathbf{X}^{(s,q)}_1={\bf x}\Big]
\times \mathbb{P}\Big[\mathbf{X}^{(s,q)}_{4+\tau+q}={\bf y}\Big | \mathbf{X}^{(s,q)}_{3+\tau+q}=(s,0,\ldots,0)\Big]\nonumber\\
& \ \ \ =\mathbb{P}\Big[\mathbf{X}^{(s,q)}_{3+\tau+q}=(s,0,\ldots,0)\Big |\mathbf{X}^{(s,q)}_1={\bf x}\Big] \times {\cal M}_1({\bf y})\nonumber\\
& \ \ \ \geq (\alpha_0)^{\tau+2}(1-p)^{\kappa_1} \min_{0\leq \kappa \leq q}\Big\{(\alpha_0)^{q-\kappa}(\alpha_1)^{\kappa}\Big\}\times {\cal M}_1({\bf y}) ={\cal M}({\bf y}).
\label{hanqin-53}
\end{align}
So the construction of ${\cal M}$ given by (\ref{hanqin-41-1}) makes that (\ref{hanqin-13}) holds.
Hence, (\ref{hanqin-13}) holds. \\[-0.15in]

Now we prove the second claim in (\ref{hanqin-11}). Let $\psi$ be the probability measure on ${\cal X}^{(s,q)}$ as follows:
\begin{align*}
\psi({\bf x})=\left\{
\begin{array}{ll}
1, &\mbox{if ${\bf x}=(s,0,0,\ldots,s)$};\\
0, &\mbox{otherwise}.
\end{array}
\right.
\end{align*}
Then from (\ref{hanqin-40-1}), for any $B \subseteq {\cal X}^{(s,q)}$ and ${\bf x}\in {\cal X}^{(s,q)}$,
\begin{align*}
\mathbb{P}\Big[\sum_{i=1}^\infty \Big\{X_i^{(s,q)}\in B \Big\} \Big| X^{(s,q)}_1={\bf x}\Big]>0 \ \ \mbox{if $(s,0,0,\ldots,0)\in B$}.
\end{align*}
This shows  that the $\psi$-irreducibility of $\{X_i^{(s,q)}=(I^{(s,q)}_i,B_{i-1}^{(s,q)},Q_{i-\tau}^{(s,q)},\ldots,Q^{(s,q)}_{i-1}): i\geq 1\}$.
Further, note that if $\mathbb{P}[D_1=0]=\alpha_0>0$, then
\begin{align*}
\mathbb{P}\Big[ X_2^{(s,q)}\in B \Big| X^{(s,q)}_1=(s,0,0,\ldots,0)\Big]\geq \alpha_0 \ \ \mbox{if $(s,0,0,\ldots,0)\in B$}.
\end{align*}
This gives us that the strong aperiodicity of  $\{X_i^{(s,q)}=(I^{(s,q)}_i,B_{i-1}^{(s,q)},Q_{i-\tau}^{(s,q)},\ldots,Q^{(s,q)}_{i-1}): i\geq 1\}$.
Thus, the second claim in (\ref{hanqin-11}) is proved.
$\lim_{i\rightarrow \infty}\mathbb{E}{\cal V}(X_i^{(s,q)})=\mathbb{E}{\cal V}(X_\infty^{(s,q)})$ directly follows from Theorem 15.0.1 in \cite{meyntweedie2012}.
\hfill$\Box$\\[-0.1in]

\noindent
{\it Proof of {\bf Proposition \ref{couple the diff for positive on-hand inv}}} :
For (i), when $0<I^{(s,q)}_i\le s$ for $i= i_0,\ldots,i_0+k$ for some $i_0\ge2$ and $k\ge 0$, we know that there is no backorders during periods from $i_0-1$ to $i_0+k-1$ under the $(s,q)$-policy. This implies that the $(s,q)$-policy and $s$-policy give us the same order quantities during the periods through $i_0+1$ to $i_0+k$. We further know that under the $(s,q)$-policy, the sales in the periods through $i_0$ to $i_0+k-1$ are equal to the demands, and $D_i < s$ for $i= i_0,\ldots,i_0+k-1$, which will then implies that the sales under the $s$-policy will be $D_i$, and $0< I^s_i\le s$ for all $i = i_0+1,\ldots, i_0+k$.
For period-$(i_0+k)$, since we have both $0 < I^{(s,q)}_{i+k}\le s$ and $0 < I^s_{i+k}\le s$, we then have the sales under these two policies should be $s\wedge D_{i_0+k}$.
In sum, at period-$i$ with $i_0+1\leq i\leq i_0+k-1$, both policies sell $D_i$, and at period-$(i_0+k)$, both policies sell $s\wedge D_{i_0+k}$.
Moreover, the sales relation gives us that $I^{(s,q)}_i = I^s_i = s- D_{i-1}$ for $i=i_0+1,\ldots, i_0+k$, and $I^{(s,q)}_{i_0+k+1}=I^s_{i_0+k+1}= (s- D_{i_0+k})^+$,   which implies that $R^s_i = R^{(s,q)}_i$ for $i=i_0+1,\ldots, i_0+k$. Thus, we have ${\cal R}^s_{[i_0,i_0+k]} -{\cal R}^{(s,q)}_{[i_0, i_0+k]}= {\cal R}^s_{i_0} -{\cal R}^{(s,q)}_{i_0} \leq r\cdot B^s_{i_0-1}$. Hence, (i) is proved.

Now consider (ii).
When $I_i^{(s,q)}=0$ for $i= i_0,\ldots, i_0+k$ with $i_0\ge2$ and $k\ge 0$, we know two possible scenarios: perfectly matches demands, and gets stock out.
At period-$i$ with $i_0\leq i\leq i_0+k-1$, the total demand $D_i+B^{(s,q)}_{i-1}\geq s+q$.
Thus the difference of the profits under these two policies at period-$i$ for $i = i_0,\ldots,i_0+k-1$ follows
    \begin{align*}
  {\cal R}^s_{i} -{\cal R}^{(s,q)}_{i}&= r\cdot\Big( B^s_{i-1}+s\wedge D_{i}\Big)-h(s-D_i)^+ -r\cdot (s+q)\\
        &= r\cdot B_{i-1}^s - (h+r)\cdot(s- D_i)^+ - rq\\
        &= r\cdot\Big(B_{i-1}^s-q\Big) - (h+r)\cdot(s- D_i)^+.
    \end{align*}
Hence we have
\begin{align}
{\cal R}^s_{[i_0,i_0+k-1]} -{\cal R}^{(s,q)}_{[i_0,i_0+k-1]}
 =\sum_{i=i_0}^{i_0+k-1}r\cdot\Big(B_{i-1}^s-q\Big)  - \sum_{i=i_0}^{i_0+k-1}(h+r)\cdot(s- D_i)^+.\label{hanqin-38}
\end{align}
For period-$(i_0+k)$,
\begin{align}
 {\cal R}^s_{i_0+k} -{\cal R}^{(s,q)}_{i_0+k}
        &= r\cdot\Big(B_{i_0+k-1}^s+s\wedge D_{i_0+k}\Big)-h(s - D_{i_0+k})^+ \nonumber\\
        &  \ \ -r\cdot\Big((s+q)\wedge(B_{i_0+k-1}^{(s,q)} + D_{i_0+k})\Big)+h\cdot \Big(s+q-B^{(s,q)}_{i_0+k}-D_{i_0+k}\Big)^+\nonumber\\
        &\le r\cdot\Big(B_{i_0+k-1}^s+s\wedge D_{i_0+k}\Big)-h(s - D_{i_0+k})^+\nonumber\\
        &  \ \  -r\cdot\Big((s+q)\wedge D_{i_0+k}\Big)+h\cdot \Big(s+q-D_{i_0+k}\Big)^+ \nonumber\\
        & \le r\cdot B^s_{i_0+k-1}+hq. \label{hanqin-39}
\end{align}
Combining (\ref{hanqin-38})-(\ref{hanqin-39}) yields that
\begin{align*}
{\cal R}^s_{[i_0,i_0+k]} -{\cal R}^{(s,q)}_{[i_0,i_0+k]}\leq r\cdot B^s_{[i_0-1,i_0+k-1]} -(h+r)\cdot I^s_{[i_0+1,i_0+k]}-krq+hq,
\end{align*}
which is what (ii) claims.

Finally we look at (iii). When $s< I_i^{(s,q)}\le s+q$ for $i= i_0,\ldots, i_0+k$ with $i_0\ge2$ and $k\ge 0$, we know that $D_i< q$ and $B_i^{s}\leq (q-s)^+$ for $i=i_0-1, \ldots,i_0+k-1$.
Based on these observations, the profit difference between these two policies at period-$i$ with $i_0\leq i\leq i_0+k$,
\begin{align*}
{\cal R}_i^s - {\cal R}_i^{(s,q)} &= r\cdot\Big(B_{i-1}^s+s\wedge D_{i}\Big)-h(s-D_i)^+ - r\cdot(I^{(s,q)}_i \wedge D_i) + h\cdot(I^{(s,q)}_i - D_i)^+\\
        &\le r\cdot\Big((q-s)^+ +s\wedge D_{i}\Big)-h(s-D_i)^+ - r\cdot(I^{(s,q)}_i \wedge D_i) + h\cdot(I^{(s,q)}_i - D_i)^+\\
        &\le r\cdot\Big((q-s)^+ +s\wedge D_{i}-I^{(s,q)}_i \wedge D_i\Big) + h\cdot q \\
        &\le h\cdot q+r\cdot(q-s)^+.
    \end{align*}
This gives (iii). Therefore, the proposition is proved.\hfill$\Box$

\noindent
{\it Proof of Proposition} \ref{ergodic-s-q-demand}: In order to use Theorems 2.1-2.4 in \cite{meyntweedie1994}, we want to prove that ${\cal V}({\bf x})$ is the Lyapunov function for one step-transition function
with the small set given by
\begin{align}
 \mathfrak{F}^{(s,q)}=\Big\{(x_1,x_2,\ldots,x_{2+\tau})\in {\cal X}^{(s,q)}: x_2\le \kappa \Big\},\label{dem-1}
 \end{align}
where $d = \mathbb{E}[D]$ given in Subsection \ref{model-1}, and $\kappa=\frac{pd}{[(1+p^\delta)/2]^{1/\delta}-p}$.
For all ${\bf x}=(I^{(s,q)}_1,B^{(s,q)}_0, Q_{1-\tau}^{(s,q)},\ldots$, $Q^{(s,q)}_0)\in {\cal X}^{(s,q)}$ but does not belong to $\mathfrak{F}^{(s,q)}$,
recalling that the unmet demand in period-1 is $U^{(s,q)}_1=(D_1+B_0^{(s,q)}-I^{(s,q)}_1-Q^{(s,q)}_{1-\tau})^+$, we have
\begin{align*}
    \mathbb{E}\Big[{\cal V}(X^{(s,q)}_2)\Big|X^{(s,q)}_1={\bf x}\Big]-{\cal V}({\bf x})&
    =p\cdot \mathbb{E}\Big(\xi_{1;[1,U_1^{(s,q)}]}\Big)^\delta- p\cdot (B^{(s,q)}_0)^\delta\\
    &\le p\cdot \Big(\mathbb{E}\xi_{1;[1,U_1^{(s,q)}]}  \Big)^\delta- p\cdot (B^{(s,q)}_0)^\delta \\
    &=p\cdot\Big(p\cdot \mathbb{E}U^{(s,q)}_1\Big)^\delta- p\cdot (B^{(s,q)}_0)^\delta\\
    &\le p^{1+\delta}\cdot\Big(B^{(s,q)}_0 +\mathbb{E}D_1 \Big)^\delta- p\cdot (B^{(s,q)}_0)^\delta \\
    &= p\cdot \Big (p^\delta (B^{(s,q)}_0+\mathbb{E}D_1)^\delta -(B^{(s,q)}_0)^\delta\Big)\\
    & \leq -\frac{1}{2}(1-p^\delta)p\cdot (B^{(s,q)}_0)^\delta\\
    & \le -\frac{1}{4}(1-p^\delta){\cal V}({\bf x}).
\end{align*}

This gives us that for any ${\bf x}\in {\cal X}^{(s,q)}$,
\begin{align}
 \mathbb{E}\Big[{\cal V}(X^{(s,q)}_2)\Big|X^{(s,q)}_1={\bf x}\Big]-{\cal V}({\bf x})\leq -\frac{1}{4}(1-p^\delta){\cal V}({\bf x})+\Big(\max_{{\bf y}\in \mathfrak{F}^{(s,q)}}
 {\cal V}({\bf y}+\mathbb{E}D_1)\Big)\cdot \mathbb{I}_{\mathfrak{F}^{(s,q)}}({\bf x}). \label{dem-2}
 \end{align}

For ${\bf x}=(I_1^{(s,q)},B_0^{(s,q)},Q_{1-\tau}^{(s,q)},\ldots,Q_0^{(s,q)})\in \mathfrak{F}^{(s,q)}$, Conditioning on event ${\cal E}_1$ defined by the demand pattern $D_1=D_2=\cdots=D_{2+\tau}=0$ and the patience pattern in which the whole unmet demands $(B^{(s,q)}_0-I^{(s,q)}_1-Q^{(s,q)}_{1-\tau})^+$ in period-1 become impatient, then
\[
X^{(s,q)}_2=\Big(I^{(s,q)}_2,0,Q^{(s,q)}_{2-\tau},\ldots,Q^{(s,q)}_1\Big) \ \mbox{with} \ Q^{(s,q)}_1=\Big(s+q\cdot \mathbb{I}_{\{I^{(s,q)}_1=0\}} -
I^{(s,q)}_1-Q^{(s,q)}_{[1-\tau,0]}\Big)^+.
\]
Further, this demand pattern and the $(s,q)$-policy give
\begin{align*}
Q_2^{(s,q)}=\Big(s+q\cdot \mathbb{I}_{\{I^{(s,q)}_2=0\}} -I^{(s,q)}_2-Q^{(s,q)}_{[2-\tau,1]}\Big)^+, \ \mbox{and} \
Q_3^{(s,q)}=\cdots=Q_{2+\tau}^{(s,q)}=0.
\end{align*}
This implies that
\[
X^{(s,q)}_{3+\tau}=\Big(I^{(s,q)}_2+Q^{(s,q)}_{[2-\tau,2]},0,0,\ldots,0\Big).
\]
 Let $\kappa_1=I^{(s,q)}_2+Q^{(s,q)}_{[2-\tau,2]}-s$. As $s\leq I^{(s,q)}_2+Q^{(s,q)}_{[2-\tau,2]}\leq s+q$, then $0\leq \kappa_1 \leq q$, and furthermore, conditioning on event ${\cal E}_2$ given by the demand pattern $D_{3+\tau}=\cdots=D_{2+\kappa+\tau}=1$ and $D_{3+\kappa+\tau}=\cdots=D_{2+\tau+q}=0$, then
\begin{align*}
X_{3+\tau+q}^{(s,q)}=(s,0,0,\ldots,0).
\end{align*}
Note that $\mathbb{P}[{\cal E}_1]\geq (\alpha_0)^{\tau+2}(1-p)^{\kappa}$ and $\mathbb{P}[{\cal E}_2]=(\alpha_0)^{q-\kappa_1}(\alpha_1)^{\kappa_1}$. Then
\begin{align}
\mathbb{P}\Big[X_{3+\tau+q}^{(s,q)}=(s,0,0,\ldots,0) \Big| X_{1}^{(s,q)}={\bf x}\Big]&\geq (\alpha_0)^{\tau+2}(1-p)^{\kappa} \min_{0\leq \kappa_1 \leq q}\Big\{(\alpha_0)^{q-\kappa_1}(\alpha_1)^{\kappa_1}\Big\}\nonumber\\
&=(\alpha_0)^{\tau+2}(1-p)^{\kappa} \Big(\alpha_0\wedge \alpha_1\Big)^q.\label{dem-3}
\end{align}
For each ${\bf y}\in {\cal X}^{(s,q)}$, let  ${\cal M}_1({\bf y})=\mathbb{P}[X_2^{(s,q)}={\bf y}|X_1^{(s,q)}=(s,0,\ldots,0)]$, and define
\begin{align}
{\cal M}(y)= (\alpha_0)^{\tau+2}(1-p)^{\kappa} \Big(\alpha_0\wedge\alpha_1\Big)^q \cdot {\cal M}_1({\bf y}), \
\mbox{and} \ i_o=4+\tau+q.
\label{dem-4}
\end{align}
It follows from (\ref{dem-3}) that
\begin{align}
&\mathbb{P}\Big[\mathbf{X}^{(s,q)}_{i_o}=y\Big |\mathbf{X}^{(s,q)}_1={\bf x}\Big]\nonumber\\
& \ \ \ =\mathbb{P}\Big[\mathbf{X}^{(s,q)}_{3+\tau+q}=(s,0,\ldots,0)\Big |\mathbf{X}^{(s,q)}_1={\bf x}\Big]
\times \mathbb{P}\Big[\mathbf{X}^{(s,q)}_{4+\tau+q}={\bf y}\Big | \mathbf{X}^{(s,q)}_{3+\tau+q}=(s,0,\ldots,0)\Big]\nonumber\\
& \ \ \ =\mathbb{P}\Big[\mathbf{X}^{(s,q)}_{3+\tau+q}=(s,0,\ldots,0)\Big |\mathbf{X}^{(s,q)}_1={\bf x}\Big] \times {\cal M}_1({\bf y})\nonumber\\
& \ \ \ \geq (\alpha_0)^{\tau+2}(1-p)^{\kappa} \Big(\alpha_0\wedge\alpha_1\Big)^q\times {\cal M}_1({\bf y}) ={\cal M}({\bf y}).
\label{dem-5}
\end{align}
So the construction of ${\cal M}$ given by (\ref{dem-4}) makes that $\mathfrak{F}^{(s,q)}$ defined by (\ref{dem-1}) is also a petite set.
Let
\begin{align}
\lambda_0=1-\frac{1}{4}(1-p^\delta), \ b_0=\max_{{\bf y}\in \mathfrak{F}^{(s,q)}}
 {\cal V}({\bf y}+d), \ L=3+\tau+q, \ \mbox{and} \ \zeta=(\alpha_0)^{\tau+2}(1-p)^{\kappa} \Big(\alpha_0\wedge\alpha_1\Big)^q.\label{12-14-2}
 \end{align}
 Then from (\ref{dem-2})-(\ref{dem-3}),
 \begin{align}
 &\mathbb{E}\Big[{\cal V}(X^{(s,q)}_2)\Big|X^{(s,q)}_1={\bf x}\Big]\leq \lambda_0 {\cal V}({\bf x}) + b_0 \cdot \mathbb{I}_{\mathfrak{F}^{(s,q)}}({\bf x}) \ \
 \mbox{for ${\bf x}\in {\cal X}^{(s,q)}$}; \label{dem-5}\\
 & \sum_{i=1}^N \mathbb{P}\Big[X^{(s,q)}_i=(s,0,\ldots,0)\Big| X_1^{(s,q)}={\bf x}\Big]\geq \zeta \ \ \mbox{for ${\bf x}\in \mathfrak{F}^{(s,q)}$};\label{dem-6}\\
 & \mathbb{P}[X_2^{(s,q)}=(s,0,\ldots,0)|X_1^{(s,q)}=(s,0,\ldots,0)]=\mathbb{P}[D_1=0]=\alpha_0. \label{dem-7}
 \end{align}
By (\ref{dem-5})-(\ref{dem-7}), it directly follows from Theorems 2.1-2.4 in \cite{meyntweedie1994} that
\begin{align*}
 \sum_{{\bf x}\in {\cal X}^{(s,q)}}\Big|\mathbb{P}\Big[X_i^{(s,q)}={\bf y}\Big|X_1^{(s,q)}={\bf x}\Big]-\mathbb{P}\Big[
 X_\infty^{(s,q)}={\bf y}\Big]\Big|\leq {\cal V}({\bf x})\frac{(1+\beta(s,q))\rho^{i+1}(s,q)}{\rho(s,q)-\theta(s,q)}
 \end{align*}
 with $\beta(s,q)=b_0\cdot \Big(1+\frac{L^2}{\zeta}\Big)^L$,
 \begin{align*}
  \lambda(s,q)&=1-\frac{1-\lambda_0}{\Pi_{k=0}^{L-1}\Big(1+b_0\frac{L^2}{\zeta^2}(L^2+\zeta^2)^k\Big)},\\
\theta(s,q)&=1-\Big(1-\lambda(s,q)\Big)^2\Big[1-\lambda(s,q)+\beta(s,q)+\beta^2(s,q)\\
& \ \ \ \ \ +\frac{32-8\alpha_0^2}{\alpha^3_0}\frac{\beta^3(s,q)}{(1-\lambda(s,q))^2}
\Big(1-\lambda(s,q)+\beta(s,q)\Big)\Big]^{-1},  \ \mbox{and} \ \
\rho(s,q) \in (\theta(s,q), 1).
 \end{align*}
Hence, the proof of the proposition is completed.\hfill$\Box$

\

\noindent
{\it Proof} of {\bf Theorem \ref{regretbound}}: First note that for any $(s,q)$-policy with $(s,q)\in [0,\underline{s}]\times[0,\underline{q}]$,
${\cal C}^{(s,q)}\leq r(\underline{s}+\underline{q})$. Hence, there exists $c$ such that for any $(s_1,q_1), (s_2,q_2)\in [0,\underline{s}]\times[0,\underline{q}]$,
\begin{align}
\Big| {\cal C}^{(s_1,q_1)}-{\cal C}^{(s_2,q_2)}\Big|\leq c \cdot \Big( |s_2-s_1|+|q_2-q_1|\Big). \label{hanqin-90}
\end{align}
Next we estimate the expected number of the valid periods from the whole epochs.
\begin{align}
\mathbb{E}\Big[\sum_{\ell=1}^{L_{\cal E}(N)}\sqrt{2T_\ell(s^*_\ell,q^*_\ell)} \Big]&=\mathbb{E}\Big[\sum_{(s,q)}\sum_{\ell=1}^{L_{\cal E}(N)}\sqrt{2T_\ell(s,q)\mathbb{I}_{\{(s^*_\ell,q^*_\ell)=(s,q)\}}} \ \Big]\nonumber\\
& = \mathbb{E}\Big[\sum_{(s,q)}\sum_{\ell=1}^{L^{(s,q)}_{\cal E}(N)}\sqrt{2^k} \Big] \leq (\sqrt 2+1)\underline{s}\underline{q}\sqrt{N},\label{hanqin-100}
\end{align}
where in the last inequality, we use  $L^{(s,q)}_{\cal E}(N)\leq \log_2N$.
According to the selection $(s^*_\ell,q^*_\ell)$ in epoch-$\ell$, we have
\begin{align}
&\mathbb{E} \Big[\sum_{\ell=1}^{L_{\cal E}(N)}\sum_{k=\nu_\ell}^{\eta_\ell} \Big({\cal C}^{(s^*,q^*)}-\breve{\cal R}_k^{(s^*_\ell,q^*_\ell)}\Big)\mathbb{I}_{\{C_\ell\}} \Big]\nonumber\\
& \ \ \ \leq \mathbb{E} \Big[\sum_{\ell=1}^{L_{\cal E}(N)}\sum_{k=\nu_\ell}^{\eta_\ell} \Big({\cal C}^{(s^*,q^*)}-F_\ell^{(s^*,q^*)}\Big)\mathbb{I}_{\{C_\ell\}} \Big]+
\mathbb{E} \Big[\sum_{\ell=1}^{L_{\cal E}(N)}\sum_{k=\nu_\ell}^{\eta_\ell} \Big( F_\ell^{(s^*_\ell,q^*_\ell)}-\breve{\cal R}_k^{(s^*_\ell,q^*_\ell)}\Big)\mathbb{I}_{\{C_\ell\}} \Big].\label{hanqin-91}
\end{align}
For the notation simplicity, the first and second terms on the right-hand side in (\ref{hanqin-91}) are named as ${\cal T}_1$ and ${\cal T}_2$, respectively.
In view of the way to produce $F_\ell^{(s,q)}$ in the algorithm and (\ref{hanqin-100}), we have
\begin{align}
{\cal T}_1&=\mathbb{E}\Big[\sum_{\ell=1}^{L_{\cal E}(N)}\sum_{k=\nu_\ell}^{\eta_\ell}
\Big({\cal C}^{(s^*,q^*)}-\frac{1}{|{\cal N}_\ell(s^*,q^*)|}\times\nonumber\\
& \ \ \ \ \ \ \ \ \ \ \ \ \ \times \sum_{(s,q)\in {\cal N}_\ell(s^*,q^*)}\Big\{ \widetilde{\cal R}^{(s,q)}_{\ell-1}+\frac{H_{\ell-1}}{T_{\ell-1}(s,q)}\Big(\sqrt{2(T_{\ell-1}(s,q)-1)\ln\frac{2}{\delta_\ell}}+1\Big)\Big\}\Big)\mathbb{I}_{\{C_\ell\}}\Big]\nonumber\\
& \le \mathbb{E}\Big[\sum_{\ell=1}^{L_{\cal E}(N)}\sum_{k=\nu_\ell}^{\eta_\ell}
\Big({\cal C}^{(s^*,q^*)}-\frac{1}{|{\cal N}_\ell(s^*,q^*)|}\sum_{(s,q)\in {\cal N}_\ell(s^*,q^*)}{\cal C}^{(s,q)}\Big)\mathbb{I}_{\{C_\ell\}}\Big]\nonumber\\
& \leq  \mathbb{E}\Big[ \sum_{\ell=1}^{L_{\cal E}(N)} T_\ell(s^*_\ell,q^*_\ell) \cdot \frac{(\varepsilon_1+\varepsilon_2)c}{1\vee\sqrt{\eta_{\ell-1}}} \mathbb{I}_{\{C_\ell\}}\Big]
\ \ \mbox{(by (\ref{hanqin-90}))} \nonumber\\
& \le(\varepsilon_1+\varepsilon_2)c\cdot \mathbb{E}\Big[\sum_{\ell=1}^{L_{\cal E}(N)}\sqrt{2T_\ell(s^*_\ell,q^*_\ell)} \Big] \ \ \mbox{(by $T_{\ell}(s,q)\leq 2\eta_{\ell-1}$)}\nonumber\\
& \leq (\sqrt 2+1)(\varepsilon_1+\varepsilon_2)c\underline{s}\underline{q}\sqrt{N},\label{hanqin-92}
\end{align}

We now turn to the term ${\cal T}_2$. Again, conditioning on $C_\ell$ holds and by the definition of $F_\ell^{(s,q)}$ in the algorithm, we have
\begin{align}
{\cal T}_2 &\leq \mathbb{E} \Big[\sum_{\ell=1}^{L_{\cal E}(N)}\sum_{k=\nu_\ell}^{\eta_\ell} \Big(\frac{1}{|{\cal N}_\ell(s^*_\ell,q^*_\ell)|}\sum_{(s,q)\in {\cal N}_\ell(s^*_\ell,q^*_\ell)}\Big\{ \widetilde{\cal R}^{(s,q)}_{\ell-1}+\frac{H_{\ell-1}}{T_{\ell-1}(s,q)}\Big(\sqrt{2(T_{\ell-1}(s,q)-1)\ln\frac{2}{\delta_\ell}}+1\Big)\Big\}\nonumber\\
& \ \ \ \ \ \ \ \ \ \ \ \ \ \ \ \ \ \ \ \ \ \ - \breve{\cal R}_k^{(s^*_\ell,q^*_\ell)}\Big)\mathbb{I}_{\{C_\ell\}}\Big]. \label{hanqin-93}
\end{align}
According to the definition of set ${\cal N}_\ell(s,q)$, we have that for $(s,q)\in {\cal N}_\ell(s^*_\ell,q^*_\ell)$,
\begin{align*}
\frac{H_{\ell-1}}{T_{\ell-1}(s,q)}\Big(\sqrt{2(T_{\ell-1}(s,q)-1)\ln\frac{2}{\delta_\ell}}+1\Big)\leq \frac{H_{\ell-1}}{T_{\ell-1}(s^*_\ell,q^*_\ell)}\Big(\sqrt{2(T_{\ell-1}(s^*_\ell,q^*_\ell)-1)\ln\frac{2}{\delta_\ell}}+1\Big).
\end{align*}
Hence, from (\ref{hanqin-90}) and (\ref{hanqin-93}),
\begin{align}
{\cal T}_2 &\leq \mathbb{E} \Big[\sum_{\ell=1}^{L_{\cal E}(N)}\sum_{k=\nu_\ell}^{\eta_\ell} \Big(\frac{1}{|{\cal N}_\ell(s^*_\ell,q^*_\ell)|}\sum_{(s,q)\in {\cal N}_\ell(s^*_\ell,q^*_\ell)}\Big\{{\cal C}^{(s,q)}+\frac{2H_{\ell-1}}{T_{\ell-1}(s,q)}\Big(\sqrt{2(T_{\ell-1}(s,q)-1)\ln\frac{2}{\delta_\ell}}+1\Big)\Big\}\nonumber\\
& \ \ \ \ \ \ \ \ \ \ \ \ \ \ \ \ \ \ \ \ \ \ - \breve{\cal R}_k^{(s^*_\ell,q^*_\ell)}\Big)\mathbb{I}_{\{C_\ell\}}\Big]\nonumber\\
&\leq \mathbb{E} \Big[\sum_{\ell=1}^{L_{\cal E}(N)}\sum_{k=\nu_\ell}^{\eta_\ell} \Big(\frac{1}{|{\cal N}_\ell(s^*_\ell,q^*_\ell)|}\sum_{(s,q)\in {\cal N}_\ell(s^*_\ell,q^*_\ell)}\Big\{  {\cal C}^{(s,q)}+\frac{2H_{\ell-1}}{T_{\ell-1}(s^*_\ell,q^*_\ell)}\Big(\sqrt{2(T_{\ell-1}(s^*_\ell,q^*_\ell)-1)\ln\frac{2}{\delta_\ell}}+1\Big)\Big\}\nonumber\\
& \ \ \ \ \ \ \ \ \ \ \ \ \ \ \ \ \ \ \ \ \ \ -\breve{\cal R}_k^{(s^*_\ell,q^*_\ell)}\Big)\mathbb{I}_{\{C_\ell\}}\Big]\nonumber\\
& \leq \mathbb{E} \Big[\sum_{\ell=1}^{L_{\cal E}(N)}\sum_{k=\nu_\ell}^{\eta_\ell}\Big({\cal C}^{( s^*_\ell,q^*_\ell)}-\breve{\cal R}_k^{(s^*_\ell,q^*_\ell)}\Big)\Big]+\mathbb{E}\Big[ \sum_{\ell=1}^{L_{\cal E}(N)} T_\ell(s^*_\ell,q^*_\ell) \cdot \frac{(\varepsilon_1+\varepsilon_2)c}{1\vee\sqrt{\eta_{\ell-1}}} \mathbb{I}_{\{C_\ell\}}\Big]\nonumber\\
& \ \ +\mathbb{E} \Big[\sum_{\ell=1}^{L_{\cal E}(N)}\sum_{k=\nu_\ell}^{\eta_\ell}\frac{2H_{\ell-1}}{T_{\ell-1}(s^*_\ell,q^*_\ell)}\Big(\sqrt{2(T_{\ell-1}(s^*_\ell,q^*_\ell)-1)\ln\frac{2}{\delta_\ell}}+1\Big)\Big].\label{hanqin-94}
\end{align}
Note that
\begin{align}
&\mathbb{E} \Big[\sum_{\ell=1}^{L_{\cal E}(N)}\sum_{k=\nu_\ell}^{\eta_\ell}\frac{2H_{\ell-1}}{T_{\ell-1}(s^*_\ell,q^*_\ell)}\Big(\sqrt{2(T_{\ell-1}(s^*_\ell,q^*_\ell)-1)\ln\frac{2}{\delta_\ell}}+1\Big)\Big]\nonumber\\
& \ \ \ \le \mathbb{E} \Big[\sum_{\ell=1}^{L_{\cal E}(N)} 4H_{\ell-1}\Big(\sqrt{2(T_{\ell-1}(s^*_\ell,q^*_\ell)-1)\ln\frac{2}{\delta_\ell}}+1\Big)\Big]\nonumber\\
& \ \ \ \le 4H_{L_{\cal E}(N)}\Big( \sqrt{\ln\frac{2}{\delta_{L_{\cal E}(N)}}} \mathbb{E} \Big[\sum_{\ell=1}^{L_{\cal E}(N)}\sqrt{2T_{\ell}(s^*_\ell,q^*_\ell)}\Big]+L_{\cal E}(N)\Big).\label{hanqin-95}
\end{align}
Now consider $\mathbb{E} \Big[\sum_{\ell=1}^{L_{\cal E}(N)}\sum_{k=\nu_\ell}^{\eta_\ell}\Big({\cal C}^{( s^*_\ell,q^*_\ell)}- \breve{\cal R}_k^{(s^*_\ell,q^*_\ell)}\Big)$. Define
\begin{align*}
C_\ell^{(s^*_\ell,q^*_\ell)}=\Big\{
 \Big|\sum_{k=\nu_\ell}^{\eta_\ell}
\Big(\breve{\cal R}_k^{(s^*_\ell,q^*_\ell)}-{\cal C}^{(s^*_\ell,q^*_\ell)}\Big)\Big|
\leq H_\ell\Big(\sqrt{2(T_\ell(s^*_\ell,q^*_\ell)-1)\ln\frac{2}{\delta_\ell}}+1\Big)\Big\},
\end{align*}
then
\begin{align}
&\mathbb{E} \Big[\sum_{\ell=1}^{L_{\cal E}(N)}\sum_{k=\nu_\ell}^{\eta_\ell}\Big({\cal C}^{( s^*_\ell,q^*_\ell)}- \breve{\cal R}_k^{(s^*_\ell,q^*_\ell)}\Big)\mathbb{I}_{\{C^{(s^*_\ell,q^*_\ell)}_{\ell}\}}\Big]\nonumber\\
& \ \ \ \leq
H_N \mathbb{E} \Big[\sum_{\ell=1}^{L_{\cal E}(N)}\Big(\sqrt{2(T_\ell(s^*_\ell,q^*_\ell)-1)\ln\frac{2}{\delta_\ell}}+1\Big)\Big]\nonumber\\
& \ \ \ \le H_{L_{\cal E}(N)}\cdot L_{\cal E}(N)+H_N \sqrt{\ln\frac{2}{\delta_{L_{\cal E}(N)}}} \times \mathbb{E} \Big[\sum_{\ell=1}^{L_{\cal E}(N)}\sqrt{2T_\ell(s^*_\ell,q^*_\ell)}\Big]. \label{hanqin-96}
\end{align}
On the other hand, using Theorem \ref{dem-bound}, similar to (\ref{hanqin-81}),
\begin{align}
&\mathbb{E} \Big[\sum_{\ell=1}^{L_{\cal E}(N)}\sum_{k=\nu_\ell}^{\eta_\ell}\Big({\cal C}^{( s^*_\ell,q^*_\ell)}- \breve{\cal R}_k^{(s^*_\ell,q^*_\ell)}\Big)\mathbb{I}_{\{C^{(s^*_\ell,q^*_\ell)}_{\ell,c}\}}\Big]\leq r\underline{s}\underline{q}(\underline{s}+\underline{q})\log_2N.\label{hanqin-96}
\end{align}
Combining (\ref{hanqin-91}), (\ref{hanqin-92}), and (\ref{hanqin-94})-(\ref{hanqin-96}) yields that
\begin{align}
&\mathbb{E} \Big[\sum_{\ell=1}^{L_{\cal E}(N)}\sum_{k=\nu_\ell}^{\eta_\ell} \Big({\cal C}^{(s^*,q^*)}-\breve{\cal R}_k^{(s^*_\ell,q^*_\ell)}\Big)\mathbb{I}_{\{C_\ell\}} \Big]\nonumber \\
& \ \ \ \leq  5 H_N \Big(\sqrt{\ln\frac{2}{\delta_N}} \mathbb{E} \Big[\sum_{\ell=1}^{L_{\cal E}(N)}\sqrt{2T_{\ell}(s^*_\ell,q^*_\ell)}\Big]+L_{\cal E}(N)\Big)\nonumber\\
& \ \ \ \ +2(\sqrt 2+1)(\varepsilon_1+\varepsilon_2)c\underline{s}\underline{q}\sqrt{N}+r\underline{s}\underline{q}(\underline{s}+\underline{q})\log_2N.\label{hanqin-101}
\end{align}
The theorem directly follows from (\ref{hanqin-80})-(\ref{hanqin-81}), (\ref{hanqin-91})-(\ref{hanqin-92}), and (\ref{hanqin-101}).
\hfill$\Box$

\end{APPENDICES}

\end{document}